\documentclass[10pt]{article}

\usepackage[a4paper, total={6.5in, 10in}, includehead]{geometry}
\usepackage[parfill]{parskip}
\usepackage{needspace}
\usepackage{microtype}

\usepackage[runin]{abstract}
\usepackage[twoside]{fancyhdr}
\usepackage[explicit]{titlesec}
\usepackage{orcidlink}

\usepackage[T1]{fontenc}
\usepackage[cal=esstix,
            scr=boondoxupr,
            frak=euler,
            bb=libus,
            scaled=0.95]{mathalpha}
\usepackage[OT1, small,
            euler-hat-accent,
            euler-digits]{eulervm}
\usepackage{Baskervaldx}
\usepackage{pifont}

\usepackage[leqno]{amsmath}
\usepackage{amssymb}
\usepackage{amsthm}
\usepackage{mathtools,thmtools}
\usepackage{stmaryrd}
\usepackage{xfrac,xparse,xspace}
\usepackage{stackrel,scalerel}

\usepackage[shortlabels]{enumitem}
\usepackage{booktabs,diagbox,makecell,multirow}
\usepackage{arydshln}

\usepackage[english]{babel}
\usepackage{csquotes}
\usepackage[backend=biber,
            style=alphabetic,
            hyperref=auto,
            useprefix=false,
            maxbibnames=99,
            maxalphanames=99,
            maxcitenames=99]{biblatex}
\usepackage{hyperref}
\usepackage{cleveref}
\usepackage{url}
\usepackage{crossreftools}

\usepackage[dvipsnames]{xcolor}
\usepackage{graphicx}
\usepackage{tikz,tikz-cd,tikz-imagelabels}
\usepackage{dynkin-diagrams}
\usepackage{rotating}
\usepackage{float}
\usepackage[font=small,
            labelfont=bf,
            labelsep=period,
            format=plain,
            singlelinecheck=false]{caption}

\hypersetup{
    hypertexnames = false,
    bookmarksdepth = 2,
    bookmarksopen = true,
    colorlinks,
    linkcolor = ., citecolor = teal, urlcolor = .,
    pdfstartview={XYZ null null 1}}

\DeclareSortingNamekeyTemplate{
  \keypart{\namepart{family}}
  \keypart{\namepart{prefix}}
  \keypart{\namepart{given}}
  \keypart{\namepart{suffix}}
}
\renewbibmacro*{begentry}{\midsentence}
\makeatletter
\AtBeginDocument{\toggletrue{blx@useprefix}}
\makeatother

\declaretheoremstyle[
    spaceabove=\parsep, spacebelow=\parsep,
    headfont=\bfseries, notefont=\normalfont, bodyfont=\itshape,
    headpunct={.}, notebraces={}{}, postheadspace={ },
    ]{basic-theorem}
\declaretheoremstyle[
    spaceabove=\parsep, spacebelow=\parsep,
    headfont=\bfseries, notefont=\normalfont, bodyfont=\normalfont,
    headpunct={.}, notebraces={(}{)}, postheadspace={ },
    ]{basic-definition}
\declaretheoremstyle[
    spaceabove=\parsep, spacebelow=\parsep,
    headfont=\itshape, notefont=\normalfont, bodyfont=\normalfont,
    headpunct={.}, notebraces={(}{)}, postheadspace={ },
    ]{basic-remark}

\theoremstyle{basic-theorem}
\newtheorem{keytheorem}{Theorem}[section]

\newtheorem{theorem}{Theorem}[section]
\newtheorem{corollary}[theorem]{Corollary}
\newtheorem{lemma}[theorem]{Lemma}
\newtheorem{proposition}[theorem]{Proposition}

\theoremstyle{basic-definition}
\newtheorem{definition}[theorem]{Definition}
\newtheorem{example}[theorem]{Example}
\newtheorem{algorithm}[theorem]{Algorithm}

\theoremstyle{basic-remark}
\newtheorem{remark}[theorem]{Remark}

\titleformat{\section}{\normalfont\bfseries\large\centering}{\S\,\thesection}{0.5em}{#1}{}
\crefname{section}{\S\kern -1pt}{\S\S \kern -1pt}
\newcommand{\fakesection}[1]{%
    \par\refstepcounter{section}%
    \sectionmark{#1}%
    \addcontentsline{toc}{section}{\protect\numberline{\thesection}#1}%
}

\titleformat{\subsection}[runin]{\normalfont\bfseries}{\S\,\thesubsection}{0.5em}{#1.}{}
\titlespacing{\subsection}{0pt}{1em}{\wordsep}[]
\crefname{subsection}{\S \kern -1pt}{\S\S \kern -1pt}
\newcommand{\fakesubsection}[1]{%
  \par\refstepcounter{subsection}%
  \subsectionmark{#1}%
  \addcontentsline{toc}{subsection}{\protect\numberline{\thesubsection}#1}%
}

\titleformat{\paragraph}[runin]{\normalfont\normalsize\bfseries}{\theparagraph}{1em}{#1.}{}

\imagelabelset{
    coarse grid color = red,
    fine grid color = gray,
    image label font = \sffamily\bfseries\scriptsize,
    image label distance = 2mm,
    image label back = none,
    image label text = black,
    coordinate label font = \normalfont\scriptsize,
    coordinate label distance = 2mm,
    coordinate label back = none,
    coordinate label text = black,
    annotation font = \normalfont\scriptsize,
    arrow distance = 0,
    border thickness = 0.2pt,
    arrow thickness = 0.1pt,
    tip size = 0.5mm,
    outer dist = 0.5cm,
}

\makeatletter
\newsavebox{\@brx}
\newcommand{\llangle}[1][]{\savebox{\@brx}{\(\m@th{#1\langle}\)}%
  \mathopen{\copy\@brx\kern-0.5\wd\@brx\usebox{\@brx}}}
\newcommand{\rrangle}[1][]{\savebox{\@brx}{\(\m@th{#1\rangle}\)}%
  \mathclose{\copy\@brx\kern-0.5\wd\@brx\usebox{\@brx}}}
\makeatother

\newcommand{\cmark}{\ding{51}}%
\newcommand{\xmark}{\ding{55}}%

\newcommand{\ul}{\underline}
\newcommand{\ol}{\overline}
\newcommand{\wt}{\widetilde}
\newcommand{\dual}{\check{}\,}
\newcommand{\smashdual}{\kern -1.5pt \check{}\;}
\newcommand{\orth}{{}^\perp}

\renewcommand{\inf}{\mathop\textnormal{inf}}
\renewcommand{\sup}{\mathop\textnormal{sup}}

\newcommand{\subsetdot}{\mathrel{\ooalign{$\subset$\cr \hidewidth\hbox{$\cdot\mkern5mu$}\cr}}}
\newcommand{\supsetdot}{\mathrel{\ooalign{$\supset$\cr \hidewidth\hbox{$\cdot\mkern5mu$}\cr}}}

\newcommand{\pair}[2]{\text{\scalebox{0.9}{\(\left(\tfrac{#1}{#2}\right)\)}}}
\newcommand{\XZ}{\pair{X}{Z}}
\newcommand{\WZ}{\pair{W}{Z}}
\newcommand{\XY}{\pair{X}{Y}}
\newcommand{\WY}{\pair{W}{Y}}

\newcommand{\into}{\hookrightarrow}
\newcommand{\onto}{\twoheadrightarrow}
\newcommand{\isoto}{\mathrel{\ooalign{
     $\to$\cr
     \hidewidth\raise.3em\hbox{$\scaleobj{.7}{\sim}\mkern7mu$}\cr
    }}
}
\newcommand{\leftrarrows}{\mathrel{\raise.75ex\hbox{\oalign{%
  $\scriptstyle\leftarrow$\cr
  \vrule width0pt height.5ex$\hfil\scriptstyle\relbar$\cr}}}}
\newcommand{\lrightarrows}{\mathrel{\raise.75ex\hbox{\oalign{%
  $\scriptstyle\relbar$\hfil\cr
  $\scriptstyle\vrule width0pt height.5ex\smash\rightarrow$\cr}}}}
\newcommand{\Rrelbar}{\mathrel{\raise.75ex\hbox{\oalign{%
  $\scriptstyle\relbar$\cr
  \vrule width0pt height.5ex$\scriptstyle\relbar$}}}}

\makeatletter
\def\leftrightarrowsfill@{\arrowfill@\leftrarrows\Rrelbar\lrightarrows}
\newcommand{\xleftrightarrows}[2][]{\ext@arrow 3399\leftrightarrowsfill@{#1}{#2}}
\makeatother

\DeclareMathOperator{\img}{\textnormal{im}}

\DeclareMathOperator{\Hom}{\textnormal{Hom}}
\DeclareMathOperator{\End}{\textnormal{End}}
\DeclareMathOperator{\Ext}{\textnormal{Ext}}
\DeclareMathOperator{\Tor}{\textnormal{Tor}}

\newcommand{\rmod}{\ensuremath{\mathcal{m\kern -0.5pt o \kern -0.5pt d}}}
\newcommand{\rmodif}{\ensuremath{\mathcal{m\kern -0.5pt o \kern -0.5pt d \kern
-0.5pt i \kern -0.8pt f}}}
\newcommand{\rproj}{\ensuremath{\mathcal{p\kern -0.5pt r \kern -0.5pt o \kern
-1pt j}}}
\newcommand{\rflmod}{\ensuremath{\mathcal{f\kern -0.5pt l\kern -1pt m\kern
-0.5pt o\kern -0.5pt d}}}

\DeclareMathOperator{\Db}{\mathbf{D}^\textnormal{b}\kern -1pt}
\DeclareMathOperator{\Dm}{\mathbf{D}^0\kern -1.5pt}
\DeclareMathOperator{\Dfl}{\mathbf{D}^\textnormal{fl}\kern -1pt}

\DeclareMathOperator{\RHom}{\ensuremath{\mathbf{R}\textnormal{Hom}}}
\newcommand{\LTensor}{\ensuremath{\otimes^\mathbf{L}}}

\DeclareMathOperator{\Spec}{\ensuremath{\textnormal{Spec}}}
\DeclareMathOperator{\Coh}{\mathcal{C\kern -2pt o\kern -1pt h\kern -1pt}}
\DeclareMathOperator{\coh}{\mathcal{c\kern -1pt o\kern -1pt h\kern -1pt}}
\newcommand{\cohsupp}[1]{\mathcal{c\kern -1pt o\kern -1pt h\kern -1pt}_{#1}}
\DeclareMathOperator{\anticoh}{\overline{\mathcal{c\kern -1pt o\kern
    -1pt}}\mathcal{h\kern -1pt}}
\DeclareMathOperator{\sheafEnd}{\mathcal{E\kern -2pt n\kern -1pt d\kern -1pt}}

\DeclareMathOperator{\zeroPer}{\mathcal{P\kern -0.5pt e\kern -0.5pt r\kern
-2pt}}
\DeclareMathOperator{\zeroper}{\mathcal{p\kern -0.5pt e\kern -0.5pt r\kern -2pt}}
\newcommand{\zeropersupp}[1]{\mathcal{p\kern -0.5pt e\kern -0.5pt r\kern -1pt}_{#1}\kern -1pt}
\DeclareMathOperator{\antizeroper}{\overline{\mathcal{\kern -0.5pt p\kern -0.5pt
e\kern -0.5pt r\kern -2pt}}}

\DeclareMathOperator{\KK}{\mathbf{K}}
\DeclareMathOperator{\CC}{\mathbf{C}}
\newcommand{\WW}{\textnormal{W}}
\newcommand{\CycleGroup}{\textnormal{Z}_1\kern -1pt}
\newcommand{\JJ}{\mathfrak{J}}
\newcommand{\II}{\mathfrak{I}}
\newcommand{\HH}{\ensuremath{\textnormal{H}}}
\newcommand{\OO}{\mathscr{O}}
\newcommand{\LL}{\mathscr{L}}
\newcommand{\VV}{\mathscr{V\kern -1 pt}}
\newcommand{\VdB}{\textnormal{VdB}}
\newcommand{\flop}{\textnormal{flop}}
\newcommand{\ulLambda}{\rotatebox[origin=c]{180}{\(\mathscr{V}\)}}

\newcommand{\Chamb}{\mathrm{C}}
\newcommand{\bbZ}{\mathbb{Z}}
\newcommand{\bbF}{\mathbb{F}}
\newcommand{\bbR}{\mathbb{R}}
\DeclareMathOperator{\Chambers}{\textnormal{Cham}}
\DeclareMathOperator{\cl}{\textnormal{cl}}

\newcommand{\curlyHom}{\ensuremath{\mathcal{H\kern -2pt o\kern -1pt m}}}
\newcommand{\curlyExt}{\ensuremath{\mathcal{E\kern -0.5pt x\kern -1pt t}}}

\newcommand{\MaxMod}{\mathop\textnormal{MM}{}^\ensuremath{N}}
\newcommand{\MaxModGen}{\mathop\textnormal{MMG}{}^\ensuremath{N}}
\newcommand{\Bir}{\mathop\textnormal{Bir}\kern -1pt}
\DeclareMathOperator{\Pic}{\textnormal{Pic}}
\DeclareMathOperator{\RealPic}{\textnormal{Pic}_\mathbb{R}}

\DeclareMathOperator{\tors}{\textnormal{tors}}
\DeclareMathOperator{\torf}{\textnormal{torf}}

\DeclareMathOperator{\tilt}{\textnormal{tilt}}
\DeclareMathOperator{\simp}{\textnormal{simp}}

\DeclareMathOperator{\sbrick}{\textnormal{sbrick}}
\DeclareMathOperator{\HFan}{\textnormal{HFan}}
\DeclareMathOperator{\faces}{\textnormal{faces}}

\DeclareMathOperator{\Supp}{\textnormal{Supp}}

\DeclareMathOperator{\ExQuiv}{\textnormal{ExQuiv}}

\DeclareMathOperator{\Arr}{\textnormal{Arr}}
\newcommand{\Root}{\textnormal{Root}}

\newcommand{\tr}{\textnormal{tr}}
\newcommand{\tf}{\textnormal{tf}}
\renewcommand{\tt}{\textnormal{tt}}
\renewcommand{\ss}{\textnormal{ss}}

\newcommand{\affE}[2]{\dynkin[extended, edge length = 5pt,
edge/.style={draw=none}, root radius=1.5pt, affine mark ={#1}, ordering=Dynkin] E{#2}}

\begin{document}

\title{\huge Torsion pairs and 3-fold flops}
\author{Parth Shimpi \orcidlink{0009-0008-2066-7474}}

\AtEndDocument{\bigskip{\footnotesize%
        \textsc{The Mathematics and Statistics Building, University
        of Glasgow, University Place, Glasgow G12 8QQ, UK.}\\
        \textit{Email address: }\texttt{parth.shimpi@glasgow.ac.uk}\\
        \textit{Web: }\texttt{https://pas201.user.srcf.net}}}

\date{}
\maketitle

\begin{abstract}
    This paper classifies t-structures on the local derived category of
    a 3-fold flopping contraction, that are intermediate with respect to the
    heart of perverse coherent sheaves. Equivalently,
    this describes the complete lattice of torsion classes for the associated
    modification algebra. The intermediate hearts are (1) categories of coherent
    sheaves on birational models and tilts thereof in skyscrapers,  (2) algebraic
    t-structures described in the homological minimal model programme, or (3)
    combinations of the above over appropriate open covers. An analogous
    classification is also proved for minimal (and partial)
    resolutions of Kleinian singularities, thus providing a description of all
    torsion pairs in the module categories of (contracted) affine preprojective
    algebras. The results have immediate applications to the classification of
    spherical modules and (semi)bricks, and are first steps towards describing
    all t-structures and spherical objects in derived categories of surfaces and
    3-folds.
\end{abstract}

\medskip
\fakesection{Introduction}
It is a truth universally acknowledged, that a mathematician in possession of a
triangulated category, must be in want of tools to understand its
autoequivalences. This truth guides the homological algebraist's pursuit of
t-structures, the representation theorist's pursuit of torsion pairs, and the
geometer's of spherical objects.

Such is also the predicament we find ourselves in;
this paper studies the derived category of a Gorenstein terminal
3-fold \(X\) appearing in a flopping contraction \(\pi\colon X\to Z=\Spec
(R,\mathfrak{m})\) over a complete local base. The assumptions on the
singularities of \(X\) are equivalent to requiring that \(R\) is an isolated
compound Du Val (cDV) singularity. The bounded derived category of coherent
sheaves \(\Db X\), as well as the full subcategory supported on the exceptional fiber
\[
    \Dm X= \left\{x\in \Db X
        \;\middle\vert\;
        \mathop{\text{Supp}}x\subseteq \pi^{-1}[\mathfrak{m}]\right\}
\]
has been of interest to birational and symplectic geometers alike~\cite{aspinwallPointsPointView2003,todaStabilityConditionsCrepant2008,hiranoFaithfulActionsHyperplane2018,hiranoStabilityConditions3fold2023,keatingSymplectomorphismsSphericalObjects2024}.
The same can be said of the analogous situation in dimension \(2\), where \(Z\)
is a canonical surface singularity and \(X\) a (partial) resolution
\cite{crawley-boeveyExceptionalFibresKleinian2000,ishiiAutoequivalencesDerivedCategories2005,bridgelandStabilityConditionsKleinian2009,ishiiStabilityConditionsAnSingularities2010,bapatSphericalObjectsStability2023}.
In either setting, an understanding of autoequivalences, t-structures, and
spherical objects is desirable.

When \(X\) is smooth, the autoequivalences of the derived category are sometimes
controlled by \textit{spherical objects}. In general they are not,
so Hara--Wemyss \cite{haraSphericalObjectsDimensions2024} argue for the study of
broader collections of objects which behave like the simples of an Abelian
category, namely \textit{(semi)bricks}. The rich interplay between semibricks
and t-structures has been a staple tool for representation theorists
\cite{marksTorsionClassesWide2017, asaiSemibricks2020}, and Hara--Wemyss use
this to study the null subcategory
\(\mathcal{C}=\ker(\mathbf{R}\pi_\ast)\subseteq \Db X\) which is known to be the
`finite-type' counterpart to \(\Dm X\), the latter exhibiting `affine' behaviour
\cite{bridgelandStabilityConditionsKleinian2009,hiranoStabilityConditions3fold2023}.

Hara--Wemyss show that a global classification of hearts (of t-structures) and bricks in
\(\mathcal{C}\) is indeed possible. The classification of spherical objects
(when they exist) comes as a corollary. Up to well-understood mutation
equivalences, the only hearts are the finitely subcategories
\(H_1,...,H_n\subset \mathcal{C}\) given by the homological minimal model
programme \cite{wemyssFlopsClustersHomological2018}, and each is described as
the module category of some finite dimensional algebra
\(\Lambda_{i,\text{con}}\) \((i=1,...,n)\). Further every brick in
\(\mathcal{C}\) is the track of some simple \(\Lambda_{i,\text{con}}\)-module.

\subsection*{Hearts on the `affine' category}
\fakesubsection{Hearts on the `affine' category}
The category \(\Dm X\), too, has algebraic hearts. Van den Bergh
\cite{vandenberghThreedimensionalFlopsNoncommutative2004} observed that there is
a module-finite \(R\)-algebra \(\Lambda\) and a derived equivalence
\(\VdB\colon \Db X\to \Db \Lambda\) that identifies \(\rflmod\Lambda\) with the
full subcategory of \textit{perverse sheaves} \(\zeroper\XZ\subset \Dm X\). This
serves as our reference heart. The homological minimal model program extends
this to a family of \(R\)-algebras \(\{{}_\nu\Lambda_\nu\;|\; \nu\text{ a
sequence of mutations}\}\) called \emph{modification algebras} that are derived
equivalent to \(\Lambda\)  (see \cref{subsec:brennerbutler}).  Tracking the natural
hearts \(\rflmod {}_\nu \Lambda_\nu\) across these equivalences then produces
more algebraic t-structures on \(\Dm X\) which we say are \emph{mutations of
\(\zeroper\XZ\)}.

However it is evident that there can be no autoequivalence which identifies
the above algebraic t-structures with the geometric (and
non-Artinian) heart \(\coh X\). Other geometric hearts can be constructed
by iteratively flopping exceptional curves in \(X\) to obtain new flopping
contractions \(\pi\colon W\to Z\), we say the 3-fold \(W\) thus obtained is a
\emph{birational model} of \(X\). Indeed flops give derived equivalences
\cite{bridgelandFlopsDerivedCategories2002,chenFlopsEquivalencesDerived2002}, so
tracking \(\coh W\) across the Bridgeland--Chen functor \(\Psi\colon \Db
W\to \Db X\) gives a geometric t-structure on \(\Dm X\).

To complicate matters further, there are hearts which are `algebro-geometric' in quite
the literal sense. The construction of the perverse heart is local, so any
crepant morphism \(\tau\colon X\to Y\) contracting some (but not all) of the
\(\pi\)-exceptional curves gives a sheaf \(\ulLambda\) of coherent
\(\OO_{Y}\)-algebras and a derived equivalence \(\Db X \to \Db \ulLambda\). This
identifies a category of suitable \(\ulLambda\)-modules with a heart
\(\zeroper\XY\subset \Dm X\) whose objects are precisely complexes that look
like finite length \(\Lambda\)-modules on the \(\tau\)-exceptional locus and
like \(\OO_X\)-modules elsewhere (\cref{thm:structureofzeroperxy}). The
category is again not equivalent to any heart examined previously, and it is
hence apparent that any classification must also account for these (semi-)geometric
hearts.

\subsection*{Classification of intermediate hearts}
\fakesubsection{Classification of intermediate hearts}
The main result of this paper shows that, at least for sufficiently small
cohomological spread, the above possibilities are in fact exhaustive.

\begin{minipage}[c]{0.7\textwidth}
This involves building semi-geometric
hearts locally using algebraic and geometric categories for `smaller' flopping
contractions (\cref{subsec:partialperverse}), thus obtaining an inductive
description of t-structures.  Indeed any partial contraction \(\tau:X\to Y\) is
an isomorphism away from the open locus \(X^\circ\) of its non-exceptional
fibers, while in any sufficiently small neighbourhood of a \(\tau\)-exceptional
fiber \(C_I\), the map \(\tau\) restricts to a flopping contraction \(X_I\to
Z_I\) with complete local base \(Z_I\) in \(Y\). This gives a (flat) cover of
\(X\). We show any intermediate t-structure (i.e.\ a t-structure whose heart is
concentrated in cohomological degrees \(0\) and \(-1\)) decomposes into purely
algebraic or purely geometric hearts with respect to some such cover.
\end{minipage}
\begin{minipage}[c]{0.3\textwidth}
    \begin{annotationimage}{width=\textwidth}{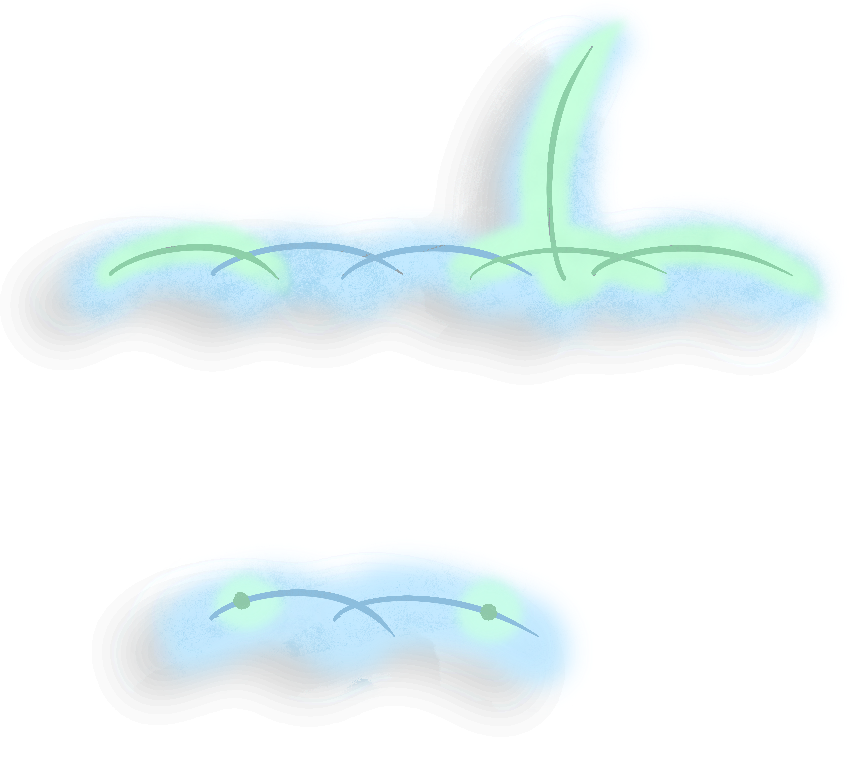}

        \imagelabelset{coordinate label font=\scriptsize}
        \draw[coordinate label={\color{cyan}\(X^\circ\) at (0.45,0.6)}];
        \draw[coordinate label={\(\color{JungleGreen}X_{1}\) at (0.2,0.6)}];
        \draw[coordinate label={\(\color{JungleGreen}X_{456}\) at (0.7,0.6)}];
        \draw[coordinate label={\color{cyan}\(X^\circ\) at (0.45,0.14)}];
        \draw[coordinate label={\(\color{JungleGreen}Z_{1}\) at (0.3,0.14)}];
        \draw[coordinate label={\(\color{JungleGreen}Z_{456}\) at (0.6,0.14)}];

        \imagelabelset{coordinate label font=\tiny}
        \draw[coordinate label={\(\color{JungleGreen}1\) at (0.22,0.7)}];
        \draw[coordinate label={\(\color{cyan}2\) at (0.35,0.7)}];
        \draw[coordinate label={\(\color{cyan}3\) at (0.5,0.7)}];
        \draw[coordinate label={\(\color{JungleGreen}4\) at (0.68,0.7)}];
        \draw[coordinate label={\(\color{JungleGreen}5\) at (0.8,0.7)}];
        \draw[coordinate label={\(\color{JungleGreen}6\) at (0.66,0.8)}];

        \draw[->] (0.45,0.5)--(0.45,0.3) node[midway, right] {\scriptsize\(\tau\)};
    \end{annotationimage}
\end{minipage}

\begin{keytheorem}[(=~\ref{thmC},~\ref{prop:partialperverseinterval})]\label{thm:main}
    Let \(K\subset \Dm X\) be the heart of a t-structure that is
    contained in \(\zeroper\XZ[-1,0]\). Then there is a birational model \(W\) of
    \(X\), and a partial contraction \(\tau\colon W\to Y\), such that \(K\) satisfies
    the following after being pulled back across the flop 
    equivalence \(\Psi\colon \Dm W \to \Dm X\).
    \begin{enumerate}[(1)]
        \item On the locus \(W^\circ\subset W\) where \(\tau\) is an
            isomorphism, \(\Psi^{-1}K\) restricts to the category of coherent
            sheaves \(\coh W^\circ\), possibly tilted in skyscraper sheaves
            \(\langle \OO_p \;|\; p\in Q\rangle\) for some subset of closed
            points \(Q\subset W^\circ\).
        \item For each \(\tau\)-exceptional fiber \(C_I\) in \(W\),
            for sufficiently small neighbourhoods \(W_I\supset C_I\) and
            \(Z_I\ni \tau(C_I)\), the restriction of \(\Psi^{-1}K\)
            to a \(W_I\) is an (algebraic) mutation of the
            category \(\zeroper\pair{W_I}{Z_I} \subset \Dm W_I\).
    \end{enumerate}
    Moreover, the above data uniquely determines \(K\).
\end{keytheorem}

In particular if \(\tau\) above contracts all exceptional curves, then \(K\) is an
algebraic mutation of \(\zeroper\WZ\) (hence also of \(\zeroper\XZ\)) and thus
\(K\) is the category of finite length modules over a modification algebra. On
the other extreme if \(\tau\) contracts no curves then \(K\) is a geometric
t-structure on the birational model \(W\).

\subsection*{Classification of bricks}
\fakesubsection{Classification of bricks}
As a corollary we get a succinct description of bricks in \(\zeroper\XZ\), i.e.\
\(\pi\)-perverse sheaves whose endomorphism algebra is one-dimensional. This
extends to modification algebras a result of Crawley-Boevey \cite[lemma
1]{crawley-boeveyExceptionalFibresKleinian2000}, who showed that the dimension
vector of any brick-module over an affine preprojective algebra is a root.

\begin{keytheorem}[(=~\ref{thm:brickclassification})]\label{thmB}
    Let \(b\in \zeroper\XZ\) be a brick. Then \(b\) is either
    \begin{enumerate}[(1)]
        \item a simple module over some modification algebra
            \({}_\nu\Lambda_\nu\), or
        \item a skyscraper sheaf on some birational model \(W\) of \(X\),
    \end{enumerate}
    tracked under appropriate equivalences induced by mutations or flops. In
    particular, the class of \(b\) in the Grothendieck group is a primitive restricted
    root of the (affine) Dynkin data associated to \(X\).
\end{keytheorem}

In sufficiently restricted settings (e.g.\ under the assumptions imposed in
\cite{todaStabilityConditionsCrepant2008}), every algebraic mutation of
\(\zeroper\XZ\) is equivalent to the category of perverse sheaves on some
birational model \(W\). Then the above result says that, up to applying  mutation
functors, every brick in \(\zeroper\XZ\) is either a point sheaf \(\OO_p\), or
the twisted structure sheaf \(\OO_{C_i}(-1)\) on some integral curve, or the
suspended canonical sheaf \(\omega_{\ul C}[1]\) of the scheme-theoretic
exceptional fiber on some birational model \(W\). In general however there are
modifying algebras not related to birational geometry in any obvious way
(these are the `hidden' t-structures alluded to
in~\cite{hiranoStabilityConditions3fold2023}), and expressing their simple
modules in geometric terms is difficult. For a single flopping curve this is
accomplished by Donovan--Wemyss \cite[\S 4]{donovanStringyKahlerModuli2024} who
show that each such simple is, up to mutation, determined by the structure sheaf
or the canonical sheaf of some thickening of the exceptional curve.

\subsection*{Canonical surfaces and preprojective algebras}
\fakesubsection{Canonical surfaces and preprojective algebras}
Suppose \(\ol X\) is a Gorenstein canonical surface that admits a crepant birational
morphism \(\pi\colon \ol X \to \ol Z = \Spec(\ol R, \mathfrak{p})\) with complete
local base. The results of this paper, though developed in the context of
3-folds, apply verbatim to the derived category of \(\ol X\) and its
\(\pi\)-perverse heart provided one suitably reinterprets the notions of
modification algebras and birational models.

Such a map \(\pi\) is necessarily a crepant (partial) resolution of the Kleinian
singularity \(\ol Z\), i.e.\ \(\ol X\) is obtained by contracting some
exceptional curves in the minimal resolution \(\wt X\to \ol Z\). It is
well-known that \(\wt X\) is derived equivalent to an affine preprojective
algebra \(\Pi\); the contraction \(\wt X \to \ol X\) determines an idempotent
\(e\in \Pi\) such that \(\ol X\) is derived equivalent to the \emph{contracted
preprojective algebra} \(e\Pi e\) \cite{kalckFrobeniusCategoriesGorenstein2015}.
This equivalence \(\Dfl(e\Pi e)\to \Dm (\ol X)\) is also recovered by
Van den Bergh's construction
\cite{vandenberghThreedimensionalFlopsNoncommutative2004}, and algebraic
mutations of the category of \(\pi\)-perverse sheaves can be read off from the
tilting theory of \(e\Pi e\) \cite[\S 7.4]{iyamaTitsConeIntersections}.

The discussion of birational models is only slightly more subtle, as curves in
surfaces do not flop. Birational transformations (including flops in dimension
3), on the other hand, are naturally induced via geometric invariant theory.
Indeed \(\wt X\) appears as a moduli space \(\mathcal{M}(\theta,\delta)\) of
stable \(\Pi\)-modules, for some generic stability parameter \(\theta\)
relative to a \(\KK\)-theory class \(\delta\)
\cite{cassensKleinianSingularitiesQuivers1998}. Variation of GIT parameter
\(\theta_1\rightsquigarrow \theta_2\) produces a birational map
\(\mathcal{M}(\theta_1,\delta)\dashrightarrow \mathcal{M}(\theta_2,\delta)\) of
\(\ol Z\)-schemes. In particular if \(\theta_1\) and \(\theta_2\) lie in
adjacent chambers then this birational map is defined away from a single
exceptional curve, and simple wall--crossing of the GIT parameter thus gives an
analogue of flops in dimension 2. The role of Bridgeland--Chen equivalences is
played by \emph{reflection functors}
\cite{sekiyaTiltingTheoreticalApproach2013}, and the results are extended to the
contracted setting by Iyama--Wemyss \cite{iyamaTitsConeIntersections}.

With this, all arguments in \cref{sec:numtors,sec:covlimiting} regarding the
partial order of algebraic intermediate hearts and its interactions with convex
geometry and line bundles on birational models hold verbatim. Thus the reader
indifferent to dimensions can safely read the remainder of this paper as if it
were written for surfaces.

To the reader interested in multiple dimensions, we remark that the
correspondence between our results for 3-folds and surfaces is natural. Indeed,
the morphism \(\ol X \to \ol Z\) can be embedded as a generic hyperplane section
(\emph{general elephant}) of some 3-fold flopping contraction \(X\to
\Spec(R,\mathfrak{m})\), and the corresponding non-commutative algebras are
related by the reduction \( e\Pi e= \Lambda\otimes_R \ol R\). By reducing both
to the fiber \(R/\mathfrak{m} = \ol R / \mathfrak{p}\), Kimura
\cite[theorem 5.4]{kimuraTiltingSiltingTheory2024} shows that the corresponding
functor \(\rflmod (e\Pi e) \into \rflmod \Lambda\) induces a bijection of
torsion classes. Likewise lemma 5.1 \emph{ibid.} gives the correspondence
between bricks.

\subsection*{Mutating v/s tilting}
\fakesubsection{Mutating v/s tilting}
To motivate why a result such as \cref{thm:main} is desirable towards a full
classification, we briefly sketch the key argument of
\cite{haraSphericalObjectsDimensions2024}. Each algebra
\(\Lambda_{i,\text{con}}\) produced by the homological minimal model programme
is silting discrete \cite{augustTiltingTheoryContraction2020}, thus can be
assigned a finite and complete hyperplane arrangement, the \emph{silting fan}.
The chambers \(\sigma_1,...,\sigma_n\) of this arrangement are in bijection with
the hearts \(H_i\subset \mathcal{C}\), and each minimal sequence of wall--crossings
\(\sigma_j\rightsquigarrow \sigma_i\) is assigned an atomic mutation functor
\(\Phi_{ij}\). Then, say, given an object \(x\in H_i[-n,0]\) with non-zero
cohomologies in degrees \(0\) and \(n(>\!0)\), Hara--Wemyss consider the set of
paths \(\sigma_j\rightsquigarrow \sigma_i\) such that \(\Phi_{ij}(x)\) lies in
\(H_j[-n,0]\) and show that for the longest such path
\(\sigma_{j_0}\rightsquigarrow \sigma_i\), the object \(\Phi_{ij_0}(x)\) in fact
lies in \(H_{j_0}[-n+1,0]\). Iterating shows \(x\) lies in some \(H_j\) after a
finite sequence of mutations.

Finiteness of the silting fan underpins the argument,
guaranteeing that any poset considered has maximal elements. The
affine curse ensures we don't have such privileges when studying
\(\Dm X\), and the silting fan of \(\Lambda\) is infinite and incomplete. So
while there still is an assignment of algebraic hearts and mutation functors to
chambers and wall-crossings (\cref{thm:MMRDynkinlabelling}), the fan has
infinite paths which culminate outside its support.

We thus replace mutation and the silting fan with \emph{tilting in torsion pairs}
and the \emph{heart fan} respectively. Given any Abelian category \(H\),
Happel--Reiten--Smal\o~\cite{happelTiltingAbelianCategories1996} show that the
set \(\tilt(H)\) of intermediate t-structures (i.e.\ t-structures with hearts
contained in \(H[-1,0]\)) can be recovered by tilting \(H\) in its torsion
subcategories, and the containment order of torsion classes enhances
\(\tilt(H)\) to a \emph{complete lattice} (i.e.\ a poset in which every subset
has an infimum and a supremum).
Broomhead--Pauksztello--Ploog--Woolf~\cite{broomheadHeartFanAbelian2024} assign
each heart \(K\in \tilt(H)\) to a cone \(\CC K\) and show that the ensemble
\(\HFan(H)\) of heart cones is a fan. This heart fan is complete and simplicial
for categories such as \(\zeroper\XZ\), and the silting fan then sits inside
\(\HFan(H)\) as a sub-fan. Relating atomic mutations to tilting in
\emph{functorially finite} torsion classes, we thus have a `completion' of
silting theory.

Giving a complete description of the lattice of torsion classes is hard, even
for finite dimensional
algebras~\cite{thomasIntroductionLatticeTorsion2021,demonetLatticeTheoryTorsion2023}.
The heart fan gives some insight into tilts satisfying numerical criteria
(\cref{thm:allheartsonacone}), but there are
examples of tilts of algebraic hearts which have trivial (i.e.\ \(\mathbf{0}\))
heart cone (see \cite[example 6.7]{broomheadHeartFanAbelian2024}). In such
situations the numerical criteria lead to tautologies, and we are left to our
own devices.

The main result of this paper, then, is that such situations do not arise when
studying tilts of \(\zeroper\XZ\). Consequently a complete description of the
lattice of torsion classes is possible.

\begin{keytheorem}[(=~\ref{cor:heartfanofzeroper})]\label{thmC}
    The heart fan of \(H=\zeroper\XZ\) is given by an intersection arrangement
    associated to the restriction of an affine Dynkin root system, i.e.\ one of
    the arrangements described in
    \textnormal{\cite{iyamaTitsConeIntersections}}. The heart cones are
    described as follows.
    \begin{enumerate}[(1)]
        \item \textnormal{(= \ref{thm:everyheartisnumerical})}\label{item:thmC1}
            The trivial cone \(\mathbf{0}\) is \emph{not} the heart cone of any
            intermediate heart.
        \item \textnormal{(= \ref{cor:heartconesinsiltingfan})}\label{item:thmC2}
            A cone outside the imaginary-root hyperplane
            \(\{\delta=0\}\) is a heart cone if and only if it is
            full-dimensional, and in this case it is the heart cone of a unique
            algebraic heart in \(\tilt(H)\). The positive and negative orthants
            \(\Chamb^+,\Chamb^-\) are the heart cones of \(H\) and
            \(H[-1]\) respectively, and the heart associated to any other full
            dimensional cone can be expressed as a mutation of \(\zeroper\XZ\)
            by choosing wall--crossing paths from \(C^\pm\).

            The hearts \(H\) and \(H[-1]\) are the maximal and minimal elements
            in \(\tilt(H)\) respectively, and the partial order on remaining
            algebraic hearts respects (atomic) wall--crossing distance from these.

        \item
            \textnormal{(=~\ref{lem:ktheoryofskyscraper},~\ref{prop:heartconeofpartialperverse},~\ref{prop:intermediacyofpartialperverse})}
            \label{item:thmC3}
            The induced finite hyperplane arrangement on \(\{\delta=0\}\) is
            naturally identified with the (real) Picard group of \(X\), in a way
            that every maximal cone \(\sigma\) in this subfan thus parametrises nef divisors on a
            unique birational model \(W\) of \(X\). Such \(\sigma\) is then the
            heart cone of \(\coh W\) tracked under the composite
            Bridgeland--Chen functor \(\Dm W \to \Dm X\).
        \item \textnormal{(= \ref{thm:nonmaximalheartcones}, \ref{prop:partialperverseinterval})}
            \label{item:thmC4}
            More generally, every non--trivial cone
            \(\sigma\subset\{\delta=0\}\) can be assigned a unique partial
            contraction \(W\to Y\) such that \(\zeroper\WY\) (tracked under
            flop equivalences) is the maximal heart in
            \(\tilt(H)\) with heart cone \(\sigma\). Every other heart with
            heart cone \(\sigma\) can
            be obtained by arbitrarily mutating the algebraic components of
            \(\zeroper\XY\) and tilting in skyscrapers in the geometric
            components.

            An algebraic intermediate heart \(K\) then satisfies
            \(K>\zeroper\WY\) if and only if there is an atomic path
            \(\Chamb^+\!\rightsquigarrow\CC K\) which can be atomically extended
            to an infinite sequence of wall--crossings
            \(\Chamb^+\!\rightsquigarrow\CC K \rightsquigarrow \CC K_1
            \rightsquigarrow \CC K_2\rightsquigarrow \cdot\cdot\cdot\)  with
            generic limit in \(\sigma\). In this case \(\zeroper\XY\) is the
            infimum of the decreasing chain \(K>K_1>K_2>...\) in \(\tilt(H)\).
    \end{enumerate}
\end{keytheorem}

In \cref{sec:numtors} we establish the numerical story, enumerating for each cone in
the intersection arrangement all the intermediate hearts and King--semistable
objects associated to it. This involves establishing the fact that coherent
sheaves on any birational model \(W\) (which are a priori only intermediate with
respect to the perverse heart on \(W\)) are in fact intermediate with respect to
the perverse heart on \(X\). For this a careful analysis of the partial order is
necessary, and convex geometric tools are needed not only to supply a
systematic enumeration scheme but also to establish crucial limiting results as
in~\cref{thmC}~\ref{item:thmC4}.

Once the heart fan is filled in, the question of whether there are any hearts
hidden away in the \(\mathbf{0}\)--cone remains. Here the affine blessing
ensures that the non-algebraic locus in the heart fan is codimension 1, so
every heart in \(\tilt(H)\) has tight algebraic bounds. 

The following example is illustrative.

\subsection*{Single-curve flops}
\fakesubsection{Single-curve flops}
Let \(R=\mathbb{C}\llbracket u,v,x,y\rrbracket / (uv-xy)\) be the
\(\text{cA}_1\) singularity, i.e.\ the base of the Atiyah flop
\(X\dashrightarrow W\) where \(X\) is the neighbourhood of a \((-1,-1)\)
rational curve \(C\). The two flopping contractions \(X,W\to Z=\Spec R\) are
obtained by blowing up the ideals \((u,x)\) and \((u,y)\) respectively.

Accordingly the modification algebra associated to
\(X\) is \(\Lambda = \End(R\oplus (u,x))\). It can be shown that all
modification algebras are isomorphic and are related by derived equivalences as
below.
\[\begin{tikzcd}[column sep=3em, row sep=3em]
    \cdot\cdot\cdot
    \arrow[r, "\Psi_1", shift left=1] &
    \Dfl \End\left(\!\!\begin{smallmatrix}L \\\oplus\\(u,y)\end{smallmatrix}\!\!\right)
    \arrow[r, "\Psi_0", shift left=1]
    \arrow[l, "\Psi_1", shift left=1]&
    \Dfl\End\left(\!\!\begin{smallmatrix}R \\\oplus\\(u,y)\end{smallmatrix}\!\!\right)
    \arrow[r, "\Psi_1", shift left=1]
    \arrow[l, "\Psi_0", shift left=1]&
    \Dfl\End\left(\!\!\begin{smallmatrix}R \\\oplus\\(u,x)\end{smallmatrix}\!\!\right)
    \arrow[r, "\Psi_0", shift left=1]
    \arrow[l, "\Psi_1", shift left=1]&
    \Dfl\End\left(\!\!\begin{smallmatrix}K \\\oplus\\(u,x)\end{smallmatrix}\!\!\right)
    \arrow[r, "\Psi_1", shift left=1]
    \arrow[l, "\Psi_0", shift left=1]&
    \cdot \cdot \cdot
    \arrow[l, "\Psi_1", shift left=1] \\
                                     &&
    \Dm W \arrow[u,"\text{VdB}"]     &
    \Dm X \arrow[u,"\text{VdB}"]     &&
\end{tikzcd}\]
Here the functors \(\Psi_i\colon \Dfl(\End N) \to \Dfl(\End M)\) are given as
\(\RHom(\Hom(N,M),-)\), and the index records which summand changes between
\(N\) and \(M\). This changed summand is computed via \emph{mutation}, thus for
instance \((u,x)\) and \((u,y)\) are (first) syzygies of each other while
\(L\subset (u,y)^{\oplus 2}\) is the kernel of the natural map \(((u,y)\into
R)\oplus((u,y)\isoto (v,x)\into R)\).

\begin{figure}[H]
    \begin{minipage}[c]{0.33\textwidth}
        \caption{
            The heart fan for a minimal resolution of the
            \(\text{cA}_1\) singularity.\newline\newline
            Hyperplanes are induced by the \(\widetilde{\text{A}}_1\) root
            system, where the simple real roots \(\alpha_0,\alpha_1\) are
            identified with the \(\KK\)-theory classes of \(\omega_{C}[1]\) and
            \(\OO_C(-1)\) respectively.
        } \label{fig:AffA1}
    \end{minipage}\qquad
    \begin{minipage}[c]{0.36\textwidth}
        \begin{annotationimage}{width=\textwidth}{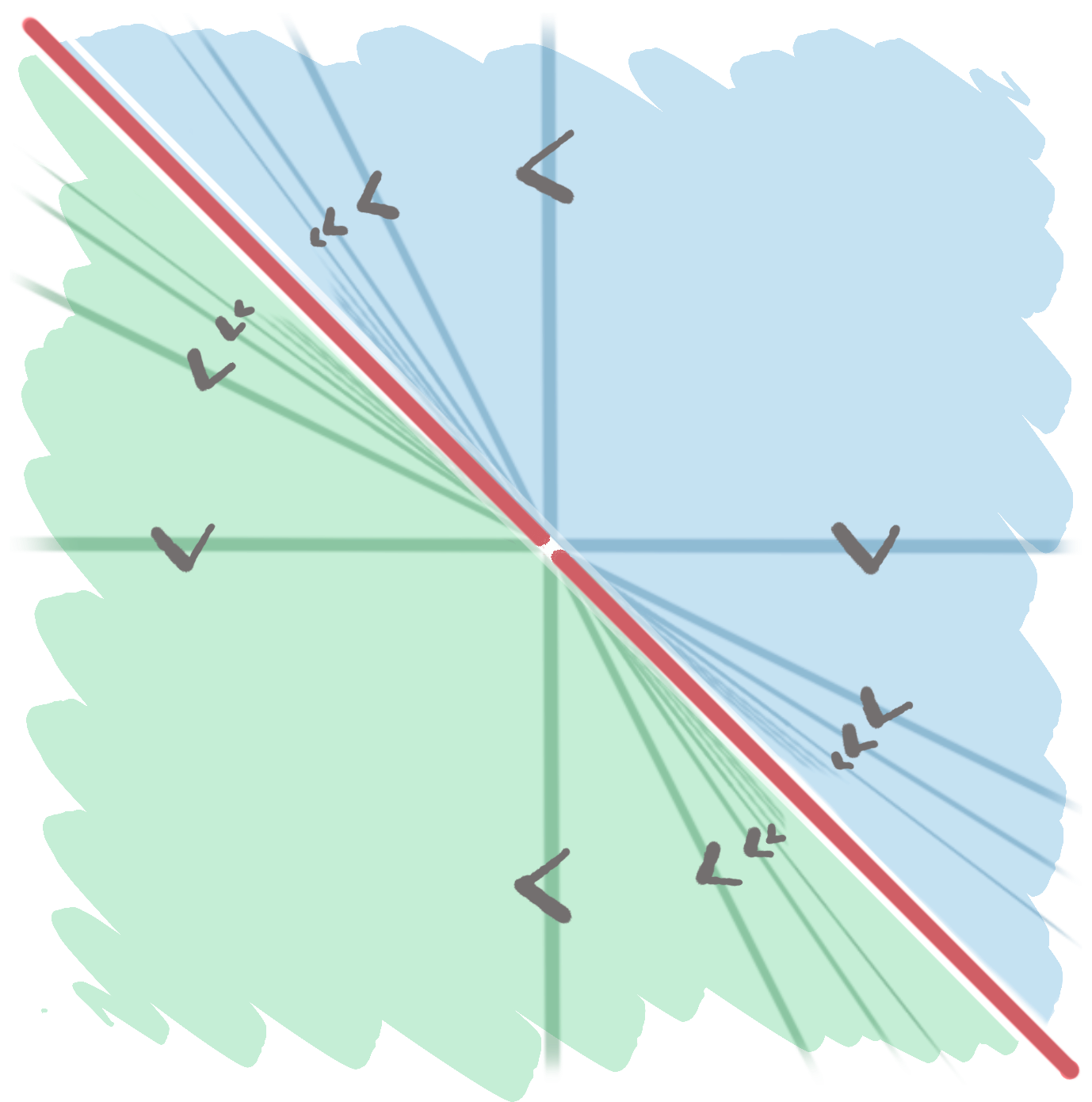}
            \draw[coordinate label = {\(H[-1]\) at (0.3,0.3)}];
            \draw[coordinate label = {\(H\) at (0.7,0.7)}];
            \draw[coordinate label = {\(\Psi_0 H\) at (0.82,0.425)}];
            \draw[coordinate label = {\(\Phi_0 H[-1]\) at (0.2,0.57)}];
            \draw[annotation left = {\(\Phi_1 H[-1]\) at 0.25}]
                to (0.57,0.25);
            \draw[annotation left = {\(\Phi_0\Phi_1 H[-1]\) at 0.15}]
                to (0.73,0.15);
            \draw[annotation left = {\(\Phi_1\Phi_0\Phi_1 H[-1]\) at 0.05}]
                to (0.84,0.05);
            \draw[annotation left = {\(\Phi_1\Phi_0 H[-1]\) at 0.7}]
                to (0.17,0.7);
            \draw[annotation left = {\(\Phi_0\Phi_1\Phi_0 H[-1]\) at 0.8}]
                to (0.09,0.8);
            \draw[annotation right = {\(\Psi_1 H\) at 0.75}]
                to (0.45,0.75);
            \draw[annotation right = {\(\Psi_0\Psi_1 H\) at 0.85}]
                to (0.3,0.85);
            \draw[annotation right = {\(\Psi_1\Psi_0\Psi_1 H\) at 0.95}]
                to (0.175,0.95);
            \draw[annotation right = {\(\Psi_1\Psi_0 H\) at 0.3}]
                to (0.85,0.3);
            \draw[annotation right = {\(\Psi_0\Psi_1\Psi_0 H\) at 0.2}]
                to (0.95,0.2);
            \draw[coordinate label = {{\color{gray} \scriptsize\(\{\alpha_0=0\}\)}
                at (0.98,0.53)}];
            \draw[coordinate label = {\rotatebox{-90}{\color{gray} \scriptsize\(\{\alpha_1=0\}\)}
                at (0.53,0.04)}];
            \draw[coordinate label = {\rotatebox{-45}{\color{gray} \scriptsize\(\{\delta=0\}\)}
                at (1,0.06)}];
        \end{annotationimage}
    \end{minipage}
\end{figure}

Any category \(\Dfl(\End M)\) has a natural heart \(H_M=\rflmod(\End M)\), and
tracking this across equivalences appearing in the above diagram (and their
inverses) produces new t-structures in \(\Dm X\). Such a t-structure is
intermediate with respect to the perverse heart when the chain of
functors \(\Dfl(\End M)\rightsquigarrow \Dm X\) is chosen as short as possible,
i.e.\ the tracks of \(H_M\) that lie in \(\zeroper\XZ[-1,0]\) are of the form
\[
    \VdB^{-1}\circ (...\circ \Psi_0\circ \Psi_1 \circ \Psi_0 \circ
    ...) H_M  \quad \text{or} \qquad
    \VdB^{-1}\circ (...\circ \Psi_0\circ \Psi_1 \circ \Psi_0 \circ
    ...) \!{}^{-1}\; H_M[-1].
\]
Note the specified chain completely determines the domain of the functor, so we
may drop the subscript for brevity. We also drop \(\VdB\) from the
notation, so for instance \(\Psi_1 H\) is the image of \(\rflmod(R\oplus
(u,y))\) under the functor \(\VdB^{-1}\circ \Psi_1\). The algebraic hearts
in \(\Dm X\) thus obtained are in natural bijection with full-dimensional cones
in the \(\widetilde{\text{A}}_1\) (\(\dynkin[arrows=false]B{o*}\)) intersection
arrangement, see \cref{fig:AffA1}.

The ray \(\{\delta=0,\alpha_0\geq 0\}\), geometrically the `limit' of the path
\(\CC(H)\dashrightarrow \CC(\Psi_0 H)\dashrightarrow \CC(\Psi_1\Psi_0
H)\dashrightarrow\cdots\),
is the heart cone of \(\coh X =\, \inf\, \{H,\,\Psi_0 H,\, \Psi_1\Psi_0 H,\,
\ldots\}\). Likewise it is also the heart cone of the \emph{reversed geometric heart}
\(\anticoh X = \sup\, \{H[-1], \Phi_1 H[-1], \Phi_0\Phi_1 H[-1],\ldots\}\) which can
be obtained by tilting \(\coh X\) in the torsion class generated by all
skyscrapers.

The ray \(\{\delta=0, \alpha_0\leq 0\}\) is similarly seen to be the heart cone of
\(\flop(\coh W)\) and other geometric hearts that live on \(W\), where
\(\flop\colon \Dm W\to \Dm X\) is the Bridgeland--Chen flop functor (in this
case isomorphic to \(\VdB^{-1}\circ \Psi_1 \circ \VdB\)).

Now the partial order on algebraic hearts (indicated in \cref{fig:AffA1}) is
such that wall--crossings correspond to simple tilts
(\cref{subsec:simplemutations}), thus \(\Psi_0 H\) is obtained by tilting \(H\)
in the torsion class \(\langle \omega_C[1] \rangle\) while \(\Phi_1H[-1]\) and
\(H[-1]\) are related by a tilt in \(\langle \OO_{C}(-1) \rangle\). Considering
how simples sit in relation to torsion pairs then allows us to `push'
non-algebraic intermediate hearts \(K\in \tilt(H)\) towards the geometric hyperplane.

\begin{algorithm}[Simple--tilting]
    \label{alg:naivetilting}
    Let \(K\) be a non-algebraic intermediate heart. Write \(H=T\ast F\) for the torsion pair corresponding to \(K(=F\ast T[-1])\).
    \begin{enumerate}[1.]
        \item  Since \(K\neq H\), one of the simples of \(H\) (say
            \(S_0=\omega_C[1]\)) lies in \(T\) and this forces the inequality
            \(K<\Psi_0 H\).

        \item The inequality \(\Psi_0 H > K\) shows that \(T\) intersects
            \(\Psi_0 H\) nontrivially, hence must contain some simple \(\{\Psi_0
            S_0, \Psi_0 S_1\}\) of \(\Psi_0 H\). Since \(\Psi_0 S_0=\omega_C\) lies
            outside \(H\), we must have \(\Psi_0 S_1 \in T\) and thus
            \(\Psi_1\Psi_0H > K\). Iterating, we have \(K< ... \Psi_0\Psi_1\Psi_0
            H\) for arbitrarily long chains, and hence \(K\leq \coh X\).

        \item On the other hand \(K\neq H[-1]\) so \(F\) contains some simple of
            \(H\), which in this case must be \(S_1=\OO_C(-1)\). This shows
            \(K>\Phi_1H[-1]\), and iterating as above gives \(K \geq \anticoh X\).

        \item  Since \(\anticoh X\) and \(\coh X\) both share the heart cone
            \(\{\delta=0,\alpha_0\geq 0\}\), \(K\) must do so too, and in
            particular \(K\) is numerical. In fact \(K\) can be expressed as a
            tilt of \(\coh X\) in \(\langle \OO_p\;|\; p\in Q \rangle\) for some
            \(Q\subset C\).
    \end{enumerate}
\end{algorithm}

Thus every heart in \(\tilt(\zeroper\XZ)\) is either algebraic, or can be
shown to be geometric via the above recipe.

\subsection*{Multi-curve flops}
\begin{figure}[H]
    \begin{minipage}[c]{0.6\textwidth}
        \begin{annotationimage}{width=\textwidth}{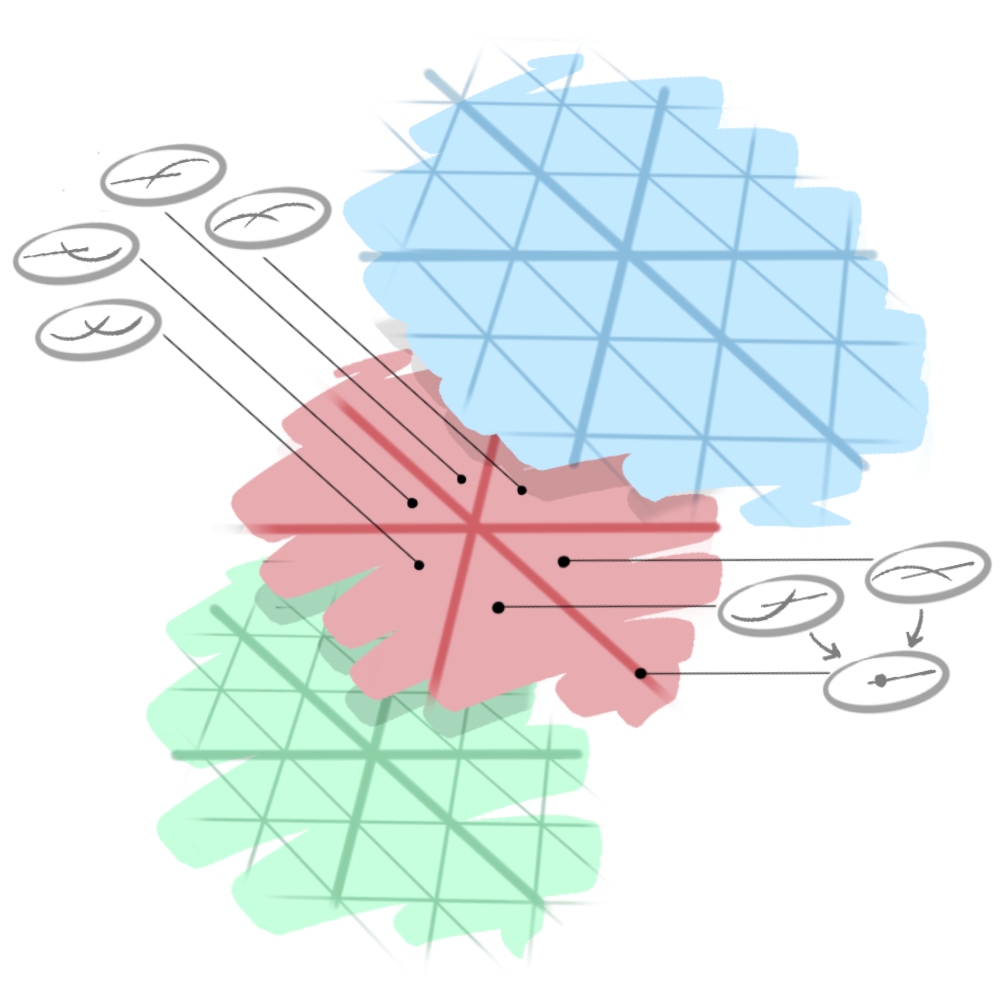}
            \draw[coordinate label = {\(\Psi_2H\) at (0.67,0.76)}];
            \draw[coordinate label = {\(\Psi_{202}H\) at (0.78,0.76)}];
            \draw[coordinate label = {\(\Psi_{212}H\) at (0.56,0.76)}];
            \draw[coordinate label = {\(H\) at (0.69,0.725)}];
            \draw[coordinate label = {\(\Psi_{21}H\) at (0.575,0.725)}];
            \draw[coordinate label = {\(\Psi_{12}H\) at (0.6,0.81)}];
            \draw[coordinate label = {\(\Psi_{02}H\) at (0.7,0.81)}];
            \draw[coordinate label = {\(\Psi_{012}H\) at (0.585,0.84)}];
            \draw[coordinate label = {\(\Psi_1H\) at (0.655,0.675)}];
            \draw[coordinate label = {\(\Psi_{0}H\) at (0.765,0.675)}];
            \draw[coordinate label = {\(\Psi_{01}H\) at (0.79,0.645)}];

            \draw[coordinate label = {\(H[-1]\) at (0.32,0.265)}];
            \imagelabelset{outer dist=-4em}
            \draw[annotation left = {\(\Phi_0 H[-1]\) at 0.29}]
                to (0.27,0.29);
            \draw[annotation left = {\(\Phi_2 H[-1]\) at 0.22}]
                to (0.335,0.22);
            \draw[annotation left = {\(\Phi_{102} H[-1]\) at 0.15}]
                to (0.32,0.15);

            \draw[coordinate label = {{\color{gray} \scriptsize\(\{\delta=0\}\)}
                at (0.2,0.5)}];
            \draw[coordinate label = {{\color{gray} \scriptsize\(\{\delta+1=0\}\)}
                at (0.14,0.4)}];
            \draw[coordinate label = {{\color{gray} \scriptsize\(\{\delta-1=0\}\)}
                at (0.4,0.9)}];
            \draw[coordinate label = {\normalsize\(\sigma\) at (0.685,0.3)}];

            \draw[coordinate label = {{\color{magenta} \scriptsize\(\nu_2X\)}
                at (0.3,0.825)}];
            \draw[coordinate label = {{\color{magenta} \scriptsize\(\nu_{12} X\)}
                at (0.2,0.875)}];
            \draw[coordinate label = {{\color{magenta} \scriptsize\(\nu_{212} X\)}
                at (0.06,0.79)}];
            \draw[coordinate label = {{\color{magenta} \scriptsize\(\nu_{21} X\)}
                at (0.1,0.62)}];

            \draw[coordinate label = {{\color{magenta} \scriptsize\(X\)}
                at (0.92,0.47)}];
            \draw[coordinate label = {{\color{magenta} \scriptsize\(\nu_{1} X\)}
                at (0.8,0.44)}];
            \draw[coordinate label = {{\color{magenta} \scriptsize\(Y\)}
                at (0.9,0.27)}];
        \end{annotationimage}
    \end{minipage}\hspace{-1em}
    \begin{minipage}[c]{0.43\textwidth}
        \vspace{7em}
        \begin{annotationimage}{width=\textwidth}{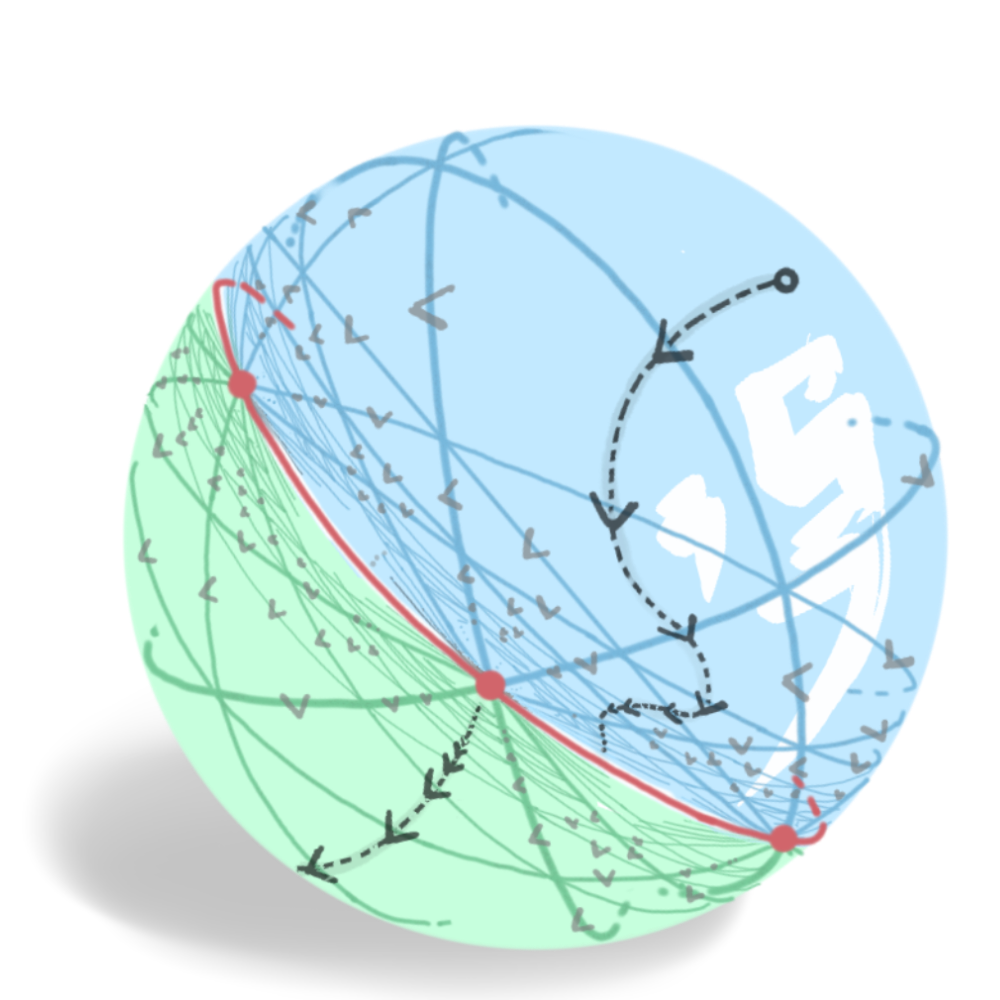}
            \draw[coordinate label = {\(H\) at (0.82,0.7)}];
            \draw[coordinate label = {\(\Psi_1H\) at (0.57,0.6)}];
            \draw[coordinate label = {\(\Psi_0H\) at (0.95,0.4)}];
            \draw[coordinate label = {\normalsize\(\sigma\) at (0.8,0.13)}];
        \end{annotationimage}
    \end{minipage}

    \caption{
        The \(3\)-dimensional heart fan for a \(\text{cA}_2\) resolution, sliced
        along affine hyperplanes (left) and the unit sphere (right). The
        mutation functors are abbreviated, e.g.\ by writing \(\Psi_{01}\) for
        \(\Psi_1\circ \Psi_0\). Cones in the hyperplane \(\{\delta=0\}\) have
        been labelled by the \(3\)-fold whose geometric hearts
        they are associated to, where e.g.\ \(\nu_{212}X\) is the flop of
        the second curve in \(\nu_{12}X\).
    }
    \label{fig:AffA2}
\end{figure}

\fakesubsection{Multi-curve flops}
After accounting for perverse hearts for partial contractions, the above
description of the heart fan and the partial order carries over to higher rank
cases. Thus for a crepant resolution \(X\to Z\) of the \(\text{cA}_2\) singularity
\(Z=\Spec \mathbb{C}\llbracket u,v,x,y \rrbracket/(uv-xy(x+y))\), the heart fan of
\(\zeroper\XZ\) is described by an \(\widetilde{\text{A}}_2\)
(\raisebox{-2pt}{\(\dynkin[edge length=6pt]A[1]2\)}) root system as in
\cref{fig:AffA2}.

Given a non-algebraic heart \(K\in \tilt(\zeroper\XZ)\), we can follow
\cref{alg:naivetilting} to produce bounding sequences
\[\ldots\Psi_{i_2}\Psi_{i_1}H \quad >\quad K \quad > \quad \ldots\Phi_{j_2}\Phi_{j_1}H[-1]\] by iterated simple tilting. However this does not suffice for `pushing \(K\) to the geometric hyperplane', and it is possible that the corresponding paths in the heart fan do not converge
to the same cone (illustrated in \cref{fig:AffA2}). In this regard the
single-curve case is deceptively simple.

Establishing \(\CC K\neq 0\) thus demands more finesse, and our strategy is
to utilise line bundles on birational models for this purpose.
Under the identification of \(\RealPic X\) with the
hyperplane \(\{\delta=0\}\), (the proper transform of) any line bundle \(\LL\in
\Pic W\) corresponds to a vector \(\theta\) in the heart fan, so for any cone
\(\sigma\) we can consider the submonoid \(\Pic W \cap \sigma\) (e.g.\ if
\(\sigma=\CC(\flop(\coh W))\) then this is the monoid of nef bundles). On
the other hand, \(\Pic W\) has a natural action on \(\Dm W\) which induces an
action \(\Pic W \circlearrowright \Dm X\) across the composite flop equivalence.
Using the abbreviated notation \(\LL\otimes H\coloneq \flop (\LL\otimes_W
(\flop^{-1}H))\), we establish the following in \cref{subsec:nefmonoid-dynamics}.

\begin{keytheorem}[(=~\ref{thm:nefintermediacy},~\ref{cor:partialorderoftwists},~\ref{thm:nefcomparisons})]\label{thmD}
    Writing \(H=\zeroper\XZ\) for the standard heart in \(\Dm X\), the
    following statements hold.
    \begin{enumerate}[(1)]
        \item Given any birational model \(W\) of \(X\) and a line bundle
            \(\LL\in \Pic W\), any heart in \(\Dm X\) of the form
            \(\LL\,\dual\otimes H\) or \(\LL\otimes H[-1]\) lies in \(\tilt(H)\)
            if and only if \(\LL\) is nef.
        \item If \(W'\) is another birational model such that
            \(\LL\in \Pic W\) and its proper transform \(\LL'\in \Pic W'\) are
            both nef, then there are equalities of t-structures
            \(\LL\,\dual\otimes H = \LL'\,\dual\otimes H\), \(\LL\otimes
            H[-1]=\LL'\otimes H[-1]\).
        \item For any cone \(\sigma\subset \CC(\textnormal{flop}\coh W)\) in
            \(\HFan(H)\), the induced actions of monoid \(\Pic W\cap \sigma\) on
            the subsets
            \[
                \sigma\text{-}\tilt^+(H)=\;
                \smashoperator[lr]{\bigcup_{\LL\in \Pic W\cap \sigma}}\,
                    \left[\LL\,\dual\otimes H, H\right], \qquad
                \sigma\text{-}\tilt^-(H)=\;
                \smashoperator[lr]{\bigcup_{\LL\in \Pic W\cap \sigma}}\,
                    \left[H[-1], \LL\otimes H[-1]\right] 
            \]
            respect the partial order inherited from \(\tilt(H)\) and the
            monoid-order on \(\Pic W \cap \sigma\).
        \item Each heart in the posets above is algebraic, and conversely every
            algebraic heart in \(\tilt(H)\) lies in some poset of the above form
            for suitable \(W,\sigma\).
        \item The infimum of \(\sigma\)-\(\tilt^+(H)\) is the maximal element of
            \(\tilt(H)\) with heart cone \(\sigma\), likewise the supremum of
            \(\sigma\)-\(\tilt^-(H)\) is the minimal element of \(\tilt(H)\)
            with heart cone \(\sigma\).
    \end{enumerate}
\end{keytheorem}

\begin{figure}[H]
    \begin{minipage}[c]{0.33\textwidth}
        \caption{
            Continuing from \cref{fig:AffA2}, the actions of \(\Pic X\) and
            \(\nu_1 X)\) on \(H\) are shown. Here
            \(\left[\!\begin{smallmatrix}i\\j\end{smallmatrix}\!\right]\)
            denotes the line bundle on \(X\) which has degrees \(i\) and
            \(j\) on the two exceptional curves respectively, and we use double
            brackets (\(\llbracket \!\begin{smallmatrix}i\\j\end{smallmatrix}\!
            \rrbracket\)) for bundles on the flop \(\nu_1 X\).
            \newline \newline
            Note that since
            \(\left[\!\begin{smallmatrix}0\\1\end{smallmatrix}\!\right]\) and
            its proper transform \(\llbracket
            \!\begin{smallmatrix}0\\1\end{smallmatrix}\! \rrbracket\) are both
            trivial on the flopped curve, their actions coincide i.e.\ \(
                \left\llbracket\!\begin{smallmatrix}0\\1\end{smallmatrix}\!\right\rrbracket
                    \dual\!\!\otimes\!H
              = \left[\!\begin{smallmatrix}0\\1\end{smallmatrix}\!\right]
                    \dual\!\!\otimes\!H\).
             \newline \newline
             The shaded region represents \(\sigma\)-\(\tilt^+(H)\) for
         \(\sigma\) shown in \cref{fig:AffA2}.}
        \label{fig:Pic-on-AffA2}
    \end{minipage}\qquad
    \begin{minipage}[c]{0.6\textwidth}
        \begin{annotationimage}{width=\textwidth}{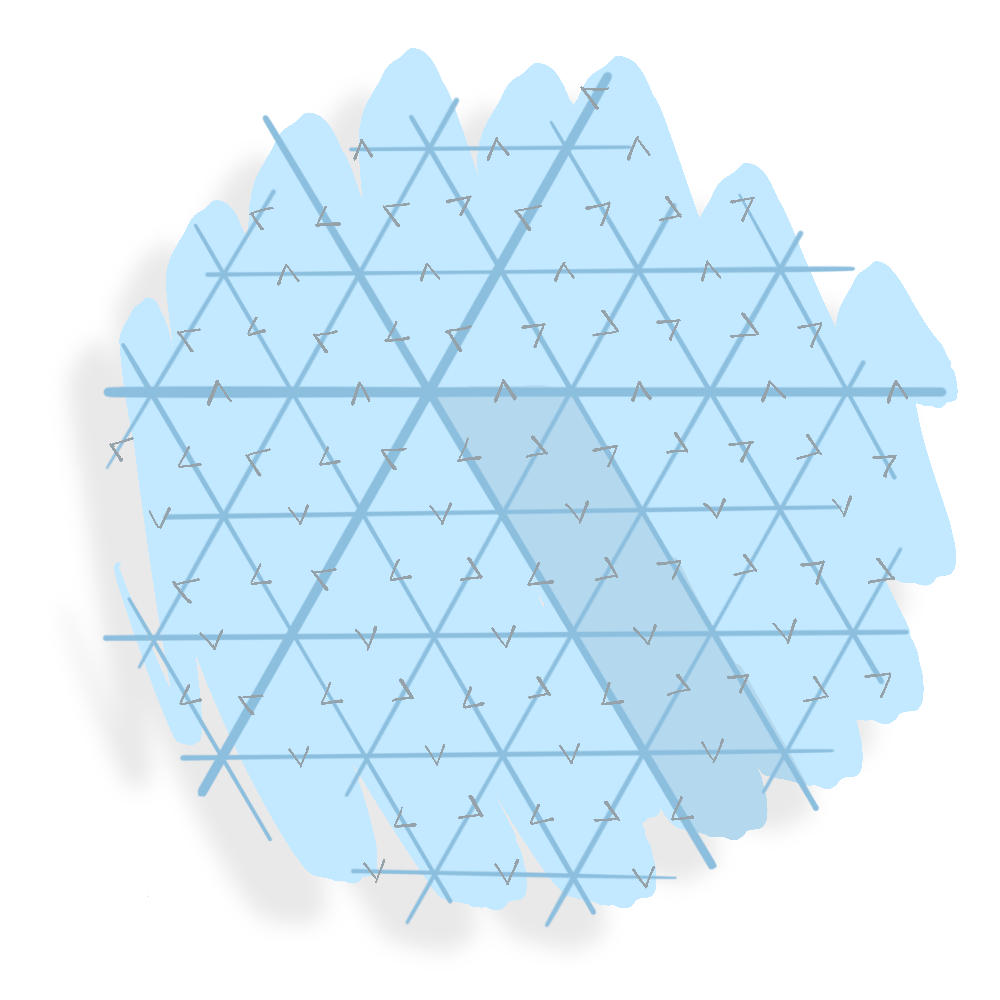}
            \draw[coordinate label = {\(H\) at (0.5,0.575)}];
            \draw[coordinate label = {\(\left[\!\begin{smallmatrix}1\\0\end{smallmatrix}\!\right]\dual\!\!\otimes\!\!H\)
                at (0.64,0.575)}];
            \draw[coordinate label = {\(\left[\!\begin{smallmatrix}2\\0\end{smallmatrix}\!\right]\dual\!\!\otimes\!\!H\)
                at (0.78,0.575)}];

            \draw[coordinate label = {\(\left\llbracket\!\begin{smallmatrix}1\\0\end{smallmatrix}\!\right\rrbracket\!\dual\!\!\otimes\!\!H\)
                at (0.44,0.455)}];
            \draw[coordinate label = {\(\left[\!\begin{smallmatrix}0\\1\end{smallmatrix}\!\right]\dual\!\!\otimes\!\!H\)
                at (0.57,0.455)}];
            \draw[coordinate label = {\(\left[\!\begin{smallmatrix}1\\1\end{smallmatrix}\!\right]\dual\!\!\otimes\!\!H\)
                at (0.72,0.455)}];
            \draw[coordinate label = {\(\left[\!\begin{smallmatrix}2\\1\end{smallmatrix}\!\right]\dual\!\!\otimes\!\!H\)
                at (0.85,0.455)}];

            \draw[coordinate label = {\(\left\llbracket\!\begin{smallmatrix}2\\0\end{smallmatrix}\!\right\rrbracket\!\dual\!\!\otimes\!\!H\)
                at (0.36,0.33)}];
            \draw[coordinate label = {\(\left\llbracket\!\begin{smallmatrix}1\\1\end{smallmatrix}\!\right\rrbracket\!\dual\!\!\otimes\!\!H\)
                at (0.5,0.33)}];
            \draw[coordinate label = {\(\left[\!\begin{smallmatrix}0\\2\end{smallmatrix}\!\right]\dual\!\!\otimes\!\!H\)
                at (0.64,0.33)}];
            \draw[coordinate label = {\(\left[\!\begin{smallmatrix}1\\2\end{smallmatrix}\!\right]\dual\!\!\otimes\!\!H\)
                at (0.79,0.33)}];

            \draw[coordinate label = {\(\left\llbracket\!\begin{smallmatrix}3\\0\end{smallmatrix}\!\right\rrbracket\!\dual\!\!\otimes\!\!H\)
                at (0.3,0.21)}];
            \draw[coordinate label = {\(\left\llbracket\!\begin{smallmatrix}2\\1\end{smallmatrix}\!\right\rrbracket\!\dual\!\!\otimes\!\!H\)
                at (0.44,0.21)}];
            \draw[coordinate label = {\(\left\llbracket\!\begin{smallmatrix}1\\2\end{smallmatrix}\!\right\rrbracket\!\dual\!\!\otimes\!\!H\)
                at (0.58,0.21)}];
            \draw[coordinate label = {\(\left[\!\begin{smallmatrix}0\\3\end{smallmatrix}\!\right]\dual\!\!\otimes\!\!H\)
                at (0.72,0.21)}];

            \imagelabelset{coordinate label font = \large}
            \draw[coordinate label = {\rotatebox{-30}{\(\rightsquigarrow\)} at (0.92,0.3)}];
            \draw[coordinate label = {\rotatebox{-60}{\(\rightsquigarrow\)} at (0.76,0.16)}];
            \draw[coordinate label = {\rotatebox{-90}{\(\rightsquigarrow\)} at (0.5,0.1)}];

            \imagelabelset{coordinate label font = \scriptsize}
            \draw[coordinate label = {\(\coh X\) at (0.96,0.28)}];
            \draw[coordinate label = {\(\zeroper\XY\) at (0.79,0.12)}];
            \draw[coordinate label = {\(\text{flop}(\coh \nu_1X)\) at (0.5,0.05)}];
        \end{annotationimage}
    \end{minipage}
\end{figure}
\pagebreak

The above result can be seen as a `fixed-point theorem' for \(\tilt(H)\). Indeed if  \(\LL\in \Pic X\) is an ample line bundle, then the orbit of
\(H\) under iterative applications of \(\LL\,\dual \otimes (-)\) limits to
\(\coh X\) which is fixed as a heart by \(\Pic X\).

We then have the following recipe, which shows every heart cone is non-zero.

\begin{algorithm}[Skyscraper--hunting]
    Suppose \(K\subset \zeroper\XZ[-1,0]\) is a non-algebraic t-structure
    corresponding to a torsion pair \(\zeroper\XZ=T\ast F\).
    \begin{enumerate}[1.]
        \item By \cref{lem:anyKlivesonX} there is a
            unique birational model \(W\in \Bir\XZ\) such that
            \(\flop^{-1}(K)\) is intermediate with respect to
            \(\zeroper\WZ\) and also contains the full subcategory \(\{w\in \coh
            W\;|\; \mathbf{R}\pi_\ast w = 0\}\) (we say \(K\) \emph{lives} on
            \(W\)). Thus replacing \(K\) by \(\flop^{-1}(K)\) if
            necessary, we may assume \(K\) lives on \(X\).

            With this hypothesis satisfied, every skyscraper sheaf \(\OO_p\in
            \zeroper\XZ\cap \coh X\) is either torsion or torsion-free with
            respect to the torsion pair \(T\ast F\)
            (\cref{lem:everyskyscraperistorsion}).

        \item Suppose \(C_i\subset X\) is an exceptional curve and \(\LL_i\in \Pic
            X\) is the line bundle with degree \(1\) on \(C_i\) and degree \(0\)
            elsewhere. In particular, the maximal and minimal elements of
            \(\tilt(H)\) whose heart cone is generated by \(\LL_i\) are the
            categories \(\zeroper\pair{X}{X_i}\) and  \(\antizeroper\pair{X}{X_i}\)
            respectively, where \(X\to X_i\) is the partial contraction of \(C_i\).

            Now if some \(p\in C_i\) satisfies \(\OO_p\in T\), then
            \cref{lem:boundonbothsides} shows \(\LL_i^{\otimes (-n)}\otimes H > K\)
            for arbitrarily large \(n\), and thus \(\zeroper\pair{X}{X_i}\geq
            K\). Likewise if some \(p\in C_i\) satisfies \(\OO_p\in F\), then we have
            the bound \(K\geq \antizeroper\pair{X}{X_i}\).

        \item Thus if there is some curve \(C_i\) with two points \(p,q\in C_i\)
            such that \(\OO_p\in T\) and \(\OO_q\in F\), then \(K\) lies in the
            interval \(\left[\zeroper\pair{X}{X_i},
            \antizeroper\pair{X}{X_i}\right]\) and in particular has non-zero
            heart cone.

        \item Otherwise by connectivity of the exceptional curves in \(X\),
            either \(T\) or \(F\) contains \emph{all} the skyscrapers. By a similar
            argument, this forces one of the bounds \(K\geq \coh X\) or \(K\leq
            \anticoh X\). The bound on the other side can then be found by chasing
            simple tilts from \(H\) or \(H[-1]\) as in \cref{alg:naivetilting},
            and this suffices to prove that \(\CC (K)\) is non-zero
            (\cref{lem:boundedbycoh}).
    \end{enumerate}
\end{algorithm}

\subsection*{Heart cones of semi-geometric hearts}
\fakesubsection{Heart cones of semi-geometric hearts}
Continuing to work with the \(\text{cA}_2\) crepant resolution above, consider a
flop \(\rho: X\dashrightarrow \nu_1 X\), and the partial contraction
\(\tau\colon X\to Y\) of the flopped curve. The heart cone of the associated
semi-geometric heart \(\zeroper\XY\) is given by the ray \(\sigma = \CC(\coh
X)\cap \CC(\flop\coh (\nu_1 X))\).

To examine other hearts that lie on \(\sigma\), note that \(Y\) has a unique
singular point (the image of the flopped curve), and \(\rho\) is simply the
Atiyah flop over a neighbourhood \(Z_1\) of this point. Computing the category of
\(\sigma\)-semistable objects then recovers the category of perverse sheaves
associated to the \(\text{cA}_1\) flopping contraction \(\tau^{-1}Z_1\to Z_1\),
and thus (up to tilting in skyscrapers in the smooth locus) every other intermediate
heart with heart cone \(\sigma\) is obtained from tilts of this smaller category
\(\zeroper\pair{\tau^{-1}Z_1}{Z_1}\).

Convex geometrically this manifests as the fact that the \(\text{cA}_2\)
intersection arrangement `looks like' the \(\text{cA}_1\) arrangement locally around
\(\sigma\) (as is apparent from \cref{fig:AffA1,fig:AffA2}). This phenomenon is
more general and a correspondence between the heart fan of \(H\) and that of its
semistable subcategories is sketched in \cref{rmk:ZoomIntoHFan}.
Broomhead--Pauksztello--Ploog--Woolf's \emph{tangent multifan} construction
\cite[see][`Further work']{broomheadHeartFanAbelian2024} investigates this to a
greater depth.

\subsection*{Notation and Conventions}
\fakesubsection{Notation and Conventions}
We work over the ground field \(\mathbb{C}\), and fix once and for all a
complete local isolated compound du Val singularity \(Z=\Spec R\) with singular
point corresponding to the maximal ideal \(\mathfrak{m}\subset R\).

All algebras we consider are finitely generated \(R\)-algebras. Given
such an algebra \(\Lambda\), we write \(\Db \Lambda\) for the bounded derived
category of right \(\Lambda\)-modules and \(\Dfl \Lambda\) for the full
subcategory of complexes whose cohomology (with respect to \(\rmod\Lambda\)) has
finite length. This forms a triangulated subcategory of \(\Db \Lambda\), and we
write \(\KK \Lambda\) for the Grothendieck group of \(\Dfl \Lambda\).

Likewise, all \(3\)-folds we consider are \(Z\)-schemes. Given such a scheme
\(\pi\colon X\to Z\), we write \(\Db X\) for the bounded derived category of coherent
sheaves on \(X\) and \(\Dm X\) for the full subcategory of \(\Db X\) containing
complexes whose cohomology sheaves are supported on \(\pi^{-1}[\mathfrak{m}]\).
Again, we write \(\KK X\) for the Grothendieck group of \(\Dm X\).

Given full subcategories \(U,V\) of an Abelian (or triangulated) category, we
write \(U\ast V\) for the full subcategory of objects \(x\) that sit in an exact
sequence \(0\to u\to x \to v \to 0\) (resp.\ an exact triangle \(u\to x\to v\to
u[1]\)) with \(u\in U\), \(v\in V\). Note the operation \((\ast)\) is
associative~\cite[lemma 1.3.10]{beilinsonFaisceauxPervers1982} so writing
\(U^{\ast n}=U\ast \cdots \ast U\) (\(n\) factors) is unambiguous, and
we write \(\langle U \rangle = \bigcup_{n\geq 0}U^{\ast n}\) for the
extension-closure of \(U\). If \(U\) sits in a triangulated category and
\(I\subset \bbZ\) is an interval, we write \(U[I]=\langle \{u[i]\;|\; u\in U, \;
i\in I\}\rangle\) and employ obvious shorthands such as
\(U[0,1]=U[[0,1]]\), \(U[\leq 0]=U[(-\infty, 0]]\) when convenient.

\subsection*{Acknowledgements} Many thanks to Michael Wemyss for his relentless
optimism and for patiently watching me triangulate his board week after week;
to Wahei Hara for helpful comments on a previous version of this paper; to Jon
Woolf, David Ploog, David Pauksztello, and Nathan Broomhead for continued
encouragement and eagerness to answer all questions heart fan; to Nick Rekuski
for explaining to me why reversed geometric hearts on curves are Artinian (see
\cref{rmk:artinNoether-geometrichearts}); to Rachael Boyd, Franco Rota, and
Marina Godinho for many helpful conversations; and to Theodoros Papazachariou
and Timothy De Deyn for periodic progress--checks which nudged this manuscript
to completion.

\subsection*{Funding} The author was supported by ERC Consolidator Grant 101001227 (MMiMMa).

\subsection*{Open Access} For the purpose of open access, the author has applied a Creative Commons Attribution (CC:BY) licence to any Author Accepted Manuscript version arising from this submission.

\section{Dynkin combinatorics and intersection arrangements}\label{sec:dynkin}

A recurring motif in our treatment is that of
mutation combinatorics, which allow us to enumerate sets over Dynkin data. We
briefly explain the construction in an abstract setting, and recognise its
various manifestations as they come up throughout the article.

\subsection{Sets with mutation} Let \(G\) be a Dynkin graph with associated Weyl
group \(\WW(G)\), which is generated by simple reflections \(\{s_i\;|\;i\in
G\}\). When \(\WW(G)\) is finite there is a unique longest element \(w_G\) in
the weak (Bruhat) order, and further~\cite[lemma
1.2]{iyamaTitsConeIntersections} there is an involution \(\text{inv}_G:G\to G\)
such that \(w_Gs_iw_G=s_{\text{inv}_G(i)}\). If \(\WW(G)\) is not finite, we
simply declare \(\text{inv}_G\) to be the identity.

For any subgraph \(J\subset G\), this defines a map \(\iota_J: G\setminus
J \to G\) given by \(\iota_J(i)=\text{inv}_{J+i}(i)\) where
\(\text{inv}_{J+i}\) is the involution for the full subgraph \(J\cup
\{i\}\). When the choice of \(J\) is clear, we simply write \(\iota\) instead of
\(\iota_J\).

\begin{definition} 
    The \textit{simple mutation of} \(J\subset G\) at \(i\in
    G\setminus J\) is given by the subgraph \(\nu_i J =   \{i\}\cup J \setminus
    \{\iota(i)\}\) in \(G\).
\end{definition}

Clearly, simple mutation is involutive in the sense that
\(\nu_{\iota(i)}(\nu_i J)=J\). 

When iterating we omit brackets, thus writing
\(\nu_{i_n}\ldots\nu_{i_1}J\) to mean \(\nu_{i_n}(\nu_{i_{n-1}}\ldots\nu_{i_1} J)\),
noting that the sequence makes sense only if \(i_1\in G\setminus J\) and
\(i_{j+1}\in G\setminus (\nu_{i_j}\ldots\nu_{i_1}J)\) for each
\(j=1,\ldots,n-1\). We say such a sequence of symbols
\(\nu=\nu_{i_n}\ldots\nu_{i_1}\) is a \emph{\(J\)-path of length \(n\)}.

\begin{definition}
    A \emph{set with \(G\)-mutation} is a set \(A\) equipped with a map
    \(\mathbb{J}\colon A\to 2^G\) that assigns each element of \(A\) to some subgraph
    of \(G\), and a collection of functions \({(\nu_i \colon \{a\in A\;|\; i \notin
    \mathbb{J}(a)\} \to A)}_{i\in G}\) (called \emph{simple mutations}) such that
    for each \(a\in A\), we have \(\mathbb{J}(\nu_i a)=\nu_i\mathbb{J}(a)\) and
    \(\nu_{\iota(i)}(\nu_i a) = a\).
\end{definition}

If \(A\) is a set with \(G\)-mutation, the \emph{exchange quiver} \(\ExQuiv(A)\)
is a quiver with vertices \(A\) and a labelled arrow \(i : a\to \nu_i a\) for
each \(a\in A\) and \(i\in G\setminus \mathbb{J}(a)\). Paths in the
exchange quiver can be described by a sequence of valid mutations from a
specified starting vertex---it is easy to see that each \(\mathbb{J}(a)\)-path
\(\nu=\nu_{i_n}\ldots\nu_{i_1}\) describes a unique path in the exchange quiver
\begin{equation}
    \label{eqn:positivepath}
    a \xrightarrow{\;i_1\;}
    \nu_{i_1}a \xrightarrow{\;i_2\;} \cdots \xrightarrow{\;i_n\;}
    \nu_{i_{n}}\nu_{i_{n-1}}\ldots \nu_{i_1}a
\end{equation}
which we call the \emph{positive path} \(\nu\) from \(a\), and further every
path in the exchange quiver corresponds to a unique pair \((a,\nu)\) in this
way. If the positive path \eqref{eqn:positivepath} has minimal length among all
paths in \(\ExQuiv(A)\) from \(a\) to \(\nu a\), we say it is \emph{minimal}.

\begin{remark}
    In what follows, we will work with both---an ADE Dynkin graph \(\Delta\) and
    its affine counterpart \(\ul\Delta\). Thus to avoid confusion when
    considering subsets \(J\subset \Delta \subset \ul\Delta\), we use the term
    \emph{spherical \(J\)-paths} for the paths corresponding to the ambient
    graph \(G=\Delta\) and reserve the unqualified term \(J\)-path to mean paths
    for \(G=\ul\Delta\). Writing \(0\in \ul\Delta\setminus\Delta\) for the
    extended vertex, it can be seen that for such \(J\) there is a natural
    bijective correspondence between spherical \(J\)-paths and \(J\)-paths in
    which the symbol \(\nu_0\) does not occur.
\end{remark}

\begin{example}
    Given a Dynkin graph \(G\), the set of all its subgraphs \(2^G\) is
    naturally a set with \(G\)-mutation. This naturally breaks up into smaller
    sets with \(G\)-mutation called \emph{mutation classes}, where the mutation
    class of \(J\subseteq G\) is the subset of \(2^G\) containing subgraphs that
    can be obtained from \(J\) by iterated mutation. Thus \(\{\emptyset\}\) and
    \(\{G\}\) are mutation classes of the trivial subgraphs \(\emptyset\) and
    \(G\) respectively. Iyama--Wemyss \cite[\S 4.2]{iyamaTitsConeIntersections}
    compute all mutation classes containing `large' subgraphs of affine Dynkin
    diagrams, i.e.\ subgraphs whose complements contain \(3\) vertices. We draw
    the mutation class \(\mathcal{E}_{7,4}\) in \cref{figure:E74mutation}.
    \begin{figure}[H]
    \begin{tikzpicture}
        \node (014) at (0,0)   {\affE{o}{*XX*XXX}};
        \node (047) at (3,0)   {\affE{o}{XXX*XX*}};
        \node (037) at (6,1)   {\affE{o}{XX*XXX*}};
        \node (237) at (9,1)   {\affE{X}{X**XXX*}};
        \node (267) at (12,0)  {\affE{X}{X*XXX**}};
        \node (256) at (15,0)  {\affE{X}{X*XX**X}};
        \node (367) at (9,-1)  {\affE{X}{XX*XX**}};
        \node (347) at (6,-1)  {\affE{X}{XX**XX*}};
        \path[->,font=\bfseries\scriptsize,>=angle 90]
            (014) edge [bend left=5] node[above]{1} (047)
            (047) edge [bend left=5] node[below]{7} (014)

            (047) edge [bend left=15] node[above]{4} (037)
            (037) edge [bend right=5] node[below]{3} (047)

            (047) edge [bend right=5] node[above]{0} (347)
            (347) edge [bend left=15] node[below]{3} (047)

            (037) edge [bend left=5] node[above]{0} (237)
            (237) edge [bend left=5] node[below]{2} (037)

            (237) edge [bend left=15] node[above]{3} (267)
            (267) edge [bend right=5] node[below]{6} (237)

            (267) edge [bend left=5] node[above]{7} (256)
            (256) edge [bend left=5] node[below]{5} (267)

            (367) edge [bend right=5] node[above]{3} (267)
            (267) edge [bend left=15] node[below]{2} (367)

            (347) edge [bend left=5] node[above]{4} (367)
            (367) edge [bend left=5] node[below]{6} (347)

            (014) edge [in= 80,out= 100, looseness=15] node[above]{0} (014)
            (014) edge [in=-80,out=-100, looseness=15] node[below]{4} (014)

            (256) edge [in= 80,out= 100, looseness=15] node[above]{2} (256)
            (256) edge [in=-80,out=-100, looseness=15] node[below]{6} (256)

            (037) edge [in= 80,out= 100, looseness=15] node[above]{7} (037)
            (237) edge [in= 80,out= 100, looseness=15] node[above]{7} (237)
            (367) edge [in=-80,out=-100, looseness=15] node[below]{7} (367)
            (347) edge [in=-80,out=-100, looseness=15] node[below]{7} (347);
    \end{tikzpicture}
    \caption{The mutation class of \(J=\{2,3,5,6,7\}\) inside the
         \(\widetilde{\text{E}}_7\) Dynkin graph
         \(\raisebox{5pt}{\dynkin[label,edge length=6pt,ordering=Dynkin, root
         radius=1.5pt] E[1]7}\), where subgraphs are indicated by marking off
         their vertices with a cross (\({\dynkin[root radius=2pt] A{X}}\)). Arrows
         indicate simple mutation, the symbol \(\nu\) has been omitted for
         brevity.}
         \label{figure:E74mutation}
    \end{figure}
\end{example}

\begin{example}
    \label{ex:trivintersectionarrangement}
    If \(A\) is a set with an action of the Weyl group \(\WW(G)\), then there is
    a natural \(G\)-mutation structure given by \(\mathbb{J}(a)=\emptyset\),
    \(\nu_i(a)=a\cdot s_i\) for all \(a\in A\), \(i\in G\). Thus for instance
    the set of chambers \(\Chambers(G)\) in the Tits cone of \(\WW(G)\) is a set
    with \(G\)-mutation, where the (right)-action is by simple wall crossings.
\end{example}

\subsection{Intersection arrangements}\label{subsec:intersectionarrangements}  
The combinatorial construction central to our exposition is a generalisation of
example~\ref{ex:trivintersectionarrangement}, where a subset \(J\subset G\)
determines a hyperplane arrangement in Euclidean space~\cite[chapter
1]{iyamaTitsConeIntersections} and \(G\)-mutation is realised as an action of
the associated Deligne groupoid (also called the \emph{\(J\)-cone groupoid}, see
\S\,2.3.1 \emph{ibid.}). The construction is most relevant to us when \(G\) is
an ADE Dynkin graph \(\Delta\), or its affine counterpart
\(\ul\Delta\) which has extended vertex \(0\in \ul\Delta\setminus \Delta\), and
we recall necessary details below.

\paragraph{Affine root systems}
We recap the unrestricted (Coxeter) setting first, following e.g.\ \cite[\S
1]{kacInfiniteRootSystems1980}. 

Given the affine Dynkin diagram \(\ul\Delta\) of rank \(n\), the \emph{root lattice}
\(\mathfrak{h}=\mathfrak{h}(\ul\Delta)\) is a free \(\mathbb{Z}\)-module with
basis given by the \emph{simple roots} \(\{\alpha_i\;|\; i\in\ul\Delta\}\). 
The Cartan matrix associated to \(\ul\Delta\) gives a degenerate symmetric bilinear
form \((-,-)\) on \(\mathfrak{h}\), and we use these to define the \emph{simple
reflections} \(s_i\colon\mathfrak{h}\to\mathfrak{h}\) (\(i\in \ul\Delta)\) given by
\begin{equation}
    \label{eqn:simplereflectiononroots}
    s_i(\alpha_j) = \alpha_j - (\alpha_i,\alpha_j)\alpha_i.
\end{equation}
These define the \emph{Weyl group} \(\WW(\ul\Delta)=\langle s_i \;|\; i\in
\ul\Delta\rangle\subset \text{GL}(\mathfrak{h})\), which is isomorphic to the
Coxeter group associated to \(\ul\Delta\) with standard generators
\(\{s_i\;|\;i\in \ul\Delta\}\). In fact each subset \(J\subseteq \ul\Delta\) determines a parabolic subgroup \(\WW(J)=\langle s_i\;|\; i\in J\rangle \subseteq
\WW(\ul\Delta)\), which is the Coxeter group associated to the full subgraph
spanned by \(J\). If \(J\neq \ul\Delta\) then this is a finite Coxeter group,
and hence has a unique Coxeter element \(w_J\) (defined to be the longest
element in the Bruhat order).

The \emph{set of real roots} \(\Root(\ul\Delta)\) is the union of
\(\WW(\ul\Delta)\)-orbits of the simple roots. We say a real root is
\emph{positive} if it can be expressed as a non-negative linear combination of
the simple roots, and write \(\Root^+(\ul\Delta)\) for the set of positive real
roots. The set \(\Root^-(\ul\Delta)\) of \emph{negative real roots} is defined
likewise. 

Vertices in the diagram \(\ul\Delta\) can be Lie-theoretically assigned
numerical labels \((\delta_i)_{i\in\ul\Delta}\) which are computed and given
for each Dynkin type in~\cite[table Z]{kacInfiniteRootSystems1980}. This is
always a tuple of positive integers, and the `extended' vertex
\(0\in\ul\Delta\setminus\Delta\) is assigned the integer \(\delta_0=1\). Then
the vector \(\delta=\sum_{i\in\ul\Delta}\delta_i\alpha_i\) is called the
\emph{primitive positive imaginary root} in \(\mathfrak{h}\), and any non-zero
multiple of \(\delta\) is called an \emph{imaginary root}. The \emph{root
system} of \(\ul\Delta\) is then the set of all (real and imaginary) roots in
\(\mathfrak{h}\), and is preserved by the action of the Weyl group which acts
transitively on \(\Root(\ul\Delta)\) and fixes \(\delta\).

\paragraph{Restricted root systems} Suppose \(\JJ\subset\ul\Delta\) is such that
\(|\ul\Delta\setminus\JJ|\geq 2\). This defines a decomposition
\[\mathfrak{h}(\ul\Delta) = \mathfrak{h}(\JJ) \oplus
\mathfrak{h}(\ul\Delta\setminus\JJ),\] where we write \(\mathfrak{h}(\JJ)\) for
the span of the simple roots \(\{\alpha_i\;|\;i\in\JJ\}\) (and define
\(\mathfrak{h}(\ul\Delta\setminus\JJ)\) likewise). We say
\(\mathfrak{h}(\ul\Delta\setminus \JJ)\) is the \emph{restricted root lattice}
associated to \(\JJ\), and the \emph{restricted root system} is the image of the
root system under the surjection \(\mathfrak{h}(\ul\Delta) \onto
\mathfrak{h}(\ul\Delta\setminus\JJ)\).

The notions of restricted (simple,
positive, negative) real roots and restricted imaginary roots are defined
likewise as the non-zero images of the corresponding objects in
\(\mathfrak{h}\). In particular there are \(|\ul\Delta\setminus\JJ|\) restricted
simple roots \(\{\alpha_i\;|\;i\in\ul\Delta\setminus\JJ\}\), and a
\emph{primitive positive restricted imaginary root} \(\delta_\JJ =
\sum_{i\in\ul\Delta\setminus\JJ} \delta_i\alpha_i\). Writing
\(\Root(\ul\Delta,\JJ)\) for the set of restricted real roots (and defining the
subsets \(\Root^\pm(\ul\Delta,\JJ)\) of positive and negative restricted real
roots accordingly), we thus see that the restricted root system is given as
\[
    \underbrace{\Root^+(\ul\Delta,\JJ) \;\sqcup\;
    \Root^-(\ul\Delta,\JJ)}_{\Root(\ul\Delta,\JJ)} \;\sqcup\;\;
    \{n\delta_\JJ\;|\;n\in\mathbb{Z}\setminus \{0\}\}.
\]

To obtain the analog of a Weyl-group action for restricted roots, it is
necessary to consider simultaneously all restricted root lattices
\(\mathfrak{h}(\ul\Delta\setminus \nu\JJ)\) ranging over \(\JJ\)-paths \(\nu\).
In the Coxeter setting, a simple reflection \(s_i\in\WW(\ul\Delta)\) was
characterised by the property of being an involution such that
\(s_i(\alpha_j)-\alpha_j\) is equal to \(-2\alpha_i\) if \(j=i\), and is a
non-negative multiple of \(\alpha_i\) otherwise.  Iyama--Wemyss suggest the
following generalisation.

\begin{proposition}
    [{\cite[lemmas 5.1 and 5.2]{nabijouGVGWInvariants2023}}]
    \label{prop:wallcrossingdefn}
    For \(\JJ\subset \ul\Delta\) and \(i\in\ul\Delta\setminus\JJ\), the linear
    map \(\mathfrak{h}\to\mathfrak{h}\) given by the action of \(w_\JJ
    w_{\JJ+i}\in\WW(\ul\Delta)\) maps the subset \(\mathfrak{h}(\nu_i
    \JJ)\) isomorphically onto \(\mathfrak{h}(\JJ)\), thus inducing an isomorphism
    \begin{equation}
        \label{eqn:wallcrossingdefn}
        \begin{tikzcd}[column sep=huge]
            \mathfrak{h} \arrow[d, twoheadrightarrow]
                \arrow[r, "w_{\JJ}w_{\JJ+i}"]
            & \mathfrak{h} \arrow[d, twoheadrightarrow]\\
            \mathfrak{h}(\ul\Delta\setminus\nu_i\JJ)
                \arrow[r, "\varphi_i", "\sim"']
                         &
            \mathfrak{h}(\ul\Delta\setminus\JJ)
        \end{tikzcd}
    \end{equation}
    This map preserves the root systems, in particular acting on the simple
    roots and the primitive positive imaginary root as
    \[
        \varphi_i(\delta_{\nu_i\JJ})=\delta_\JJ,
        \qquad \varphi_i (\alpha_{\iota(i)})=-\alpha_i,
        \qquad \varphi_i(\alpha_j)-\alpha_j \;\in\; \mathbb{Z}_{\geq 0}\cdot
        \alpha_i \quad (j\neq i).
    \]
    Further, defining
    \(\varphi_{\iota(i)}:\mathfrak{h}(\ul\Delta\setminus\JJ)\to
    \mathfrak{h}(\ul\Delta\setminus\nu_i \JJ)\) analogously, we have
    \( \varphi_{\iota(i)} = (\varphi_i)^{-1}\).
\end{proposition}

\paragraph{Hyperplane arrangements} Let \(\mathfrak{h}\) be the root lattice associated
to \(\ul\Delta\) as above, with the action of
\(\WW(\ul\Delta)\) as described. The dual action on the Euclidean space
\(\Theta=\Hom_\mathbb{Z}(\mathfrak{h},\mathbb{R})\) preserves each of the
subsets \(\{\delta>0\}\), \(\{\delta<0\}\), and \(\{\delta<0\}\), and the
fundamental domains for the respective actions have closures given by the cones
\begin{gather*}
    \Chamb^+=\{ \theta\in \Theta \;|\; \theta(\alpha_i)\geq 0 \text{ for all
    }i\in\ul\Delta\}, \qquad \Chamb^-= -\Chamb^+,\\
    \Chamb^0=\{\theta\in\Theta \;|\; \theta(\alpha_i)\geq 0 \text{ for all
    }i\in\Delta, \quad \theta(\delta)=0\}.
\end{gather*}
Further the faces of \(w\Chamb^+\) are all given by \(w\Chamb_J^+\) for some
subset \(J\subset\ul\Delta\), where \[\Chamb_J^+=\Chamb^+\cap \bigcap_{i\in
J}\{\alpha_i=0\}.\] The faces \(w\Chamb^-_J \subseteq w\Chamb^-\)
(\(J\subseteq\ul\Delta\)) and \(w\Chamb^0_J \subseteq w\Chamb^0\)
(\(J\subseteq\Delta\)) are given likewise, and as a consequence we have a
complete simplicial fan in \(\Theta\) called the \emph{Weyl arrangement} given
as
\begin{equation}
    \label{eqn:titsdecomposition}
    \Arr(\ul\Delta) =
    \underbrace{
        \left\{  
            wC_J^+ \;\middle\vert\; w\in \WW(\ul\Delta),\quad J\subseteq\ul\Delta
        \right\}
        }_{\Arr^+(\ul\Delta)} \;\;\cup\;
    \underbrace{
        \left\{  
            wC_J^0 \;\middle\vert\; w\in \WW(\Delta),\quad J\subseteq\Delta
        \right\}
        }_{\Arr(\Delta)} \;\;\cup\;
    \underbrace{
        \left\{  
            wC_J^-  \;\middle\vert\; w\in \WW(\ul\Delta),\quad J\subseteq\ul\Delta
        \right\}
        }_{\Arr^-(\ul\Delta)}.
\end{equation}
The Weyl arrangement is induced by the root system for \(\ul\Delta\), in the
sense that every face can be described as the intersection of half-spaces
defined by root hyperplanes \(\{\alpha=0\}\). Further each such hyperplane
is corresponds to a unique positive root in \(\Root^+(\ul\Delta)\cup\{\delta\}\).

The subfan \(\Arr^+(\ul\Delta)\) is supported on the rational subset
\(\{\delta>0\}\cup\{0\}\), and is called the \emph{Tits cone}. The subfan
\(\Arr(\Delta)\) is supported on \(\{\delta=0\}\), and induces a complete
simplicial fan on this hyperplane. Further it has finitely many cones which are
indexed over the parabolic subgroup \(\WW(\Delta)\). This is possible because
\(\WW(\ul\Delta)\) decomposes as a semi-direct product of \(\WW(\Delta)\) and
the coroot lattice, and the action on \(\{\delta=0\}\) is such that
\(\WW(\Delta)\) acts faithfully while the coroot lattice fixes the hyperplane
pointwise.

Now a subset \(\JJ\subset\ul\Delta\) with \(|\ul\Delta\setminus\JJ|\geq 2\)
defines the restricted root lattice \(\mathfrak{h}(\ul\Delta\setminus\JJ)\),
which we view as a split quotient of \(\mathfrak{h}\). Accordingly, the dual
space \(\Theta(\ul\Delta\setminus\JJ) =
\Hom_\mathbb{Z}(\mathfrak{h}(\ul\Delta\setminus\JJ),\mathbb{R})\) can be
identified with the subspace \(\bigcap_{i\in \JJ}\{\alpha_i=0\}\subset \Theta\).
Since the subspace is defined by root hyperplanes, every cone \(\sigma\in
\Arr(\ul\Delta)\) intersects \(\Theta(\ul\Delta\setminus\JJ)\) in a face of
\(\sigma\), whence the Weyl arrangement induces a complete simplicial fan
\begin{align*}
    \Arr(\ul\Delta,\JJ)
    &= \underbrace{
        \left\{
            \sigma\cap \Theta(\ul\Delta\setminus \JJ) \;\middle\vert\;
            \sigma\in\Arr^+(\ul\Delta)
        \right\}}_{\Arr^+(\ul\Delta,\JJ)} \;\cup\;
    \underbrace{
        \left\{
            \sigma\cap \Theta(\ul\Delta\setminus \JJ) \;\middle\vert\;
            \sigma\in\Arr(\Delta)
        \right\}}_{\Arr(\Delta,\JJ)} \;\cup\;
    \underbrace{
        \left\{
            \sigma\cap \Theta(\ul\Delta\setminus \JJ) \;\middle\vert\;
            \sigma\in\Arr^-(\ul\Delta)
        \right\}}_{\Arr^-(\ul\Delta,\JJ)}
\end{align*}
called the (full) \emph{\(\JJ\)-cone arrangement}. This is precisely the subfan of
\(\Arr(\ul\Delta)\) supported on \(\Theta(\ul\Delta\setminus\JJ)\). It can be
seen that the \(\JJ\)-cone arrangement is induced by the restricted root system
associated to \(\JJ\subset \ul\Delta\).

\begin{figure}[h]
    \begin{center}
    \begin{annotationimage}{width=0.8\textwidth}{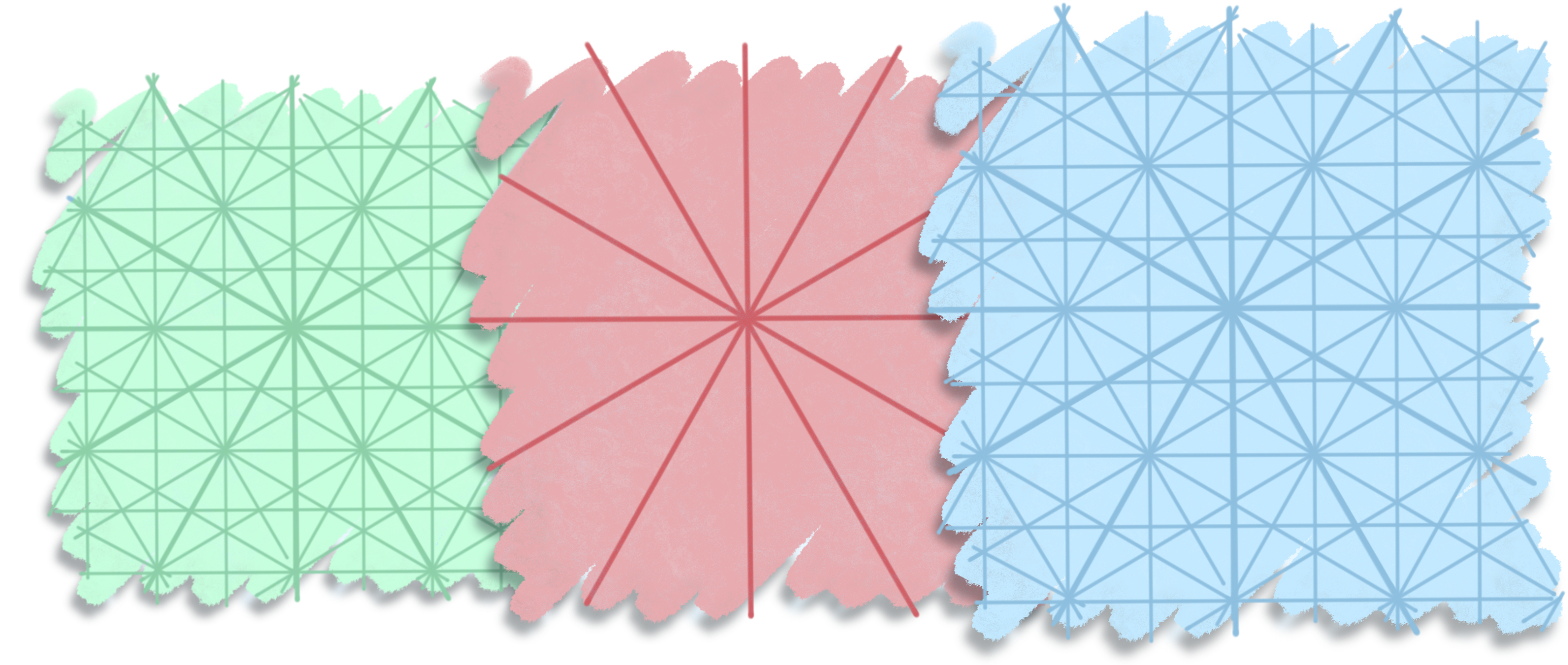}
        \imagelabelset{coordinate label font=\small\bfseries}
        \draw[coordinate label = {{+} at (0.825,0.557)}];
        \draw[coordinate label = {{0} at (0.55,0.557)}];
        \draw[coordinate label = {{--} at (0.157,0.485)}];
        \draw[coordinate label = {{\color{gray} \scriptsize\(\{\delta_J=0\}\)}
        at (0.48,0.02)}];
        \draw[coordinate label = {{\color{gray} \scriptsize\(\{\delta_J+1=0\}\)}
        at (0.18,0.02)}];
        \draw[coordinate label = {{\color{gray} \scriptsize\(\{\delta_J-1=0\}\)}
        at (0.8,0.02)}];
    \end{annotationimage}
    \end{center}
    \caption{Representative affine slices of the \(\mathcal{E}_{7,4}\)
        arrangement, associated to the Dynkin data \(\affE{o}{*XX*XXX}\)
        (\(J=\{2,3,5,6,7\}\)) as in \cref{figure:E74mutation}. The principal
        chambers \(\Chamb_J^+,\Chamb_J^0,\Chamb_J^-\) are indicated by the symbols
        \textbf{+,0,--} respectively.}
    \label{fig:E74Arr}
\end{figure}

The subfan \(\Arr^+(\ul\Delta,\JJ)\), called the \emph{Tits cone}, contains
infinitely many cones and is supported on the subset \(\{\delta_\JJ > 0\}\cup
\{0\}\). The subfan \(\Arr(\Delta,\JJ)\) on the other hand is a complete
simplicial fan on the hyperplane \(\{\delta_\JJ=0\}\), with finitely many cones
whose combinatorics are controlled by the ADE Dynkin diagram \(\Delta\). We
write \(\Chambers(\ul\Delta,\JJ)\) for the subset of maximal cones in
\(\Arr^+(\ul\Delta,\JJ)\), and \(\Chambers(\Delta,\JJ)\) for the subset of
maximal cones in \(\Arr(\Delta,\JJ)\).

\paragraph{Simple wall crossings} Any chamber \(wC^+\in\Chambers(\ul\Delta)\)
has \(|\ul\Delta|\) codimension-one faces (\emph{walls}), where the \(i\)th
wall \(w\Chamb^+_{\{i\}}\) spans the root hyperplane \(\{w\alpha_i=0\}\). The
wall \(w\Chamb^+_{\{i\}}\) is a face of precisely one other chamber, which
lies on the `other side' of this hyperplane. Examining
\eqref{eqn:simplereflectiononroots} shows that this chamber is given by
\(ws_i\Chamb^+\), and we say \(ws_i\Chamb^+\) is obtained from \(w\Chamb^+\) by
a \emph{simple wall crossing}. Iyama--Wemyss~\cite[\S
1]{iyamaTitsConeIntersections} show that while there may not be an underlying
group action on the chambers of the \(\JJ\)-cone arrangement, one can navigate
between them by analogous wall-crossing combinatorics. 

To explain the construction, we first note that if \(w\Chamb_J^+ =
w'\,\Chamb_{J'}^+\), then we must have \(J=J'\) and \(w \WW(J)=w'\,\WW(J')\) as
cosets.  Thus for each \(\sigma\in\Chambers(\ul\Delta,\JJ)\) there is a unique
subset \(\mathbb{J}(\sigma)\subset\ul\Delta\) such that \(\sigma\) can be
expressed as \(\sigma=w\Chamb^+_{\mathbb{J}(\sigma)}\). Further the
codimension-one faces of \(\sigma\) are in bijection with
\(\ul\Delta\setminus\mathbb{J}(\sigma)\), where \(i\in
\ul\Delta\setminus\mathbb{J}(\sigma)\) determines the face
\(\sigma_i=w\Chamb_{\mathbb{J}(\sigma)+i}\). We then have the following.

\begin{lemma}[{\cite[lemma 1.23]{iyamaTitsConeIntersections}}]
    \label{lem:simplewallcrossingfromcone}
    Given \(\sigma\in\Chambers(\ul\Delta,\JJ)\) and
    \(i\in\ul\Delta\setminus\mathbb{J}(\sigma)\), there is a unique cone
    \(\nu_i\sigma\in \Chambers(\ul\Delta,\JJ)\) which satisfies
    \(\nu_i\sigma\cap \sigma = \sigma_i\), which we call the \emph{simple wall
    crossing of \(\sigma\) at \(i\)}. Explicitly, if
    \(\sigma=w\Chamb^+_J\) then this cone is given by
    \[
        \nu_i \cdot w\Chamb_J^+ = w\cdot w_J w_{J+i} \cdot \Chamb^+_{\nu_iJ}.
    \]
    This gives \(\Chambers(\ul\Delta,\JJ)\) the structure of a set with
    \(\ul\Delta\)-mutation, i.e.\ we have
    \(\mathbb{J}(\nu_i\sigma)=\nu_i\mathbb{J}(\sigma)\) and
    \(\nu_{\iota(i)}\nu_i\sigma = \sigma\). Moreover, this set has a connected
    exchange quiver (i.e.\ any two chambers in \(\Chambers(\ul\Delta,\JJ)\) are
    connected by a finite sequence of simple wall crossings.)

    Analogously, \(\Chambers(\Delta,\JJ)\) equipped with the data of simple wall
    crossings is a set with \(\Delta\)-mutation, and this has a finite and
    connected exchange quiver.
\end{lemma}

\begin{figure}[H]
    \begin{minipage}[c]{0.4\textwidth}
        \caption{Continuing \cref{fig:E74Arr}, some wall crossings in the
        \(\mathcal{E}_{7,4}\) Tits cone are indicated. Within
        each chamber \(\sigma\), the walls are labelled by indices \(i\in
        \ul\Delta\setminus \mathbb{J}(\sigma)\) such that crossing the
        \(i^\text{th}\) wall lands in the adjacent chamber \(\nu_i \sigma\).}
        \[\ul\Delta = \dynkin[label,edge length=6pt,ordering=Dynkin, root radius=1.5pt] E[1]7\]
    \end{minipage}\qquad
    \begin{minipage}[c]{0.4\textwidth}
        \begin{annotationimage}{width=\textwidth}{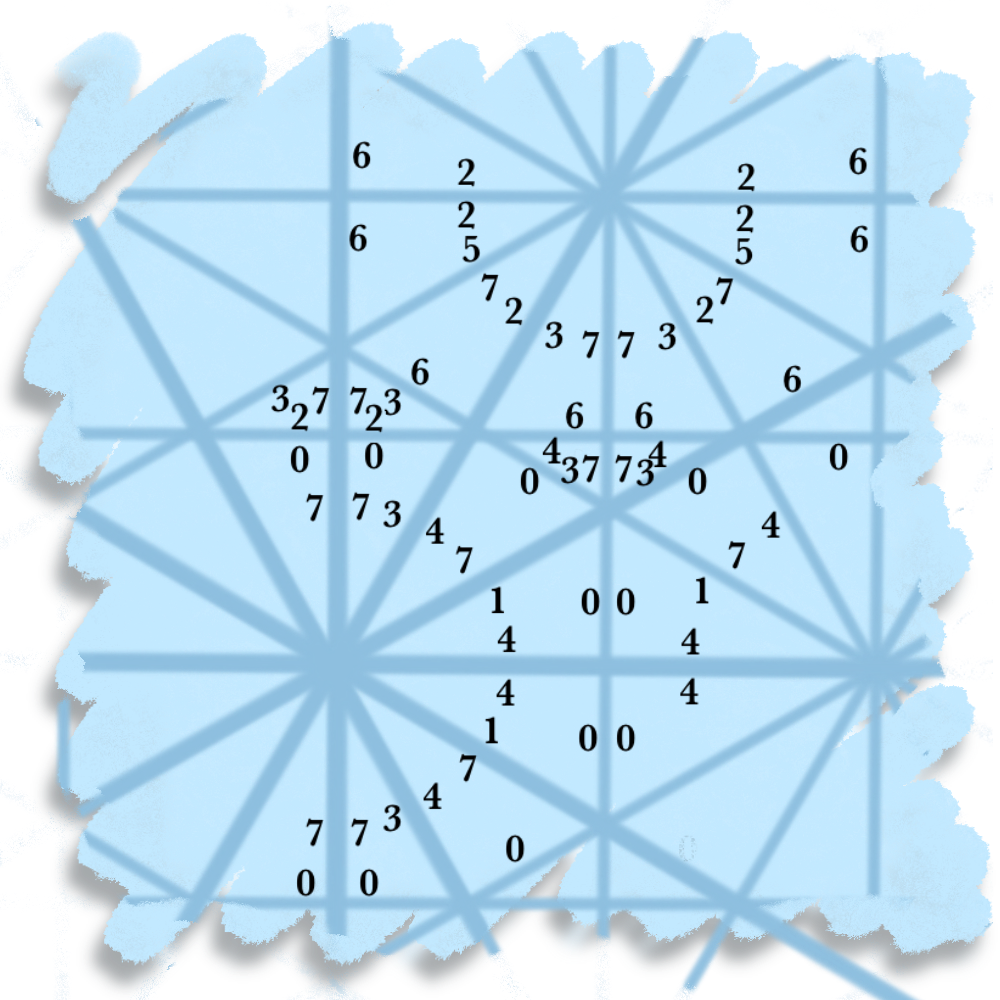}
            \imagelabelset{coordinate label font=\small\bfseries}
            \draw[coordinate label={{+} at (0.55,0.4)}];
            \draw[annotation right = {\(\nu_4\nu_1\nu_0\Chamb_J^+\) at 0.5}]
                to (0.84,0.5);
            \draw[annotation right = {\(\nu_4\nu_1\nu_0\Chamb_J^+\) at 0.5}]
                to (0.84,0.5);
            \draw[annotation right = {\(\nu_4\nu_0\nu_1\Chamb_J^+\) at
            0.62}] to (0.56,0.62);
        \end{annotationimage}
    \end{minipage}
    \label{fig:E74wallcrossings}
\end{figure}

If we write \(\Chamb_\JJ^0=\Chamb^0\cap
\Theta(\ul\Delta\setminus\JJ)\) (even when \(\JJ\) does not lie in \(\Delta\)),
we see that every \(\sigma\in\Chambers(\Delta,\JJ)\) can be written as
\(\nu\Chamb^0_\JJ\) for some spherical
\(\mathbb{J}\left(\Chamb^0_\JJ\right)\)-path \(\nu\). 
It follows that if \(\JJ\subset\Delta\), then the \(\JJ\)-cone arrangement
\(\Arr(\ul\Delta,\JJ)\) can be explicitly described as
\begin{equation}
    \label{eqn:intersectiondecomposition}
    \underbrace{
        \bigcup_{\substack{\nu \text{ an \(\JJ\)-path}\\\phantom{\JJ}}}
        \faces\left(\nu C_\JJ^+\right)}_{\Arr^+(\ul\Delta,\JJ)} \quad\cup\;
    \underbrace{\bigcup_{\substack{\nu \text{ a
        spherical}\\\text{\(\JJ\)-path}}}
        \faces\left(\nu C_\JJ^0\right)}_{\Arr(\Delta,\JJ)} \quad\cup\;
    \underbrace{\bigcup_{\substack{\nu \text{ an \(\JJ\)-path}\\\phantom{\JJ}}}
        \faces\left(\nu C_\JJ^-\right)}_{\Arr^-(\ul\Delta,\JJ)}.
\end{equation}

\paragraph{The weak order on chambers} Fixing the base chamber \(\Chamb_\JJ^+\)
in the \(\JJ\)-cone arrangement endows the set of chambers
\(\Chambers(\ul\Delta,\JJ)\) with a partial order, called the \emph{weak order},
where two chambers \(\sigma,\sigma'\) satisfy \(\sigma\geq \sigma'\) if and only
if every hyperplane separating \(\sigma\) from \(\Chamb_\JJ^+\) also separates
\(\sigma'\) from \(\Chamb_\JJ^+\). Evidently the poset then has a unique
greatest element, \(\Chamb_\JJ^+\).

The following lemma summarises key properties
of this partial order, following \cite{iyamaTitsConeIntersections,salvettiTopologyComplementReal1987}.

\begin{lemma}
    \label{lem:atomicproperties}
    Fix a chamber \(\sigma\in \Chambers(\ul\Delta,\JJ)\) and a positive
    \(\mathbb{J}(\sigma)\)-path \(\nu=\nu_{i_n},...,\nu_{i_1}\).
    \begin{enumerate}[(1)]
        \item \label{lem:coveringrelationsonchambers}
            If \(\nu\) has length \(1\), then one of \(\sigma <
            \nu\sigma\) or \(\sigma > \nu\sigma\) must hold. Moreover, the
            relation is covering, i.e.\ there is no chamber \(\sigma'\)
            satisfying \(\sigma>\sigma'>\nu\sigma\) (or
            \(\nu\sigma>\sigma'>\sigma\)).
        \item If \(\sigma\geq \nu\sigma\), then the path \(\nu\) is
            \emph{minimal} (i.e.\ has minimal length among all positive paths
            from \(\sigma\) to \(\nu\sigma\)) if and only if it is
            \emph{reduced} (i.e.\ crosses each hyperplane at most once), if and
            only if it is \emph{atomic} (i.e.\ successive truncations of \(\nu\)
            give a monotone sequence \(\sigma > \nu_{i_1}\sigma
            >\nu_{i_2}\nu_{i_1} ... > \nu \sigma\)).
        \item If \(\sigma\geq \nu\sigma\), then the interval
            \([\nu\sigma,\sigma]\) contains finitely many elements.
    \end{enumerate}
\end{lemma}

Now any two chambers in a locally finite hyperplane arrangement can be connected
by a reduced path. It follows that every relation in the poset
\(\Chambers(\ul\Delta,\II)\) is realised by the Hasse quiver, i.e.\ we have
\(\sigma < \sigma'\) if and only if there is a chain of covering relations
\(\sigma=\sigma^0\lessdot...\lessdot\sigma^n=\sigma'\). This also shows that
every covering relation in the poset must arise from simple wall crossing, so
that the Hasse quiver is a subquiver of the exchange quiver in the obvious way
and a positive path is atomic if and only if it lies in the Hasse quiver.

Analogous statements hold for \(\Chambers(\Delta,\JJ)\), with maximal
chamber \(\Chamb_\JJ^0\).

\section{Algebra, geometry, and everything in between}\label{sec:functors}

As in the introduction, let \(Z=\Spec R\) be a complete local isolated
compound du Val (cDV) singularity, equivalently the spectrum of a
three-dimensional complete local Gorenstein \(\mathbb{C}\)-algebra \(R\) such
that the unique singularity at the maximal ideal \(\mathfrak{m}\) is at worst
terminal. In particular, \(R\) is a normal domain.

A \emph{flopping contraction} over \(Z\) \cite[definition
6.10]{kollarBirationalGeometryAlgebraic1998} is a projective birational
morphism \(\pi\colon X\to Z\) such that \(X\) is normal, the exceptional locus is
codimension \(2\), and the canonical divisor on \(X\) is numerically \(\pi\)-trivial (equivalently
\(\pi^\ast\omega_Z=\omega_X\), i.e.\ \(\pi\) is \emph{crepant}). 

Fixing such a flopping contraction \(\pi\colon X\to Z\) we look at the
triangulated category \(\Dm X\) containing complexes of coherent sheaves
supported within \(\pi^{-1}[\mathfrak{m}]\), and seek to examine \emph{hearts of (bounded)
t-structures} on \(\Dm X\)---an example being \(\coh X =\Coh X \cap \Dm X\).
More generally any Abelian category \(H\) is naturally the heart of a
t-structure inside its bounded derived category \(\Db H\). 

Thus in this section
we enumerate derived equivalences of \(X\) with a multitude of algebras,
varieties, and more generally, non-commutative schemes, thereby enumerating
hearts in \(\Dm X\).

\subsection{Van den Bergh's equivalence} 
\label{subsec:VdBequiv}
By Kawamata's
vanishing theorem, the map \(\pi\) is acyclic i.e.\ \(\mathbf{R}\pi_\ast
\OO_X=\OO_Z\) and thus the exceptional fiber \(\ul{C}=\pi^{-1}[\mathfrak{m}]\)
is such that the reduced subscheme \(C=\ul{C}_\text{red}\) is a finite
collection of rational curves~\cite[lemma
3.4.1]{vandenberghThreedimensionalFlopsNoncommutative2004}.
In this context, for each integral exceptional curve \(C_i\subset C\) Van den
Bergh~\cite[\S 3.5]{vandenberghThreedimensionalFlopsNoncommutative2004} gives
the construction of a distinguished indecomposable vector bundle \(\mathscr{N}_i\) on \(X\),
such that the bundle \(\VV\XZ=\OO_X\oplus \bigoplus_i\mathscr{N}_i\) is tilting.
Consequently \(X\) is derived equivalent to the basic \(R\)-algebra \(\Lambda =
\End_X\VV\XZ\) via the functor
\begin{equation}
    \label{eqn:zeroperforpi}
    \textnormal{VdB}\colon
    \Db\Lambda\xrightarrow{\quad(-)\LTensor_\Lambda \VV\XZ\quad} \Db X.
\end{equation}
This equivalence maps the natural heart \(\rmod \Lambda\subset \Db \Lambda\) to
the \textit{category of perverse sheaves}
\(\zeroPer\XZ=\text{VdB}(\rmod\Lambda)\), which
should be thought of as the `standard' heart in \(\Db X\) owing to its amiable
cohomological properties. By~\cite[corollary
3.2.8]{vandenberghThreedimensionalFlopsNoncommutative2004},
\begin{equation}
    \label{eqn:zeroperdefn}
    \zeroPer\XZ = \left\{x\in \Coh X[0,1] \;\middle\vert\;
        \begin{array}{l}
            \mathbf{R}^1\pi_\ast (\HH^0x)=0,
            \quad \pi_\ast (\HH^{-1}x)=0,\\
            \Hom(c,\HH^{-1}x)=0 \text{ for all }c\in\Coh X
            \\ \text{satisfying }\mathbf{R}\pi_\ast c=0. 
        \end{array}
    \right\}.
\end{equation}

The equivalence~\eqref{eqn:zeroperforpi}, being \(R\)-linear, identifies the
subcategory \(\Dm X\subset \Db X\) of complexes with cohomology supported within
\(C\) with the subcategory \(\Dfl \Lambda\subset \Db \Lambda\) of complexes
whose cohomology modules have finite length over \(\Lambda\) (equivalently, the
thick subcategory generated by simple \(\Lambda\)-module.)

\begin{definition}
    Write \(\zeroper\XZ\) for the subcategory \(\zeroPer\XZ\cap \Dm X\) of \(\Db
    X\), which under the equivalence~\eqref{eqn:zeroperforpi} corresponds to the
    natural heart \(\rflmod\Lambda\subset \Dfl\Lambda\) of finite length
    \(\Lambda\)-modules. 
\end{definition}

\paragraph{Dynkin data for flopping contractions} Returning to the contraction
\(\pi\colon X\to Z\), the integral components of the reduced exceptional fiber \(C\) can be
naturally indexed over a subgraph of a Dynkin diagram by considering a generic
hyperplane section (\emph{general elephant}) \(\overline{Z}=\Spec
\overline{R}\hookrightarrow \Spec R\) and the pullback
\(\overline{X}=X\times_Z\overline{Z}\). Indeed by assumptions on \(R\), the
variety \(\overline{Z}\) is the germ of a canonical surface singularity and the
map \(\overline{X}\to \overline{Z}\) is a crepant partial resolution which is
independent of the choice of general elephant~\cite[theorem
1.14]{reidMinimalModelsCanonical}. In particular \(\overline{X}\) is obtained by
contracting a subset of exceptional curves in the minimal resolution of
\(\overline{Z}\). By the McKay
correspondence~\cite{mckayGraphsSingularitiesFinite1980}, the exceptional curves
in the minimal resolution are indexed over some ADE Dynkin graph \(\Delta\).
Recording the indices of the contracted curves in a subset \(\JJ\subset\Delta\)
shows that the integral components of \(C\) are naturally indexed over the
Dynkin subgraph \(\Delta\setminus\JJ\) and we can write the reduced exceptional
fiber as \(C=\bigcup_{i\in \Delta\setminus\JJ}C_i\) where each component \(C_i\)
is a \(\mathbb{P}^1\).

The indexing naturally carries over to indecomposable summands of Van den Bergh's bundle
\(\VV\XZ\) in the following fashion.

\begin{lemma}
    \label{lem:indexing-of-mmg-summands}
    The bundle \(\VV\XZ\) has precisely one indecomposable summand isomorphic to
    \(\OO_X\). Further, for each irreducible exceptional curve \(C_i\subset
    C\) (\(i\in \Delta\setminus\JJ\)) in \(X\), the bundle \(\VV\XZ\) has
    an indecomposable summand \(\mathscr{N}_i\) such that
    the closed subscheme \(\textnormal{c}_1(\mathscr{N}_i\dual)\subset X\)
    transversely intersects \(C_i\) once, and is disjoint from \(C_j\) for \(j\neq
    i\). 

    This assignment \(C_i\mapsto \mathscr{N}_i\) is a bijection between integral
    \(\pi\)-exceptional curves in \(X\) and non-free indecomposable summands of
    \(\VV\XZ\).
    \begin{proof}
        This is implicit in the construction of \(\VV\XZ\), see
        \cite[\S\S\,3.4 and
        3.5]{vandenberghThreedimensionalFlopsNoncommutative2004}.
    \end{proof}
\end{lemma}

It is convenient to view \(\Delta\) as sitting inside the associated extended
(affine) Dynkin diagram \(\underline{\Delta}\) with extended vertex \(0\in
\underline\Delta\setminus \Delta\). 

We reserve the index \(0\) for the trivial summand of \(\VV\XZ\), writing
\(\mathscr{N}_0=\OO_X\). Thus indecomposable summands \(\mathscr{N}_i\subset \VV\XZ\), and hence also finitely generated indecomposable projective
\(\Lambda\)-modules \(P_i=\Hom_X\left(\VV\XZ,\mathscr{N}_i\right)\), are
naturally indexed over \(\ul\Delta\setminus \JJ\). The indexing carries
over to simple \(\Lambda\) modules via the standard simple-projective duality.
By~\cite[proposition 3.5.7]{vandenberghThreedimensionalFlopsNoncommutative2004},
the simple \(\Lambda\)-module \(S_i\) dual to \(P_i\) is given by
\[
    \textnormal{VdB}(\;S_i\;)=
    \begin{cases}
        \omega_{\ul C}[1], &i=0 \\
        \OO_{C_i}(-1), &i\in \Delta\setminus\JJ.
    \end{cases}
\]

\subsection{Tilting equivalences}
\label{subsec:brennerbutler}
The approach we take to describing algebraic t-structures on
\(\Dfl\Lambda\simeq \Dm X\) is via Iyama--Wemyss's
observation~\cite{iyamaTitsConeIntersections,iyamaMaximalModificationsAuslander2014}
that the tilting theory of \(\Lambda\) is controlled by the combinatorics of
certain reflexive \(R\)-modules.

\begin{definition}
    An \(R\)-module \(N\) is said to be \textit{modifying} if it is basic,
    reflexive, and its endomorphism algebra \(\End_R N\) is Cohen--Macaulay.
    We say \(N\) is an \emph{modifying generator} if in addition it contains
    \(R\) as a direct summand (equivalently if it is Cohen--Macaulay).
\end{definition}

The set of (isomorphism classes of) modifying \(R\)-modules is naturally
partitioned into \emph{mutation classes}~\cite[corollary
9.29]{iyamaTitsConeIntersections}, where two modifying modules \(M,N\) lie in
the same mutation class if they have the same number of indecomposable summands
and if \(M\) admits a two-term approximation by \(\text{add}(N)\). Each mutation
class then furnishes a family of derived-equivalent algebras via the following
Auslander--McKay type correspondence.

\begin{theorem}
    [{\cite[theorem 9.25]{iyamaTitsConeIntersections}}]
    \label{prop:IWAuslanderMcKay}
    If \(M,N\) are modifying \(R\)-modules in the same mutation class then
    \(\Hom_R(N,M)\) is a (classical) reflexive tilting \(\End_R N\)-module, and
    further every reflexive tilting \(\End_R N\)-module arises in this way from
    a unique module in the mutation class of \(N\).
\end{theorem}

There is a natural isomorphism \(\End_{\End_R N}(\Hom_R(N,M))\cong\End_R M\)
given by reflexive equivalence \cite[lemma
2.5]{iyamaMaximalModificationsAuslander2014}. Then as is standard, the tilting
module \(\Hom_R(N,M)\) in the above setup induces quasi-inverse derived
equivalences
\begin{equation}
     \Db \left(\End_RN\right)
    \xleftrightarrows[\quad\RHom(\Hom_R(N,M),-)\quad]{(-)\LTensor \Hom_R(N,M)}
    \Db \left(\End_RM\right).
    \label{eqn:brennerbutler}
\end{equation}

Being \(R\)-linear, the functors restrict to equivalences between the categories
\(\Dfl(\End_R N)\) and \(\Dfl(\End_R M)\) containing complexes with finite
length cohomology.

\paragraph{Mutating modifying modules} The key example of a modifying
generator is the module \(\pi_\ast \VV\XZ\), which Cohen--Macaulay with a
Cohen--Macaulay endomorphism algebra by~\cite[proposition
3.2.10]{vandenberghThreedimensionalFlopsNoncommutative2004}. Given the setup of
our problem, we fix once and for all this choice \(N=\pi_\ast \VV\XZ\), which
has endomorphism algebra \(\Lambda =\End_R N\cong \End_X\VV\XZ\). Note that the
indecomposable summands of \(N\) are naturally indexed over
\(\ul\Delta\setminus\JJ\) by \cref{lem:indexing-of-mmg-summands}.

Write \(\MaxMod R\) for the mutation class of \(N\), and \(\MaxModGen R\subset
\MaxMod R\) for the subset of modifying generators in the mutation class.
Iyama--Wemyss's correspondence above is complemented by the following result
which enumerates the entire mutation class of \(N\).

\begin{theorem}
    \label{thm:MMRDynkinlabelling}
    There is a bijection \(\CC\colon \MaxMod R\to
    \Chambers(\ul\Delta,\JJ)\) such that the following statements hold.
    \begin{enumerate}[(1)]
        \item \label{item:MMRDynkin1} Given \(M\in \MaxMod R\) and \(i\in
            \ul\Delta\setminus\mathbb{J}(\CC M)\), let \(M'\in \MaxMod R\) be
            such that \(\CC(M')\) is the adjacent chamber \(\nu_i\CC(M)\). Then
            there is a unique indecomposable summand \(M_i\subset M\) such that
            \(M/M_i\) is a summand of \(M'\); this defines a bijection between
            \(\ul\Delta\setminus \mathbb{J}(\CC M)\) and the indecomposable
            summands of \(M\). 
        \item If \(M\) above is a modifying generator, then
            the extended vertex \(0\in \ul\Delta\) lies outside \(\mathbb{J}(\CC
            M)\) and corresponds to the summand \(M_0=R\).
        \item For \(N=\pi_\ast\VV\XZ\), we have \(\CC(N)=\Chamb_\JJ^+\) and the
            bijection in \ref{item:MMRDynkin1} coincides with the
            \((\ul\Delta\setminus\JJ)\)-indexing of summands given by
            \cref{lem:indexing-of-mmg-summands}.
        \item Given \(M,M'\in \MaxMod R\) corresponding to the tilting \(\Lambda\)-modules
            \(T=\Hom(N,M)\) and \(T'=\Hom(N,M')\), we have \(\CC(M)\geq
            \CC(M')\) in the weak order on \(\Chambers(\ul\Delta,\JJ)\) if and
            only if \(\Ext^1(T,T')=0\).
    \end{enumerate}
    \begin{proof}
        If \(X\) is \(\mathbb{Q}\)-factorial (i.e.\ a minimal model) then this
        is~\cite[corollary 9.8]{iyamaTitsConeIntersections}. We sketch the
        general argument.

        By considering a general elephant \(\Spec \ol R \into \Spec R\) and the
        corresponding reduction functor \(\mathbb{F}{(-)}=(-)\otimes_R \ol R\),
        Iyama--Wemyss \cite[corollary 9.29]{iyamaTitsConeIntersections}
        construct the required bijection as a composite 
        \begin{equation}
            \label{eqn:conesfromMMR}
            \MaxMod R \xrightarrow{\;\Hom_R(N,-)\;}
            \{\text{basic ref.\ tilting \(\Lambda\)-modules}\}
            \xrightarrow{\;\bbF(-)\;}
            \{\text{basic tilting \(\bbF\Lambda\)-modules}\}
            \xrightarrow{\;g\text{-cone}\;}
            \Chambers(\ul\Delta,\JJ) 
        \end{equation}
        where the final map assigns each basic tilting module over \(\mathbb{F}\Lambda\)
        (a contracted preprojective algebra) to the cone spanned by its
        associated \(g\)-vectors. Note this requires an identification of the
        vector space \(\Theta(\ul\Delta\setminus\JJ)\) with
        \(\KK_\text{split}(\rproj \mathbb{F}\Lambda)\), for this we map the class of
        the indecomposable projective \(\mathbb{F}\Lambda_i=\mathbb{F}{\Hom_R(N,N_i)}\) to the
        vector \(\theta_i\in \Theta(\ul\Delta\setminus\JJ)\) dual to
        \(\alpha_i\in \mathfrak{h}(\ul\Delta\setminus\JJ)\).

        Now the assignment \(g\text{-cone}(-)\) has each of the required
        properties since wall-crossing is compatible with tilting mutation and
        \(g\text{-cone}(\mathbb{F} T)\geq g\text{-cone}(\mathbb{F}T')\) in the
        weak order if and only if \(\Ext^1(\mathbb{F}T,\mathbb{F}T')=0\)
        \cite[theorem 5.2]{iyamaReductionTriangulatedCategories2018}. 

        Since \(\Hom_R(N,-)\) and \(\mathbb{F}(-)\) both preserve
        indecomposability, and since \(\mathbb{F}\) preserves the tilting order
        \cite[proposition 4.2]{kimuraTiltingSiltingTheory2024}, it follows that
        the composite \(\CC(-)= g\text{-cone}(\mathbb{F}\Hom(N,-))\) has the
        required properties too.
    \end{proof}
\end{theorem}

The upshot is that \(\MaxMod R\) inherits both---a partial order (where
\(M\geq M'\) if and only if \(\CC(M)\geq \CC(M')\) in the weak order), and the
structure of a set with \(\ul\Delta\)-mutation (with exchange quiver isomorphic to that of
\(\Chambers(\ul\Delta,\JJ)\)). In particular for \(M,M'\in \MaxMod R\) there
exists a positive \(\mathbb{J}(M)\)-path \(\nu\) such that
\(M'=\nu M\), and further if \(M\) and \(M'\) are modifying generators then this
path \(\nu\) can be chosen to be spherical. Thus we have
\begin{equation}
    \label{eqn:mm-via-paths}
    \MaxMod R=\{\nu N\;|\; \nu \text{ a \(\JJ\)-path}\}, \quad
    \MaxModGen R=\{\nu N\;|\; \nu \text{ a spherical \(\JJ\)-path}\}.
\end{equation}

This extends to an enumeration of a large class of algebras and their tilting
modules. Given two positive \(\JJ\)-paths \(\nu,\mu\) we write \[{}_{\mu} \Lambda_{\nu}=\Hom_R(\nu N, \mu
N)\] and note that this is a reflexive tilting \(\End_R(\nu N)\)-module
by~\cref{prop:IWAuslanderMcKay}. We omit the empty path \(\emptyset\) from the
notation, so e.g.\ \(\Lambda_{\nu}=\Hom_R(\nu N,N)\) and
\({}_\nu\Lambda=\Hom_R(N,\nu N)\).  Happily, this is consistent with
\(\Lambda=\Hom_R(N,N)\).

Then each module \(\nu N \in \MaxMod R\) can be used to obtain two equivalences
to \(\Dfl\Lambda\). Indeed the constructions in \eqref{eqn:brennerbutler}
are symmetric in the input data so we obtain a tilting \(\Lambda\)-module
\({}_\nu \Lambda\) and a tilting \({}_\nu\Lambda_\nu\)-module \(\Lambda_\nu\),
which in turn give two \emph{mutation functors} \(\Dfl {}_\nu \Lambda_\nu
\rightrightarrows \Dfl \Lambda\) defined as
\begin{equation}
    \label{eqn:twomutationfunctors}
    \Phi_\nu(-)=(-)\LTensor {}_\nu \Lambda,
    \qquad \Psi_\nu(-)=\RHom(\Lambda_\nu,-).
\end{equation}

\paragraph{Composing mutation functors}
Of course for the purpose of defining mutation functors there is nothing special
about the choice of \(N\) and \(\Lambda=\End_R N\), and we can define mutation
functors starting at \(\Dfl(\End_R M)\) for any modifying module \(M\). That is
to say if \(M=\nu N\) for some positive \(\JJ\)-path \(\nu\), and if \(\mu\) is
a positive \(\nu\JJ\)-path, then the reflexive tilting modules
\({}_{\mu\nu}\Lambda_\nu\in\rmod({}_\nu\Lambda_\nu)\) and
\({}_{\nu}\Lambda_{\mu\nu}\in\rmod({}_{\mu\nu}\Lambda_{\mu\nu})\)
define two functors 
\[
    (\;\text{-}\;)\LTensor {}_{\mu\nu} \Lambda_\nu, \quad \RHom({}_\nu\Lambda_{\mu\nu},\;\text{-}\;)\colon\Dfl{}_{\mu\nu}\Lambda_{\mu\nu}\rightrightarrows \Dfl {}_\nu\Lambda_\nu
\] 
which we again denote by \(\Phi_\mu(\;\text{-}\;)\) and
\(\Psi_\mu(\;\text{-}\;)\) respectively. 

The following lemmas record composition relations between such functors. It will
suffice to consider functors corresponding to paths of length \(1\), and for
simplicity of notation we write \(\Phi_i\) (resp.\ \(\Psi_i\)) for the functor
\(\Phi_{\nu_i}\) (resp.\ \(\Psi_{\nu_i}\)).

\begin{lemma}
    \label{lem:mutationcomposition}
    Let \(M=\nu N\) be a modifying module determined by \(\JJ\)-path \(\nu\), and
    fix a vertex \(i\in \ul\Delta\setminus \mathbb{J}(M)\). 
    \begin{enumerate}[(1)]
        \item We have natural isomorphisms \(\Phi_i\circ \Psi_{\iota(i)} \simeq \Psi_i\circ
            \Phi_{\iota(i)} \simeq \textnormal{id}\) of functors
            \(\Dfl{}_\nu\Lambda_\nu\to\Dfl{}_\nu\Lambda_\nu\).
        \item If \(M>\nu_i M\), we have natural isomorphisms
            \(\Phi_{\nu_i\nu}\simeq \Phi_\nu\circ \Phi_i\) and
            \(\Psi_{\nu_i\nu}\simeq \Psi_{\nu}\circ \Psi_i\) of functors
            \(\Dfl{}_{\nu_i\nu}\Lambda_{\nu_i\nu}\to \Dfl\Lambda\).
    \end{enumerate}

    \begin{proof}
        The first statement is immediate from definitions. For the second, if
        \(\nu N>\nu_i \nu N\) in \(\MaxMod R\), then by~\cite[theorem
        4.6]{hiranoFaithfulActionsHyperplane2018} there is an isomorphism of
        \({}_{\nu_i\nu}\Lambda_{\nu_i\nu}-\Lambda\) bimodules
        \({}_{\nu_i\nu}\Lambda_\nu \LTensor {}_\nu\Lambda \simeq
        {}_{\nu_i\nu}\Lambda\) (where the tensor product is over
        \({}_\nu\Lambda_\nu\)). The isomorphism \(\Psi_{\nu_i\nu}\simeq
        \Psi_\nu\circ \Psi_i\) is then immediate, while the isomorphism
        \(\Phi_{\nu_i\nu}\simeq \Phi_\nu\circ \Phi_i\) follows from the
        \(\LTensor\)--\(\RHom\) adjunction~\cite[lemma
        2.10]{iyamaFominZelevinskyMutationTilting2008}.
    \end{proof}
\end{lemma}

\begin{lemma}
    \label{cor:mutationdecomp}
    Given a positive \(\JJ\)-path \(\nu=\nu_{i_n}...\nu_{i_1}\), we 
    can write \(\Phi_\nu = \Upsilon_{{i_1}}\circ \Upsilon_{{i_2}}\circ
    ...\circ \Upsilon_{{i_n}}\) for some sequence of simple mutations
    \(\Upsilon_{i}\in \{\Phi_{i}, \Psi_{i}\}\). Further if the path \(\nu\) is
    atomic (i.e.\ \(N>\nu_{i_1} N > \nu_{i_2}\nu_{i_1} N> ... > \nu N\)), then
    each \(\Upsilon_i\) above is of the form \(\Phi_i\) and we have \(\Phi_\nu =
    \Phi_{i_1}\circ \Phi_{i_2}\circ ... \circ\Phi_{i_n}\).

    The analogous statement holds for \(\Psi_\nu\).

    \begin{proof}
        For \(n\geq 1\) we may write \(\nu=\nu_{i}\nu'\) where \(i=i_n\)
        and \(\nu'=\nu_{i_{n-1}}...\nu_{i_1}\) is a shorter path. By
        \cref{lem:atomicproperties} we have either \(\nu'N > \nu N\) or
        \(\nu N > \nu' N\). In the first case
        \cref{lem:mutationcomposition} gives us \(\Phi_\nu= \Phi_{\nu'}\circ \Phi_{i}\),
        while in the latter case we write \(\nu'N=\nu_{\iota(i)}\nu N\) to get
        \(\Phi_\nu=\Phi_{\nu'}\circ (\Phi_{\iota(i)})^{-1}=\Phi_{\nu'}\circ
        \Psi_{i}\). The result follows by induction.
    \end{proof}
\end{lemma}

\subsection{Birational models}\label{subsec:geometrichearts}
Geometric t-structures naturally arise from birational
modifications of \(X\), where surgery operations described below replace
exceptional curves in a way that preserves the derived category.

\begin{theorem}[{\cite[theorem 1.1]{chenFlopsEquivalencesDerived2002}}]
    \label{thm:chenflops}
    Suppose \(Y\) is a normal 3-fold with at worst terminal Gorenstein
    singularities, and the maps \(\tau\colon X\to Y\), \(\tau'\colon W\to Y\) are
    flopping contractions over \(Y\) such that whenever \(D\) is a \(\tau\)-nef
    divisor on \(X\), the proper transform of \(-D\) across the birational map
    \(X\dashrightarrow W\) is \(\mathbb{Q}\)-cartier and \(\tau'\)-nef. Then
    there is a Fourier--Mukai type equivalence of derived categories \(\Db
    W\to \Db X\).
\end{theorem}

A map \(\tau'\) satisfying the given properties, if it exists, is unique and is
called the \emph{flop} of \(\tau\). We say the equivalence \(\Db W\to \Db X\) is
the \emph{Bridgeland--Chen flop functor}.

\begin{remark}
    There are various equivalent definitions of a flop in the literature, see
    for example \cite{kollarFlipsFlopsMinimal1990}. Given
    a flopping contraction \(\tau:X\to Y\) between normal 3-folds, one typically
    needs to choose a \(\mathbb{Q}\)-cartier \(\mathbb{Q}\)-divisor \(D\) on
    \(X\) such that \(-D\) is \(\tau\)-ample in order to define the \(D\)-flop
    (i.e.\ the \((K_X+D)\)-flip) of \(\tau\), which is unique and independent of
    the choice of \(D\) by \cite[corollary
    6.4]{kollarBirationalGeometryAlgebraic1998}. In the special case that \(X\)
    (equivalently \(Y\)) is terminal, the existence of the flop is proved
    in theorem 6.14 \emph{ibid.}
\end{remark}

Now the flopping contraction \(\pi:X\to Z\) of interest to us is over a complete
local base, and therefore we can freely contract exceptional curves~\cite[\S
2]{schroerCharacterizationSemiamplenessContractions2001}. That is to say, given any collection of
exceptional curves \(C_I=\bigcup_{i\in I}C_i\subset X\), the map \(\pi\) admits
a unique factorisation
\begin{equation}
    \label{eqn:partialcontraction}
    \begin{tikzcd}
        X \arrow[r, "\tau"] \arrow[rr, "\pi"', bend right=30]
        & Y \arrow[r, "\varpi"] & Z
    \end{tikzcd}
\end{equation}
where \(\tau\) contracts \(C_i\) to a point if and only if \(i\in I\). We say
\(\tau\) is a \emph{partial contraction over \(Z\)}. Since the singularity on
\(Z\) is assumed to be isolated, it follows that both \(\tau\) and \(\varpi\)
are flopping contractions. In particular we can consider the flop \(\tau'\colon
W\to Y\) of \(\tau\) to obtain a derived equivalence \(\Db W \to \Db X\); further the map
\(\varpi\circ \tau'\colon W\to Z\) is itself a flopping contraction over \(Z\) so
we can iterate.

Wemyss's homological minimal model programme
\cite{wemyssFlopsClustersHomological2018} addresses the question of classifying
\emph{birational models} of \(X\), i.e.\ varieties \(W\) obtained by iteratively
flopping subsets of \(\pi\)-exceptional curves. Abusing notation to write
\(\pi\colon W\to Z\) for the corresponding flopping contraction, we can again consider
the associated vector bundle \(\VV\WZ\), the modifying \(R\)-module generator
\(M=\pi_\ast \VV\WZ\), and the equivalence
\[
    \text{VdB}\colon \Db(\End M)\rightarrow \Db W
\]
constructed analogously to~\eqref{eqn:zeroperforpi}. In
this context, Wemyss shows that flops obey mutation-combinatorics following the
modifying generators.

\begin{theorem}[\cite{wemyssFlopsClustersHomological2018}]
    \label{thm:birMMGbijection}
    If \(\pi\colon W\to Z\) is a birational model of \(\pi\colon X\to Z\), then the
    associated modifying \(R\)-module generators \(M=\pi_\ast\VV\WZ\) and
    \(N=\pi_\ast\VV\XZ\) lie in the same mutation class, and every
    module in \(\MaxModGen R\) arises in this way from a unique birational model
    of \(X\).
\end{theorem}

Thus the set \(\Bir\XZ\) of birational models of \(X\) is in bijection with
\(\MaxModGen (R)\) via the map \(W\mapsto \pi_\ast\VV\WZ\).
Since~\cref{thm:MMRDynkinlabelling} shows \(\MaxModGen R\) is a set with
\(\Delta\)-mutation, \(\Bir\XZ\) inherits this structure. That is to say for any
\(W\in \Bir\XZ\) with associated modifying generator \(\pi_\ast\VV\WZ=M\) we
have \(\mathbb{J}(W)=\mathbb{J}(M)\), and for any \(\mathbb{J}(W)\)-path \(\nu\)
the birational model \(\nu W\) is the unique element of \(\Bir\XZ\) satisfying
\(\pi_\ast \VV\pair{\nu W}{Z} = \nu M\).

This mutation structure on \(\Bir\XZ\) can be realised intrinsically as arising
from flops---the indecomposable summands of
\(M\) are indexed over \(\ul\Delta\setminus\mathbb{J}(M)\) in a way that
applying \(\nu_i\) corresponds to mutating the \(i\)th summand, and the only
free summand is \(M_0 = R\). It follows that the non-free indecomposable
summands of \(\VV\WZ\) are indexed over \(\Delta\setminus\mathbb{J}(W)\). By
\cref{lem:indexing-of-mmg-summands} these summands are naturally in bijection
with the integral exceptional curves of \(W\). Therefore we see that the reduced
exceptional fiber of \(\pi\colon W\to Z\) can be written as a union of integral
curves \(\bigcup_{i\in \Delta\setminus\mathbb{J}(W)} C_i\). The following result
of Wemyss realises the operation \(\nu_i\) on \(\Bir\XZ\) as a flop of the
\(i\)th curve.

\begin{theorem}
    [{\cite[theorem 4.2]{wemyssFlopsClustersHomological2018}}]
    \label{thm:flopmutation}
    Fix \(i\in \Delta\setminus\JJ\), and suppose \(\nu_i X=W\to Z\) is the
    birational model of \(X\) corresponding to the modifying \(R\)-module
    generator \(M=\nu_i N\) with reduced exceptional fiber
    \(\bigcup_{i\in\Delta\setminus \nu_i \JJ}C_i\) indexed as above. Then \(W\)
    is precisely the flop of \(X\) in the curve \(C_i\), and
    \(C_{\iota(i)}\subset W\) is the proper transform of \(C_i\subset X\).
    Further, in this case the Bridgeland--Chen flop functor \(\Db W \to \Db X\)
    coincides with the composite equivalence
    \[
        \Db W
        \xrightarrow{\quad\textnormal{VdB}^{-1}\quad}
        \Db (\End M)
        \xrightarrow{\quad \Psi_i(-)\quad}
        \Db (\End N)
        \xrightarrow{\quad\textnormal{VdB}\quad}
        \Db X,
    \]
    and restricts to an equivalence \(\Dm W \to \Dm X\) between the full
    subcategories of complexes supported on exceptional curves.
\end{theorem}

In what follows, we omit \(\textnormal{VdB}\) from notation whenever convenient,
thus for example writing \(\Psi_{i}\) for both the mutation equivalence and
the flop functor in the above theorem.

By examining the connectedness of \(\ExQuiv(\MaxModGen R)\) we also conclude
that any two birational models of \(X\) are connected by a chain of
single-curve flops, i.e.\ given birational models \(W,W'\in \Bir\XZ\) there is a \(\mathbb{J}(W)\)-path
\(\nu=\nu_{i_n}\ldots\nu_{i_1}\) and a \emph{sequence of flops}
\[
    W \;\dashrightarrow\;
    \nu_{i_1}W \;\dashrightarrow\;
    \nu_{i_2}\nu_{i_1}W \;\dashrightarrow\;
    \cdots \;\dashrightarrow\;
    \nu_{i_n}\ldots\nu_{i_1}W=W'.
\]
Extending the notion of Bridgeland--Chen functors to chains of flops, in this
case we say the \emph{flop functor} associated to \(\nu\) is the equivalence
\[
    \Db \nu W
    \xrightarrow{\quad\textnormal{VdB}^{-1}\quad}
    \Db (\End \nu M)
    \xrightarrow{\quad \Psi_\nu(-) \quad}
    \Db (\End M)
    \xrightarrow{\quad\textnormal{VdB}\quad}
    \Db W
\]
which is evidently independent of the choice of \(\nu\). 

\subsection{Semi-geometric categories}
\label{subsec:partialperverse}
When the flopping contraction \(\pi\colon X\to Z\) is not irreducible (i.e.\ the
reduced exceptional fiber has multiple components), there are intermediate partial
contractions \(X\to Y\to Z\) to be considered. In this subsection we establish
basic structural results about perverse sheaves arising from such morphisms.

Thus fix a subset \(I\subset \Delta\setminus \JJ\), equivalently a collection of
exceptional curves \(C_I=\bigcup_{i\in I}C_i\) in \(X\) which can be
contracted via the map \(\tau\colon X\to Y\) as in \eqref{eqn:partialcontraction}. As
noted previously, \(\tau\) is a flopping contraction and in particular
\cite[theorem A]{vandenberghThreedimensionalFlopsNoncommutative2004} applies.
Thus there is a vector bundle \(\VV\XY\) on \(X\), furnishing a sheaf of
\(\OO_Y\)-algebras \(\ulLambda=\tau_\ast \sheafEnd\VV\XY\) such that there is a
derived equivalence
\begin{equation}
    \label{eqn:semigeomperverseequivalence}
    \textnormal{VdB}\colon 
    \Db\ulLambda
    \xrightarrow{\quad\tau^{-1}(-)\;\LTensor\; \mathscr{V}\XY\quad}
    \Db X
\end{equation}
where the tensor product is over \(\tau^{-1}\ulLambda\). We write
\(\zeroPer\XY\) for the image of \(\Coh \ulLambda\) under
this equivalence. By \cite[proposition
3.3.1]{vandenberghThreedimensionalFlopsNoncommutative2004}, this category admits
the description
\begin{equation}
    \label{eqn:zeroperdefn1}
    \zeroPer\XY = \left\{x\in \Coh X[0,1] \;\middle\vert\;
        \begin{array}{l}
            \mathbf{R}^1\tau_\ast (\HH^0x)=0,
            \quad \tau_\ast (\HH^{-1}x)=0,\\
            \Hom(c,\HH^{-1}x)=0 \text{ for all }c\in\Coh X
            \\ \text{satisfying }\mathbf{R}\tau_\ast c=0. 
        \end{array}
    \right\}.
\end{equation}
Since we are concerned with complexes supported on exceptional fibers,
we also define the full subcategories
\begin{gather*}
    \Dm\ulLambda  = \left\{y\in \Db\ulLambda \;|\;
    \mathop{\text{Supp}}y\subseteq \varpi^{-1}[\mathfrak{m}]\right\},\\
    \coh\ulLambda = \Coh\ulLambda \;\cap\; \Dm\ulLambda, \qquad
    \zeroper\XY = \zeroPer\XY\;\cap\; \Dm X.
\end{gather*}
The functor \(\text{VdB}\) clearly restricts to an equivalence
\(\Dm\ulLambda\to \Dm X\), thus identifying \(\coh\ulLambda\) with
\(\zeroper\XY\). The truncation functors associated to the heart \(\Coh
\ulLambda\subset \Db\ulLambda\) restrict to \(\Dm\ulLambda\),
showing that \(\coh\ulLambda \subset \Dm\ulLambda\) (and hence
\(\zeroper\XY\subset \Dm X\)) is the heart of a bounded t-structure.

\begin{remark}
    \label{rmk:extremeperversity}
    The following calibration is helpful---when \(I=\emptyset\) the map \(\tau\)
    is an isomorphism and we have \(\VV\XY=\ulLambda=\OO_X\), so that
    \(\zeroPer\XY=\Coh \ulLambda = \Coh X\). In this case we have \(\Dm
    \ulLambda=\Dm X\) by definition. On the other extreme, when
    \(I=\Delta\setminus\JJ\) the map \(\varpi\) is an isomorphism and the sheaf
    \(\ulLambda\) on \(Y=\Spec R\) corresponds to the \(R\)-algebra \(\Lambda\). Since
    \((R,\mathfrak{m})\) is complete local, a complex of \(R\)-modules has
    support in \([\mathfrak{m}]\) if and only if its total cohomology
    has finite length, i.e.\ \(\Dm\ulLambda=\Dfl\Lambda\).
\end{remark}

\paragraph{Characterising perversity} We now show that the `algebraic part' of
\(\zeroper\XY\subset \Dm X\) is generated by finitely many simple perverse
sheaves which depend only on a neighbourhood of the contracted curves, and
further can be read off as a subcategory of \(\zeroper\XZ\). For simplicity we
first consider the case when \(C_I\subset X\) is connected, i.e.\ \(\tau\) is an
isomorphism away from a point \(p=\tau(C_I)\in Y\). Write \(\ul
C_I=\tau^{-1}(p)\) for the scheme theoretic exceptional fiber of \(\tau\), this
has underlying reduced subscheme \(C_I\). Then we have the following.

\begin{theorem}
    \label{prop:perversitycheck}
    Let \(C_I\subset X\) be a connected component of the reduced exceptional
    fiber in \(X\), and \(\tau:X\to Y\) the crepant contraction of \(C_I\). For
    any complex of coherent sheaves \(x\in\Dm X\) supported within \(C_I\), the
    following are equivalent.
    \begin{enumerate}
        \item \label{item:perversitycheck1}
            The complex \(x\) lies in \(\zeroper\XZ\), i.e.\ \(x\) is perverse
            with respect to the contraction \(\pi\colon X\to Z\).
        \item \label{item:perversitycheck2}
            The complex \(x\) lies in \(\zeroper\XY\), i.e.\ \(x\)
            is perverse with respect to the contraction \(\tau\colon X\to Y\).
        \item \label{item:perversitycheck3}
            The complex \(x\) is filtered by the sheaves \(\OO_{C_i}(-1)\) for
            \(i\in I\) and the complex \(\omega_{\ul C_I}[1]\) (the
            suspended canonical sheaf of scheme-theoretic exceptional fiber of
            \(\tau\)).
    \end{enumerate}
\end{theorem}

To prove this, we first establish the following lemma which reduces
conditions involving the null-category \(\{c\in \Coh X\;|\; \mathbf{R}\tau_\ast
c=0\}\) (such as those appearing in the descriptions \eqref{eqn:zeroperdefn},
\eqref{eqn:zeroperdefn1}) to checks on a finite collection of objects.

\begin{lemma}
    \label{lem:zeropushforward}
    Let \(\tau:X\to Y\) be the crepant contraction of \(C_I\) and \(c\in \Coh
    X\) be such that \(\mathbf{R}\tau_\ast c=0\). Then \(c\) is filtered by the
    sheaves \(\OO_{C_i}(-1)\) for \(i\in I\).
    \begin{proof}
        Note if \(c\) is as given, then we have \(\mathbf{R}\pi_\ast c=0\) and
        hence \(c\in \zeroPer\XZ\) from the description \eqref{eqn:zeroperdefn}.
        Further \(c\) is supported within \(\pi^{-1}[\mathfrak{m}]\). Thus
        across the equivalence \eqref{eqn:zeroperforpi}, we see that \(c\)
        thus corresponds to some finite length \(\Lambda\)-module (which we
        again denote by \(c\)). In particular \(c\) is filtered by the simples
        of \(H=\zeroper\XZ\) and the problem is reduced to showing that the only
        simples which can occur in a composition series for \(c\) are the
        \(S_i\) for \(i\in I\).

        But by \cite[theorem 2.15]{wemyssFlopsClustersHomological2018},
        \(\Lambda\) can be expressed as the quotient of the path algebra of a
        quiver with vertex set \(\ul\Delta\setminus \JJ\), and the simple
        representation supported on vertex \(i\) is precisely the module
        \(S_i\). Further writing \(e_i\in \Lambda\) for the vertex idempotent at
        \(i\in \ul\Delta\setminus\JJ\), \cite[proposition
        2.14]{wemyssFlopsClustersHomological2018} shows that a
        \(\Lambda\)-module \(x\) is annihilated by \(\sum_{i\notin I}e_i\) if
        and only if the corresponding complex \(x\in \zeroPer\XZ\) satisfies
        \(\mathbf{R}\tau_\ast x=0\). In particular, the \(\Lambda\)-module \(c\)
        is annihilated by \(\sum_{i\notin I}e_i\), and hence is filtered only by
        the vertex simples \(S_i\) for \(i\in I\) as required.
    \end{proof}
\end{lemma}

We then have the following.

\begin{proof}[Proof of \cref{prop:perversitycheck},
    \ref{item:perversitycheck1}\(\iff\)\ref{item:perversitycheck2}]
    Let \(x\) be a complex of coherent sheaves with support in \(C_I\), and for
    brevity write \(x_i=\HH^{-i}(x)\) for its cohomology sheaves.
    Since the objects of both \(\zeroper\XY\) and \(\zeroper\XZ\) can be
    described as two-term complexes of coherent sheaves on \(X\), we can assume
    \(x_i=0\) unless \(i=0,1\).

    If \(x\in \zeroper\XY\), then by the description
    \eqref{eqn:zeroperdefn1} we see that \(\tau_\ast (x_1)=0\) and hence
    \(\pi_\ast (x_1)=0\).

    Likewise we have \(\Hom(\OO_{C_j}(-1),x_1)=0\) whenever \(C_j\) is
    contracted by \(\tau\), because \(\mathbf{R}\tau_\ast \OO_{C_j}(-1)=0\).
    On the other hand if \(C_j\) is not contracted by \(\tau\) and we have a
    non-zero morphism \(f:\OO_{C_j}(-1)\to x_1\), then the image \(\img f\)
    is a non-zero subsheaf of \(x_1\) supported within the finite collection
    of points \(C_j\cap (\bigcup_{i\in I}C_i)\). But in that case,
    \(\pi_\ast (\img f)\) is a non-zero subsheaf of \(\pi_\ast (x_1)=0\), a
    contradiction. Thus in fact \(\Hom(\OO_{C_j}(-1),x_1)=0\) for all
    exceptional curves \(C_j\), and hence by \cref{lem:zeropushforward} we
    have \(\Hom(c,x_1)=0\) whenever \(c\in \Coh X\) satisfies
    \(\mathbf{R}\pi_\ast c=0\).

    Lastly, examining the Leray spectral sequence
    \(\mathbf{R}^p\varpi_\ast\circ \mathbf{R}^q\tau_\ast(x_0)
    \Rightarrow \mathbf{R}^{p+q}\pi_\ast(x_0)\)
    (which degenerates since all maps have fiber dimension \(\leq 1\)) shows
    that \(\mathbf{R}^1\pi_\ast(x_0)\) is filtered by
    \(\varpi_\ast\circ \mathbf{R}^1\tau_\ast (x_0)\) (\(=0\) since
    \(\mathbf{R}^1\tau_\ast (x_0)\) vanishes) and
    \(\mathbf{R}^1\varpi_\ast\circ\tau_\ast(x_0)\) (\(=0\) since
    \(\tau_\ast (x_0)\) is supported on a zero-dimensional subset of \(Y\)).
    Thus we have \(\mathbf{R}^1\pi_\ast(x_0)=0\) as well, and hence \(x\in
    \zeroper \XZ\).

    Conversely if \(x\in \zeroper\XZ\), we see that \(\tau_\ast(x_1)\) is a
    coherent sheaf on \(Y\) with zero-dimensional support such that \(\varpi_\ast\circ
    \tau_\ast(x_1)=\pi_\ast(x_1)=0\). This is possible only if
    \(\tau_\ast(x_1)=0\). Likewise since
    \(\mathbf{R}^1\pi_\ast(x_0)=0\), the Leray sequence shows we have
    \(\varpi_\ast \circ \mathbf{R}^1\tau_\ast(x_0)=0\) and hence
    \(\mathbf{R}^1\tau_\ast(x_0)=0\). Lastly if \(c\in \Coh X\) is such that
    \(\mathbf{R}\tau_\ast c=0\) then we have
    \(\mathbf{R}\pi_\ast c = \mathbf{R}\varpi_\ast \circ \mathbf{R}\tau_\ast
    c = 0\), and hence \(\Hom(c,x_1)=0\). This shows \(x\in \zeroper\XY\).
\end{proof}

To prove \ref{item:perversitycheck2}\(\iff\)\ref{item:perversitycheck3}, we
reduce the problem to a formal neighbourhood \(Z_I\into Y\) of \(p\), i.e.\
the spectrum of the complete local ring \(\widehat{\OO}_{Y,p}\). The scheme
\(Z_I\) has an isolated cDV singularity at its closed point, and the restriction
of \(\tau\) is a flopping contraction \(X_I\to Z_I\) with reduced exceptional
fiber \(C_I\subset X_I\). We can therefore consider the full subcategory \(\Dm
X_I\subset \Db X_I\) of complexes supported on \(C_I\).

\begin{lemma}
    The restriction (i.e.\ pullback) functor \(\Dm X \to \Dm X_I\) gives an
    equivalence of \(\Dm X_I\) with the full subcategory of
    complexes in \(\Dm X\) supported within \(C_I\).

    \begin{proof}
        The schemes \(X\) and \(X_I\) have isomorphic completions along \(C_I\),
        and we write \(\mathfrak{X}\) for the associated Noetherian formal
        scheme. Now \cite[lemma 2.1]{orlovFormalCompletionsIdempotent2011} shows
        that the category \(\Dm X_I\) is the bounded derived category of the
        Abelian category
        \[
            \coh (X_I) = \{x\in \Coh X_I\;|\; \Supp(x)\subset C_I\},
        \]
        whilst the subcategory of \(\Dm X\) containing complexes supported
        within \(C_I\) is the bounded derived category of
        \[
            \Coh_{C_I}(X) = \{x\in \Coh X\;|\; \Supp(x)\subset C_I\}.
        \]
        But proposition 2.8 \emph{ibid}.\ shows
        that both categories above are equivalent (via pullback along
        the canonical map associated to completion) to the category of torsion coherent
        sheaves on \(\mathfrak{X}\). In particular, the restriction functor
        \(\Coh_{C_I}(X)\to \coh(X_I)\) is an exact equivalence and passing to bounded
        derived categories gives the result.
    \end{proof}
\end{lemma}

Thus we can identify \(\Dm X_I\) and its subcategories \(\coh X_I\),
\(\zeroper\pair{X_I}{Z_I}\) with their images in \(\Dm X\). In particular,
\(\zeroper\pair{X_I}{Z_I}\) (as a subcategory of \(\Dm X\)) is an algebraic
Abelian category equal to the extension-closure of its simple objects
\(\omega_{\ul C_I}[1]\), \(\OO_{C_i}(-1)\) for \(i\in I\). The following is then
immediate.

\begin{proof} [Proof of \cref{prop:perversitycheck},
    \ref{item:perversitycheck2}\(\iff\)\ref{item:perversitycheck3}]
    By \cite[proposition
    3.1.4]{vandenberghThreedimensionalFlopsNoncommutative2004}, the membership
    of any complex in \(\zeroper\XY\) can be checked locally with respect to the
    flat topology on \(Y\), in particular with respect to the flat cover
    \(Y=(Y\setminus \{p\})\cup Z_I\). Now if \(x\) is supported within contracted
    curves, the restriction of \(x\) to \(Y\setminus\{p\}\) vanishes and thus
    \(x\) lies in \(\zeroper\XY\) if and only if its restriction to \(X_I\)
    (i.e.\ \(x\) viewed as an object of \(\Dm X_I\)) lies in
    \(\zeroper\pair{X_I}{Z_I}\), equivalently if \(x\) is filtered by the
    simples of \(\zeroper\pair{X_I}{Z_I}\).
\end{proof}

It is straightforward to generalise this to partial
contractions \(\tau\colon  X\to Y\) associated to possibly disconnected collections of
exceptional curves \(C_I\subset X\), since any complex in \(\Dm X\) decomposes
into a direct sum of complexes with connected support. In particular we can
define the category 
\[ 
\zeropersupp{I}\XZ = \left\{x\in \zeroper\XZ \;\middle\vert\; \Supp x \subseteq
\bigcup_{i\in I}C_i\right\} 
\] 
and note that if \(C_I\) has connected components \(C_{I_1},...,C_{I_K}\), then this decomposes as
\[ 
\zeropersupp{I}\XZ =
\smashoperator[r]{\bigoplus_{J\in\{I_1,...,I_k\}}}\;\;\zeropersupp{J}\XZ. 
\] 
But following \cref{prop:perversitycheck}, each \(\zeropersupp{J}\XZ\) is
equivalent to the category \(\zeroper\pair{X_J}{Z_J}\) associated to a formal
neighbourhood \(Z_J\) of the point \(p_J=\tau(C_J)\in Y\), and in particular is
generated by finitely many complexes which depend only on the scheme theoretic
exceptional fiber \(\ul C_J\) overlying \(C_J\). Further, considering the flat
cover \(Y=Z_{I_1}\cup \cdots\cup Z_{I_k}\cup
(Y\setminus\{p_{I_1},\ldots,p_{I_k}\})\) shows that a complex \(x\) supported
within \(C_I\) lies in \(\zeroper\XY\) if and only if its restriction to each
\(X_J\) lies in \(\zeroper\pair{X_J}{Z_J}\), and this gives us the description
\begin{align*}
    \zeropersupp{I}\XZ
    = \bigg\langle \{\OO_{C_i}(-1)\;|\; i\in I\} \;\cup\; \{\omega_{\ul
    C_{J}}[1]\;|\; J = I_1,\ldots,I_k\} \bigg\rangle.
\end{align*}

A consequence of the equivalence is that the property of a two-term complex
\(x\) of coherent sheaves being `perverse' can be checked with respect to any
partial contraction which maps \(\mathop\text{Supp}(x)\) to a (finite collection
of) points. That is to say, if \(\tau\colon X\to Y\) and \(\tau'\colon X\to Y'\)
are two crepant partial contractions which contract all the curves containing
\(\mathop\text{Supp}(x)\), then we have \(x\in \zeroper\XY\) if and only if
\(x\in \zeroper\pair{X}{Y'}\).

Since \(\zeroper\XY\) away from the contracted curves should mimic
\(\zeroper\pair{X}{X}=\coh X\), we then have the following result which reduces the
construction of \(\zeroper\XY\) to a binary choice between \(\coh X\) and
\(\zeroper\XZ\) on each \(C_i\), effectively eliminating the need to consider
the geometry of the non-commutative scheme \((Y,\ulLambda)\).

\begin{theorem}
    \label{thm:structureofzeroperxy}
    If \(\tau\colon X\to Y\) is the crepant contraction of the exceptional subset
    \(C_I\subset X\), then the associated heart \(\zeroper\XY\) is the smallest
    extension-closed subcategory of \(\Dm X\) containing the full subcategory
    \(\zeropersupp{I}\XZ\) and all coherent sheaves on \(X\) that are supported within the
    uncontracted curves \(\bigcup_{i\notin I}C_i\). In other words,
    \[
        \zeroper\XY =
        \bigg\langle
            \underbrace{\left\{ x\in \zeroper\XZ\;\middle\vert\; \Supp x \subseteq
            \bigcup_{i\in I}C_i\right\}}_{\zeropersupp{I}\XZ}
            \;\cup\;
            \left\{ x\in \coh X\;\middle\vert\; \Supp x \subseteq \bigcup_{i\notin I}C_i\right\}
        \bigg\rangle.
    \]
    \begin{proof}
        The category \(\zeropersupp{I}\XZ\) clearly lies in \(\zeroper\XY\),
        while if \(x\in \coh X\) is supported within \(\bigcup_{i\notin I}C_i\)
        then \(\mathbf{R}\tau_\ast x\) is a sheaf on \(Y\) and hence \(x\in
        \zeroper\XY\). To furnish the required description of \(\zeroper\XY\),
        it thus suffices to show that every object in the category is an
        extension of complexes of the two given forms.

        Now an arbitrary complex \(x\in \zeroper\XY\) can be written as an
        extension of its cohomology objects with respect to \(\coh X\), namely
        \(x_0=\HH^0(x)\) and \(x_1[1]=\HH^{-1}(x)[1]\). It is clear from
        \eqref{eqn:zeroperdefn1} that \(x_0\), \(x_1[1]\) are themselves
        contained in \(\zeroper\XY\). Further we have \(\tau_\ast(x_1)=0\) and
        hence \(x_1\) must be supported within the contracted curve
        \(\bigcup_{i\in I}C_i\), from which it follows that \(x_1[1]\in
        \zeropersupp{I}\XZ\).

        On the other hand, writing \(\mathscr{I}\subset \OO_X\) for the ideal
        sheaf of the closed subscheme \(\bigcup_{i\in I}C_i\), note that
        \(y=\mathscr{I}^n x_0\) (i.e.\ the image of the natural map
        \(\mathscr{I}^n\otimes x_0\to x_0\)) is a subsheaf of \(x_0\)  (in
        particular, an object of \(\coh X\)) supported within \(\bigcup_{i\notin
        I}C_i\) for \(n\gg 0\). The quotient \(x_0/y\) is clearly supported
        within \(\bigcup_{i\in I} C_i\). Further, applying \(\tau_\ast\) to this
        quotient and examining the long exact sequence shows
        \(\mathbf{R}^1\tau_\ast(x_0/y)\) vanishes since
        \(\mathbf{R}^1\tau_\ast(x_0)\) does, i.e. \(x_0/y\) lies in
        \(\zeroper\XY\) (and hence in \(\zeropersupp{I}\XZ\)). 

        This concludes our proof that \(x\) is an extension of objects of the
        required form.
\end{proof}
\end{theorem}

It is immediate from \cref{thm:structureofzeroperxy} that an object of
\(\zeroper\XY\) is simple if and only if it is simple in some summand
\(\zeropersupp{J}\XZ\subset \zeropersupp{I}\XZ\) or a simple coherent sheaf
supported on the uncontracted locus, i.e.\ we have the following.

\begin{corollary}
    \label{cor:simplesofpartialperverse}
    For \(\tau\colon X\to Y\) the crepant contraction of \(C_I\subset X\), the simple
    objects of \(\zeroper\XY\) are the skyscrapers \(\OO_p\) at closed
    points \(p\notin \bigcup_{i\in I}C_i\), the sheaves \(\OO_{C_j}(-1)\) for
    \(j\in I\), and suspended canonical sheaves \(\omega_{\ul C_J}[1]\) for each
    scheme-theoretic exceptional fiber \(\ul C_J\) of \(\tau\).
\end{corollary}

\section{Numerical intermediate hearts}\label{sec:numtors}

We begin the systematic investigation of t-structures on \(\Dm X\) that are
\emph{intermediate} with respect to the perverse t-structure, i.e.\
whose heart is contained in \(\zeroper\XZ[-1,0]\). Noting hearts of t-structures
have a natural partial order \((\leq)\) given by containment of co-aisles,
the set of intermediate hearts we are concerned with is precisely the closed interval
\[\tilt(\zeroper\XZ)=\left\{K\;\middle\vert\; \zeroper\XZ[-1]\leq K \leq
\zeroper\XZ\right\}.\]
The partial order on \(\tilt(\zeroper\XZ)\) corresponds to containment of torsion-free classes across the
Happel--Reiten--Smal\o\ bijections \cite{happelTiltingAbelianCategories1996}.
That is to say, if we write \(\tors(\zeroper\XZ)\) (resp.\
\(\torf(\zeroper\XZ)\)) for the set of torsion (resp.\ torsion-free) classes in
\(\zeroper\XZ\), then we have order-preserving bijections
\begin{equation}
    \label{eqn:posetisoms1}
    \begin{tikzcd}[row sep=huge]
        \left(\tors(\zeroper\XZ),\subseteq\right)^\text{op}
            \arrow[rr,"(-)\orth",shift left]
            \arrow[ddr, bend right=15, shift right,
                   "{(-)\orth\;\ast\;(-)[-1]}"']
        & & \left(\torf(\zeroper\XZ),\subseteq\right)
            \arrow[ll,"\orth(-)", shift left]
            \arrow[ldd, bend left=15, shift left,
                   "{(-)\;\ast\;\orth(-)[-1]}"]\\\\
        & \left(\tilt(\zeroper\XZ),\leq\right)
            \arrow[luu, bend left=15, shift right,
                   "{(-)[1]\cap \zeroper\XZ}"']
            \arrow[uur, bend right=15, shift left,
                   "{(-)\cap \zeroper\XZ}"] &
    \end{tikzcd}
\end{equation}
where the notation \(U\orth\) (resp.\ \(\orth U\)), for \(U\subset
\zeroper\XZ\), denotes the full subcategory of all objects in \(\zeroper\XZ\)
that have no non-zero morphisms from (resp.\ to) objects of \(U\).

In this section we furnish a complete description of the sub-poset of hearts in
\(\tilt(\zeroper\XZ)\) which are \emph{numerical}, i.e.\ `detected by
\(\KK\)-theory' in the sense of \cref{def:numerical} below.

\subsection{The heart fan} Fix a triangulated category \(\mathscr{T}\),
and a surjection \(\KK \mathscr{T}\onto \mathfrak{h}\) of the Grothendieck group of
\(\mathscr{T}\) onto a free Abelian group of finite rank. Taking
\(\mathbb{R}\)-linear duals, note that the finite dimensional Euclidean space
\(\Theta = \Hom_\mathbb{Z}(\mathfrak{h}, \mathbb{R})\) thus injects into
\(\Hom_\mathbb{Z}(\KK\mathscr{T},\mathbb{R})\) so that any vector \(\theta\in
\Theta\) can be regarded as an \(\mathbb{R}\)-linear functional on
\(\KK\mathscr{T}\).

\begin{definition}
    \label{def:numerical}
    Given the heart of a t-structure \(H\subset \mathscr{T}\), the \emph{heart
    cone of \(H\)} (with respect to \(\mathfrak{h}\)) is defined to be the
    subset 
    \( \CC(H)=\{ \theta\in \Theta \;|\; \theta[h]\geq 0\text{ for all }h\in H\}
    \).

    We say \(H\) is \emph{numerical} if \(\CC(H)\) contains a non-zero vector.
\end{definition}

For instance, if \(H\) is \emph{algebraic} (i.e.\ has Jordan--H\"older
filtrations and finitely many simple objects) then \(\CC(H)\) is an orthant in
\(\Theta\), spanned by the dual vectors to the basis of \(\mathfrak{h}\) given
by classes of simples in \(H\). 

More generally \(\CC(H)\subset\Theta\) is a
\emph{rational cone}, i.e.\ a closed convex subset defined by inequalities with
coefficients in \(\mathfrak{h}\). If a rational cone \(\sigma\subset \Theta\)
lies in the half space \(\{\alpha\geq 0\}\) defined by
\(\alpha\in\mathfrak{h}\), we say the cone \(\sigma\cap \{\alpha=0\}\) is a
\emph{face} of \(\sigma\), and write \(\textnormal{faces}(\sigma)\) for the set
of all faces of \(\sigma\). A collection of cones \(\Sigma\) is a \emph{fan} in
\(\Theta\) if it is closed under taking faces and the intersection of any two
cones in \(\Sigma\) is a face of both. We say the fan is \emph{complete} if its
cones cover \(\Theta\), i.e.\ \(\bigcup_{\sigma\in \Sigma}\sigma=\Theta\).

\begin{theorem}[{\cite[theorem A]{broomheadHeartFanAbelian2024}}]
    \label{thm:heartfandefn}
    For \(H\), \(\mathscr{T}\), \(\mathfrak{h}\), and \(\Theta\) as above, the
    collection of cones
    \[
        \HFan(H) =
            \bigcup_{K\in \tilt(H)}\mathop\textnormal{faces}\left(\CC K \right)
    \]
    is a fan in \(\Theta\), called the \emph{heart fan} of \(H\) (with respect
    to \(\mathfrak{h}\)). If \(H\) is algebraic, then this fan is complete. If
    in addition the map \(\KK\mathscr{T}\onto \mathfrak{h}\) is an isomorphism,
    then an intermediate heart \(K\) is algebraic if and only if \(\CC K\) is a
    full-dimensional cone, and in this case \(K\) is the unique intermediate
    heart with heart cone \(\CC K\).
\end{theorem}

The existence of the heart fan, our primary bookkeeping tool, is complemented by
the following observation that computes the `fibers' of the map \(\CC\colon
\tilt(H)\to \HFan(H)\). To introduce some notation, if \(H\) is algebraic, then
note that any \(\theta\in\Theta\) induces a tripartition \[H=H_\tr(\theta)\ast
H_\ss(\theta)\ast H_\tf(\theta)\] where \(H_\tr(\theta)\) (resp.\
\(H_\tf(\theta)\)) are the \emph{numerical torsion} (resp.\ \emph{torsion-free})
classes
\begin{align*}
    H_\tr(\theta)&=\{h\in H\;|\; \theta[f]<0 \text{ for all non-zero
    factors }h\twoheadrightarrow f\neq 0\}, \\
    H_\tf(\theta)&=\{h\in H\;|\; \theta[s]>0 \text{ for all non-zero
    sub-objects }0\neq s\hookrightarrow h\} \\ 
    \intertext{introduced by Baumann--Kamnitzer--Tingley in
    \cite{baumannAffineMirkovicVilonenPolytopes2014}, while \(H_\ss(\theta)\) is
    the wide subcategory of \emph{\(\theta\)-semistable objects} \'a la King
    \cite{kingModuliRepresentationsFinite1994} given as}
    H_\ss(\theta)&= \{h\in H\;|\; \theta[h]=0, \; \theta[s]\geq 0 \text{ for
    all sub-objects }s\into h\}.
\end{align*}

Hearts \(K\in \tilt(H)\) with \(\theta\in\CC (K)\) arise from torsion
pairs that coarsen this tripartition.

\begin{theorem}[{\cite[corollary 3.3]{broomheadHeartFanAbelian2024}}]
    \label{thm:allheartsonacone}
    Suppose \(H\) is the heart of an algebraic t-structure on the triangulated
    category \(\mathscr{T}\), and \(\theta\) is a vector in \(\Theta\) as above. For any torsion pair
    \(H_\ss(\theta)=U\ast V\) on the wide (Abelian) subcategory of
    \(\theta\)-semistable objects in \(H\), the expression 
    \begin{equation}
    \label{eqn:numtiltsbyUV}
        K=V\ast H_\tf(\theta)\ast H_\tr(\theta)[-1]\ast U[-1]
    \end{equation}
    defines a heart in \(\tilt(H)\) that satisfies \(\theta\in \CC K\). Each
    \(K\in \tilt(H)\) whose heart cone contains \(\theta\) arises in this
    fashion, in particular there is a poset isomorphism
    \[
        \tilt(H_\ss\theta)\longrightarrow \{K\in \tilt(H)\;|\; \theta\in
        \CC K\}
    \]
    where the inverse map sends \(K\in \tilt(H)\) to \(K_\ss(\theta)=\{k\in K\;|\;
    \theta[k]=0\}\subset \tilt(H_\ss\theta)\).
\end{theorem}

Setting \(U=0, V=H_\ss(\theta)\) in~\eqref{eqn:numtiltsbyUV} yields the maximal
element \(H^\tt(\theta)\) in the subposet \(\{K\;|\; \theta\in \CC K\}\subset
\tilt(H)\), and it can be characterised as follows.

\begin{corollary}
    \label{lem:characterisationsofHtt}
    Given an intermediate heart \(K\in \tilt(H)\) and a vector
    \(\theta\in \CC (K)\), the following statements are equivalent.
    \begin{enumerate}[(1)]
        \item The heart \(K\) is maximal among tilts of \(H\) with \(\theta\) in
            their heart cone, i.e.\ \(K=H^\tt(\theta)\).
        \item We have the equality \(H_\ss(\theta)=K_\ss(\theta)\).
        \item The category \(H_\ss(\theta)\) lies in \(K\).
        \item The category \(K_\ss(\theta)\) lies in \(H\).
    \end{enumerate}
\end{corollary}

Likewise setting \(U=H_\ss(\theta),V=0\) in~\eqref{eqn:numtiltsbyUV} yields the
minimal heart \(H_\tt(\theta)\in \tilt(H)\) containing \(\theta\) in its heart
cone. Evidently \(\{K\in \tilt(H)\;|\;
\theta\in \CC (K)\}\) coincides with the interval \([H_\tt\theta,
H^\tt\theta]\subseteq \tilt(H)\).

We say a vector \(\theta\) is \emph{generic} in a cone \(\sigma\) if it
lies in \(\sigma\) but not in any proper face of \(\sigma\). If
\(\theta,\theta'\) are both generic in some cone \(\sigma\in \HFan(H)\), the
above results imply that all their associated numerical categories coincide
i.e.\ \(H_\tr(\theta)=H_\tr(\theta')\), \(H_\tf(\theta)=H_\tf(\theta')\), and so
forth. Thus we use the notations \(H_\tr(\sigma)\), \(H_\tf(\sigma)\),
\(H_\ss(\sigma)\), \(H_\tt(\sigma)\), and \(H^\tt(\sigma)\) to denote the
respective categories associated to any generic vector \(\theta\in \sigma\).

\begin{remark}[Zooming into the heart fan]
    \label{rmk:ZoomIntoHFan}
    Given a parameter \(\theta\in \Theta\), the correspondence
    \ref{thm:allheartsonacone} can be used to essentially `read off' the heart
    fan of \(H_\ss(\theta)\) from that of \(H\). To see this, note that the
    exact inclusion \(H_\ss(\theta)\into H\) induces a map of Grothendieck
    groups which can be composed with the surjection onto \(\mathfrak{h}\) to
    obtain a commuting square
    \[\begin{tikzcd}
        \KK H_\ss \rar[hookrightarrow] \dar[twoheadrightarrow] & \KK H \dar[twoheadrightarrow] \\ 
        \mathfrak{h}_\ss \rar[hookrightarrow]& \mathfrak{h},
    \end{tikzcd}\]
    where \(\mathfrak{h}_\ss\) is again a free Abelian group of finite rank.
    The heart fan of \(H_\ss(\theta)\) with respect to \(\mathfrak{h}_\ss\) thus
    lives in the vector space
    \(\Theta_\ss=\Hom_\mathbb{Z}(\mathfrak{h}_\ss,\mathbb{R})\) which is
    naturally a quotient of \(\Theta\). It can be shown that the surjection
    \(\Theta\onto\Theta_\ss\) is such that the image of (the
    sub-fan generated by)
    \[
        \{\sigma\in\HFan(H)\;|\; \theta \text{ lies in some face of }\sigma\}
    \]
    is a sub-fan of \(\HFan(H_\ss\theta)\). More precisely, given a heart \(K\in
    [H_\tt\theta,H^\tt\theta]\) and the corresponding category
    \(K_\ss(\theta)\in \tilt(H_\ss\theta)\), it can be shown that the image of
    \(\CC (K)\) is always a face of \(\CC(K_\ss\theta)\) and is in fact equal to
    \(\CC(K_\ss\theta)\) if \(K\) is algebraic.

    To see an example where the image of \(\CC K\) is a proper face of
    \(\CC(K_\ss)\), one should examine the category
    \(H=\mathop\text{Rep}(\text{\small{2}}\rightrightarrows\text{\small{1}})\)
    associated to the \(2\)-Kronecker quiver \cite[example
    4.6]{broomheadHeartFanAbelian2024}, and the stability parameter
    \(\theta=(-1,1)\) given in the basis dual to vertex simples.
\end{remark}

Fixing the algebraic heart \(H=\zeroper\XZ\simeq \rflmod\Lambda\) in \(\Dm
X\simeq \Dfl\Lambda\), we will use the above results to compute the complete
subposet of numerical hearts in \(\tilt(H)\). The \(\JJ\)-cone arrangement, a
complete simplicial fan arising form the Dynkin data associated to \(X\),
naturally acts as the heart fan in this setting.

\begin{theorem}\label{cor:heartfanofzeroper}
    Upon identifying the Grothendieck group \(\KK X\) with the restricted
    root lattice \(\mathfrak{h}(\ul\Delta\setminus \JJ)\) in a way that each
    simple root \(\alpha_i\) (\(i\in\ul\Delta\setminus\JJ\)) corresponds to the \(i\)th simple of
    \(H=\zeroper\XZ\) (indexed as in \cref{subsec:VdBequiv}), the heart fan of
    \(H\) is given by the \(\JJ\)-cone arrangement \(\Arr(\ul\Delta,\JJ)\) in
    \(\Theta(\ul\Delta\setminus\JJ)\simeq \Hom(\KK X,\mathbb{R})\).
\end{theorem}

In order to prove \cref{cor:heartfanofzeroper} it will suffice to identify, for
each cone \(\sigma\) of the form \(\nu\Chamb_\JJ^\pm\) or \(\nu
\Chamb_\JJ^0\), a heart \(K\in \tilt(H)\) with heart cone \(\sigma\). This is
accomplished in \cref{thm:heartconeofmutated,prop:heartconeofpartialperverse}.
We then examine the categories \(H_\ss(\sigma)\) for each (not necessarily
maximal) cone \(\sigma\) in the \(\JJ\)-cone arrangement
(\cref{subsec:semistablecategories}) to obtain
a complete list of numerical intermediate hearts.

\subsection{Algebraic intermediate hearts}\label{subsec:algebraicheartcones} 
Each modifying module \(\nu N \in
\MaxMod R\) yields a tilting \(\Lambda\)-module \({}_\nu \Lambda\) and a
tilting \({}_\nu\Lambda_\nu\)-module \(\Lambda_\nu\), defining the 
functors \(\Dfl {}_\nu \Lambda_\nu \rightrightarrows \Dfl \Lambda\) given as in
\eqref{eqn:twomutationfunctors}. Applying the Brenner--Butler theorem~\cite[chapter 4, remark
3.10]{angelerihugelHandbookTiltingTheory2007} to these equivalences 
gives two torsion pairs
\begin{equation}
    \label{eqn:mutatedtorsion}
    \rflmod \Lambda = T_\nu\ast F_\nu =U_\nu \ast V_\nu,
\end{equation}
\vspace{-2em}
\begin{alignat*}{3}
    T_\nu&=\{x\in \rflmod\Lambda\;|\; \Ext^1({}_\nu\Lambda,x)=0\}
         & \qquad
    U_\nu&=\{x\in \rflmod\Lambda\;|\; x\otimes \Lambda_\nu=0   \}\\
    F_\nu&=\{x\in \rflmod\Lambda\;|\; \Hom({}_\nu\Lambda,x)=0  \}
         & \qquad
    V_\nu&=\{x\in \rflmod\Lambda\;|\; \Tor_1(x,\Lambda_\nu)=0  \}
\end{alignat*}
such that the squares below commute.
\[
    \begin{tikzcd}[column sep=0em]
                 & \Dfl {}_\nu\Lambda_\nu &\hspace{3em}& \Dfl\Lambda\\
                 & \rflmod {}_\nu\Lambda_\nu && F_\nu[1]\ast T_\nu
                 \arrow["\Phi_\nu","\sim"', from=1-2, to=1-4]
                 \arrow["\sim", from=2-2, to=2-4]
                 \arrow[hook, from=2-2, to=1-2]
                 \arrow[hook, from=2-4, to=1-4]
    \end{tikzcd} \qquad\qquad
    \begin{tikzcd}[column sep=0em]
                 & \Dfl {}_\nu\Lambda_\nu &\hspace{3em}& \Dfl\Lambda\\
                 & \rflmod {}_\nu\Lambda_\nu && V_\nu\ast U_\nu[-1]
                 \arrow["{\Psi_\nu}","\sim"', from=1-2, to=1-4]
                 \arrow["\sim", from=2-2, to=2-4]
                 \arrow[hook, from=2-2, to=1-2]
                 \arrow[hook, from=2-4, to=1-4]
    \end{tikzcd}
\]

Abusing notation to write \(H\) for \(\rflmod\Lambda\subset \Dfl\Lambda\) as well as
\(\rflmod{}_\nu\Lambda_\nu\subset \Dfl{}_\nu\Lambda_\nu\), each \(\JJ\)-path
\(\nu\) thus gives us two hearts in \(\Dfl\Lambda\) that are intermediate with
respect to \(\rflmod\Lambda\), namely
\begin{equation}
    \label{eqn:intermediatehearts}
    \Phi_\nu H[-1]=F_\nu\ast T_\nu[-1], \qquad \Psi_\nu H = V_\nu\ast
    U_\nu[-1].
\end{equation}

To compute the heart cones of these, we now show that the
\(\KK\)-theoretic maps induced by the mutation functors
\(\Phi_\nu,\Psi_\nu\) obey the same combinatorial rules which govern
mutations of root lattices. Note that for any \(\JJ\)-path
\(\nu\), \cref{thm:MMRDynkinlabelling} indexes the summands
of \(\nu N\) (and hence the simple \({}_\nu\Lambda_\nu\)-modules) over the
vertices of \(\ul\Delta\setminus \nu\JJ\). This identifies
\(\KK({}_\nu\Lambda_\nu)\) with the restricted root lattice
\(\mathfrak{h}(\ul\Delta\setminus\nu\JJ)\) by mapping the classes of simples to
corresponding simple restricted roots. We then have the following.

\begin{lemma}
    \label{lem:mutationtocombirelation}
    For \(i\in\ul\Delta\setminus\JJ\), the simple mutation functors
    \(\Phi_i,\Psi_i\colon\Dfl{}_{\nu_i}\Lambda_{\nu_i}\to \Dfl\Lambda\) both
    induce the same map on \(\KK\)-theory and this agrees with the map
    \(\varphi_i\colon \mathfrak{h}(\ul\Delta\setminus\nu_i\JJ)\to
    \mathfrak{h}(\ul\Delta\setminus\JJ)\) given as in
    \eqref{eqn:wallcrossingdefn}.

    \begin{proof}
        This is essentially~\cite[remark 5.2]{nabijouGVGWInvariants2023}, we
        spell out the details for clarity. By
        \cite[lemma 3.7]{hiranoStabilityConditions3fold2023}, both \(\Phi_i\)
        and \(\Psi_i\) induce the same map on \(\KK\)-theory which we
        momentarily call \(\psi_i\colon \mathfrak{h}(\ul\Delta\setminus \nu_i\JJ)\to
        \mathfrak{h}(\ul\Delta\setminus\JJ)\). 

        As in the proof of
        \cref{thm:MMRDynkinlabelling}, the connection to mutation-combinatorics comes by
        reducing to surfaces. Thus we consider a general elephant \(\Spec
        \overline{R}\into \Spec R\) and the reduction functor
        \(\bbF(-)=(-)\otimes_R\overline{R}\). We have the identifications of
        vector spaces
        \[
            \KK_\text{split}(\rproj \bbF\Lambda)\otimes\bbR \quad\simeq\quad
            \KK_\text{split}(\rproj \Lambda)\otimes\bbR \quad\simeq\quad
            \Theta(\ul\Delta\setminus\JJ), 
        \]
        where the first is an identification of \(g\)-fans~\cite[proposition
        B]{vangarderenStabilityCDVSingularities2023} and the second arises from
        the duality between simples and indecomposable projectives. Here
        Iyama--Wemyss~\cite[theorem 5.2]{iyamaTitsConeIntersections} show that
        the \(g\)-fan of the contracted preprojective algebra \(\bbF\Lambda\),
        and hence also that of the modification algebra \(\Lambda\), coincides with the fan \(\Arr^+(\ul\Delta,\JJ)\cup
        \Arr^-(\ul\Delta,\JJ)\).

        Analogously, the \(g\)-fan of \({}_{\nu_i}\Lambda_{\nu_i}\) sits in the
        \(\nu_i\JJ\)-cone arrangement. 

        Being induced by a derived equivalence, the map
        \[ 
        (\psi_i\dual)^{-1}\colon \KK_\text{split}(\rproj
        {}_{\nu_i}\Lambda_{\nu_i})\otimes\mathbb{R} \to \KK_\text{split}(\rproj
        \Lambda)\otimes\mathbb{R}, 
        \] 
        respects silting theory and is a map of \(g\)-fans. In particular, the
        image of \(\Chamb^+_{\nu_i\JJ}\) under \((\psi_i\dual)^{-1}\) is a cone
        in the intersection arrangement of \(\Theta(\ul\Delta\setminus\JJ)\). 
        Reading off the matrix \((\psi_i\dual)^{-1}\) from \cite[lemma
        5.2]{hiranoStabilityConditions3fold2023}, we see that 
        \((\psi_i\dual)^{-1}\Chamb_{\nu_i\JJ}^+\) intersects (and hence is equal to) the cone
        \(\nu_i\Chamb_\JJ^+\), and the images of primitive integral generators
        are primitive integral.

        Evidently from \cref{prop:wallcrossingdefn} the
        map \((\varphi_i\dual){}^{-1}\) admits an identical description, so that 
        \(\psi_i=\varphi_i\) as required.
    \end{proof}
\end{lemma}

Combined with \cref{cor:mutationdecomp},
the above lemma shows that if \(\nu=\nu_{i_n}\ldots\nu_{i_1}\) is a
\(\JJ\)-path then the map on \(\KK\)-theory induced by \(\Phi_\nu\) and
\(\Psi_\nu\) agrees with the composite \(\varphi_{i_n}\circ \cdots \circ
\varphi_{i_1}\). Computation of the heart cones \(\CC(\Psi_\nu H)\) and
\(\CC(\Phi_\nu H[-1])\) is then immediate.

\begin{proposition}
    \label{thm:heartconeofmutated}
    For each \(\JJ\)-path \(\nu\), the heart \(\Psi_\nu H\) has heart cone
    \(\nu\Chamb_\JJ^+\) while the heart \(\Phi_\nu H[-1]\) has heart cone
    \(\nu\Chamb_\JJ^-\). 
\end{proposition}

Since the complement of \(\bigcup_\nu \nu\Chamb_\JJ^\pm\) has positive codimension
in \(\Theta(\ul\Delta\setminus\JJ)\), we deduce from the above proposition and
\cref{thm:heartfandefn} that every algebraic heart in \(\tilt(H)\) must be of
the form \(\Psi_\nu H\) or \(\Phi_\nu H[-1]\) for some \(\nu\).

\paragraph{Faces of algebraic heart cones} We show that non-maximal cones in
\(\Arr^\pm (\ul\Delta\setminus\JJ)\) cannot arise as heart cones. This can be
accomplished by observing that the corresponding subcategories of semistable
modules are module categories over \(\tau\)-tilting finite
algebras~\cite{vangarderenStabilityCDVSingularities2023}, we take a
more direct approach instead.

\begin{theorem}\label{cor:heartconesinsiltingfan}
    A non-zero cone \(\sigma \in\Arr^\pm(\ul\Delta,\JJ)\) is the heart cone of some
    \(K\in\tilt(H)\) if and only if \(\sigma\) is full-dimensional. In this
    case, \((K,\sigma)\) are of the form \((\Psi_\nu H, \nu\Chamb_\JJ^+)\) or
    \((\Phi_\nu H[-1], \nu\Chamb_\JJ^-)\) for some positive \(\JJ\)-path \(\nu\). 
    \begin{proof}
        If \(\sigma\) is full dimensional it must be of the form
        \(\nu\Chamb_\JJ^\pm\), so it is a heart cone of some \(K\in \tilt(H)\)
        determined uniquely by \cref{thm:heartconeofmutated}.

        To show the converse, suppose for the sake of contradiction that there
        is a heart \(K\) in \(\tilt(H)\) whose heart cone \(\sigma=\CC K\) is non-zero but non-maximal in (without loss of generality) \(\Arr^+(\ul\Delta,\JJ)\). Consider the sets
        \begin{gather}
            \label{eqn:Kcomparable}
            \{H'\in \tilt(H)\;|\; H' \text{ algebraic},\quad \sigma\subset \CC
            (H'),\quad H'\geq K \} \\ 
            \label{eqn:Kcomparable2}
            \{H'\in \tilt(H)\;|\; H' \text{ algebraic},\quad \sigma\subset \CC
            (H'),\quad H'\leq K \}
        \end{gather}
        of algebraic hearts which are comparable with \(K\) and have \(\sigma\)
        in their heart cone. Since the \(\JJ\)-cone arrangement is locally finite away from the
        origin, these are finite.

        Replacing \((K,\sigma)\) by another pair
        if necessary, we may assume one of the sets~\eqref{eqn:Kcomparable}
        or~\eqref{eqn:Kcomparable2}, without loss of generality the former, is
        non-empty. Indeed if not, then in particular the heart
        \(K'=H_\tt(\sigma)\) is non-algebraic so that \(\sigma'=\CC(K')\) is a
        non-maximal cone containing \(\sigma\), and for each chamber
        \(\nu\Chamb_\JJ^+=\CC(\Psi_\nu H)\) surrounding \(\sigma'\) we evidently
        have \(K'=H_\tt(\sigma)\leq \Psi_\nu H\).
        
        Thus there is a minimal algebraic heart \(H'\) satisfying \(H'\geq K\) and \(\sigma\subset\CC(H')\).
        Examining the inequalities \[H\;\geq\;  H' \; \geq \; K \; \geq \;
        H[-1] \;\geq\; H'\] 
        shows that \(K\) is a tilt of \(H'\). Since \(K\neq H'\) the
        corresponding torsion class \(K[1]\cap H'\) is non-zero and in
        particular contains a simple object \(s\) of \(H'\). But then the simple
        tilt \[H''=\{h\in H'\;|\; \Hom(s, h)=0\}\ast \langle s \rangle [-1]\]
        satisfies \(H'>H''\geq K\) (in particular \(H''\in \tilt(H)\) and
        \(\sigma \subset \CC(H'')\)) and moreover is algebraic~\cite[lemma
        5.5]{hiranoStabilityConditions3fold2023}---thus contradicting the
        minimality of \(H'\).
    \end{proof}
\end{theorem}

We remark that the key ingredient in the above proof, namely covering relations
arising from simple tilts, are addressed in greater detail in
\cref{sec:covlimiting} below and in particular we identify explicitly all simple
tilts of \(\Psi_\nu H\) and \(\Phi_\nu H[-1]\) as arising from simple mutation.

\subsection{The subfan of movable cones}\label{subsec:movablecones} 
Setting aside momentarily the question
of whether the geometric and semi-geometric hearts of
\cref{subsec:geometrichearts,subsec:partialperverse} are intermediate with
respect to \(\zeroper\XZ\), we show that each non-zero cone in the finite
subfan \(\Arr(\Delta,\JJ)\subset \Arr(\ul\Delta,\JJ)\) is the heart cone of some
\(\Psi_\nu \zeroper\pair{\nu X}{Y}\) arising from a sequence of flops
\(X\dashrightarrow \nu X\) and crepant contraction \(\nu X\to Y\).

We fix, as in \cref{subsec:algebraicheartcones}, the identification \(\KK
X\simeq\mathfrak{h}(\ul\Delta\setminus \JJ)\) given by writing
\(\alpha_0= \left[\omega_{\ul C}[1]\right]\) and \(\alpha_i = [\OO_{C_i}(-1)]\)
for \(i\in \Delta\setminus \JJ\). For birational models \(W=\nu X\), there is
the analogous identification \(\KK W\simeq \mathfrak{h}(\ul\Delta\setminus\nu
\JJ)\).

\paragraph{Divisors and root lattices} By~\cite[lemma
3.4.3]{vandenberghThreedimensionalFlopsNoncommutative2004} the inclusion of the
reduced exceptional fiber \(C\into X\) induces an isomorphism of Picard groups,
in particular \(\Pic X\) is a free Abelian group of rank
\(|\Delta\setminus\JJ|\) that is naturally dual to the group \(\CycleGroup(X)\) of \(\pi\)-relative \(1\)-cycles via the perfect pairing
\begin{align}
    \label{eqn:lbcurveintersection}
    \Pic X\;\otimes_\bbZ\; \CycleGroup(X) \rightarrow \bbZ, \qquad 
    (\LL, C_i)\mapsto \deg(\LL|_{C_i}).
\end{align}
Since objects in \(\coh X\)  have proper support of dimension \(\leq 1\),
Snapper--Kleiman intersection theory~\cite[\S
VI]{kollarRationalCurvesAlgebraic1996} extends the above map to an injection
\(\Pic X\into \Hom(\KK X,\bbZ)\) by associating a line bundle \(\LL\in \Pic X\)
to the functional \((\LL\cdot -)\colon \KK X\to \bbZ\) given on (classes of)
coherent sheaves \(\mathscr{F}\in\coh(X)\) as
\[
    [\mathscr{F}] \longmapsto \chi(\mathscr{F})-\chi(\LL\,\dual \otimes
    \mathscr{F}).
\]
Note we have \((\LL\,\cdot\, \OO_{C_i}(n))=\deg(\LL|_{C_i})\) for all \(n\in
\bbZ\) and \(i\in \Delta\setminus \JJ\), consequently \((\LL \,\cdot\, \OO_p)=0\)
for any closed point \(p\in C\).

\begin{lemma}
    \label{lem:ktheoryofskyscraper}
    For any closed point \(p\) in the exceptional fiber of \(\pi\colon
    X\to Z\), the skyscraper sheaf \(\OO_p\) has \(\KK\)-theory class
    \(\delta_{\JJ}\). Consequently the assignment \(\LL\mapsto (\LL\,\cdot\,-)\)
    extends to an isomorphism of \(\RealPic(X)=\Pic X\,\otimes_\bbZ\, \bbR\)
    with the hyperplane
    \(\{\delta_\JJ=0\}\) in \(\Theta(\ul\Delta\setminus\JJ)\).
    \begin{proof}
        By \cite[theorem 5.2.4]{karmazynQuiverGITVarieties2017}, the
        \(\KK\)-theory class of \(\OO_p\) is determined by the ranks of summands
        of \(\VV\XZ\) as \([\OO_p]=\sum \text{rk}(\mathscr{N}_i)\alpha_i\). By
        \cite[proposition 9.4]{iyamaTitsConeIntersections}, this vector is
        precisely the imaginary root \(\delta_\JJ\).

        It follows that the injection \(\RealPic(X)\to
        \Hom(\KK X, \bbR)\simeq\Theta(\ul\Delta,\JJ)\) lands in the subspace
        \(\{\delta_\JJ=0\}\), and counting dimensions shows it is an isomorphism
        onto this hyperplane.
    \end{proof}
\end{lemma}

This explicit description of the functionals allows us to readily compute heart
cones of geometric and semi-geometric hearts.

\begin{proposition}\label{prop:heartconeofpartialperverse}
    Let \(\tau\colon X\to Y\) be the crepant contraction of the collection of
    exceptional curves \(\{C_i\;|\; i\in I\}\) defined by
    \(I\subset\Delta\setminus \JJ\), and suppose \(Y\not\simeq Z\) (i.e. not all
    \(\pi\)-exceptional curves are contracted). Then the heart
    \(\zeroper\XY\subset \Dm X\) of \(\tau\)-perverse sheaves has heart cone
    \(\Chamb_{\JJ\cup I}^0\subset\Theta(\ul\Delta\setminus\JJ)\).

    \begin{proof}
        If \(\theta\) lies in the heart cone of \(K=\zeroper\XY\), then
        considering the sheaves \(\OO_{C_i}(n)\) (\(n\in \mathbb{Z}\)) on
        uncontracted curves and the sheaves \(\OO_{C_i}(-1)\) on contracted
        curves shows \(\theta\) lies in the chamber \(\Chamb_\JJ^0\). In
        particular, \(\theta\) necessarily vanishes on the class \(\delta_\JJ\)
        of skyscraper sheaves. But for a connected component \(C_J\subset C_I\)
        and a closed point \(p\in C_J\), we see from
        \cref{lem:ktheoryofskyscraper} that \(\OO_p\) admits a filtration
        by the simples of \(\zeropersupp{J}\XZ\) in which every simple appears
        at least once. Thus in \(\KK\)-theory we have
        \begin{equation}
            \label{eqn:ktheoryofpartialcontraction}
            [\OO_p] = \left[\omega_{\ul C_J}[1]\right] + \sum_{i\in
            J}m_i\,[\OO_{C_i}(-1)]
            \quad = \left[\omega_{\ul C_J}[1]\right] + \sum_{i\in
            J}m_i\,\alpha_i
        \end{equation}
        for positive integers \(m_i\). Since \(\theta\) is non-negative
        on each class appearing in the above expression and in particular
        vanishes on \([\OO_p]\), it must vanish on every class on the right-hand
        side. Repeating the argument on all connected components of \(C_I\)
        shows \(\theta(\alpha_i)=0\) whenever \(i\in I\), and thus the heart
        cone \(\CC K\) is contained in \(\Chamb_{\JJ\cup I}^0\).

        Conversely if \(\theta\) lies in \(\Chamb_{\JJ\cup I}^0\), then in
        particular \(\theta\) lies in the image of the nef cone in \(\Pic X\)
        under the map of \cref{lem:ktheoryofskyscraper}, and thus 
        \(\theta[x]\geq 0\) whenever \(x\) is a coherent sheaf. On the other
        hand if \(x\in\zeroper\XY\) is supported in \(C_I\), then we may assume
        \(x\) is supported in some connected component \(C_J\subset C_I\) so
        that \(x\) lies in \(\zeropersupp{J}\XZ\). But then
        \cref{prop:perversitycheck} and the expression
        \eqref{eqn:ktheoryofpartialcontraction} show that the \(\KK\)-theory
        class \([x]\) can be expressed as some integral linear combination
        \begin{align*}
            [x] = n_\delta\, [\OO_p] + \sum_{i\in J} n_i\, [\OO_{C_i}(-1)] \quad
                = n_\delta\,\delta_\JJ + \sum_{i\in J} n_i \,\alpha_i,
        \end{align*}
        so that \(\theta[x]=0\). Since an arbitrary complex \(x\in \zeroper\XY\)
        is an extension of objects of the above forms
        (\cref{thm:structureofzeroperxy}), it follows that \(\theta[x]\geq 0\)
        and thus \(\Chamb_{\JJ\cup I}^0\) is equal to the heart cone of
        \(\zeroper\XY\).
    \end{proof}
\end{proposition}

Likewise if \(W=\nu X\) is a birational model of \(X\), we have an injection of
\(\RealPic(W)\) into \(\Theta(\ul\Delta\setminus \nu \JJ)\), and hence
into \(\Theta(\ul\Delta\setminus\JJ)\) via
the wall-crossing isomorphisms in \cref{prop:wallcrossingdefn}. It follows that
for \(I\subsetneq \Delta\setminus \nu \JJ\) and the crepant contraction
\(\tau\colon W\to Y\) of the curves \(\{C_i\;|\; i\in I\}\) in \(W\), the heart
\(\Psi_\nu \zeroper\WY\subset \Dm X\) has heart cone
\[
    (\varphi_\nu\smashdual)^{-1}\Chamb_{\nu\JJ\cup I}^0 \quad
    = \nu\Chamb_\JJ^0 \;\cap\;
      \bigcap_{i\in I}\{\varphi_\nu \alpha_i = 0\}.
\]

\paragraph{Intermediacy of (semi-)geometric hearts} 
It is evident from the description~\eqref{eqn:zeroperdefn} of \(H=\zeroper\XZ\)
that we have \(\coh X\subset H[-1,0]\), and from \cref{thm:structureofzeroperxy}
we deduce that the same can be said of all the hearts \(\zeroper\XY\) arising from
crepant contractions \(X\to Y\).

We will now show that for any birational model \(W=\nu X\) and any crepant
contraction \(W\to Y\), the heart \(\Psi_\nu\zeroper\WY\) is intermediate with
respect to \(\zeroper\XZ\). This will in particular show that every non-zero
cone in \(\Arr(\Delta,\JJ)\) is of the form \(\CC K\) for some \(K\in
\tilt(H)\), thus (with \cref{thm:heartconeofmutated}) completing the proof of
\cref{cor:heartfanofzeroper}.

\begin{proposition}
    \label{prop:intermediacyofpartialperverse}
    For any sequence of flops \(\nu\) from \(X\) and any crepant contraction
    \(W\to Y\) involving the birational model \(W=\nu X\), the heart \(\Psi_\nu
    \zeroper\WY\subset \Dm X\) is intermediate with respect to \(\zeroper\XZ\).

    \begin{proof}
        Again, \cref{thm:structureofzeroperxy} allows us to reduce to the case
        of \(\Psi_\nu \coh(W)\). Since \(\zeroper\XZ\geq \Psi_\nu \zeroper\WZ\)
        (see \cref{subsec:algebraicheartcones}) and \(\zeroper\WZ\geq \coh(W)\),
        we have \(\zeroper\XZ\geq \Psi_\nu \coh(W)\) and it thus suffices to exhibit
        the lower bound
        \[
            \Psi_\nu \coh(W)\geq \zeroper\XZ[-1].
        \]
        Applying the equivalence \({(\Psi_\nu)}^{-1}\) and swapping \(X\) and \(W\)
        for simplicity, we equivalently show
        \[ 
            \coh(X)\geq \Phi_{\nu}\zeroper\pair{\nu X}{Z}[-1]
        \]
        for any spherical and (without loss of generality) atomic \(\JJ\)-path
        \(\nu\). 

        Noting there are only finitely many atomic spherical \(\JJ\)-paths
        (since the set \(\Chambers(\Delta,\JJ)\) is finite), we can
        further reduce to the case when \(\nu\) a maximal spherical atomic
        \(\JJ\)-path. Indeed if \(\nu\) and \(\nu_i \,\nu\) are both atomic
        spherical \(\JJ\)-paths, then we have
        \(
            \Phi_i \zeroper\pair{\nu_i W}{Z}\geq \zeroper\WZ
        \),
        and hence (by \cref{lem:mutationcomposition})
        \[
            \Phi_{\nu_i\,\nu}\zeroper\pair{\nu_i\nu X}{Z}=\Phi_\nu\Phi_i \zeroper\pair{\nu_i W}{Z} \geq \Phi_\nu\zeroper\WZ.
        \]
        Noting \(\Chambers(\Delta,\JJ)\) is the face poset of a finite
        hyperplane arrangement, any longest spherical \(\JJ\)-path \(\nu\)
        defines the chamber \(\nu\Chamb_\JJ^0=-\Chamb_\JJ^0\) that is furthest
        from the principal chamber \(\Chamb_\JJ^0\). Likewise \(\nu\Chamb_\JJ^-\) contains the ray
        \(\Chamb_{\Delta}^-\) and is the furthest
        such chamber from \(\Chamb_\JJ^-\). 

        It follows that a generic vector
        \(\theta\in\nu\Chamb_\JJ^-\) satisfies \(\theta(\alpha_i)>0\) for all
        \(i\in \Delta\setminus\JJ\), and thus (up to rescaling) can be written as
        \(\theta=\theta_0 - \alpha_0^\ast\) for some \(\theta_0\in
        \Chamb_\JJ^0\) and \(\alpha_0\) defined as
        \begin{equation}
            \label{eqn:deltastar}
            \alpha_0^\ast(\alpha_i)=\begin{cases}
                1, &i=0,\\
                0, &i\in\Delta\setminus\JJ
            \end{cases}.
        \end{equation}

        Since \(\theta\) lies in the heart cone of \(\Phi_\nu\zeroper\WZ[-1]\)
        (\cref{thm:heartconeofmutated}), by
        \cref{thm:allheartsonacone} we see that the corresponding torsion-free
        class satisfies
        \begin{align*}
            \Phi_\nu \zeroper\WZ[-1]\cap H \quad&=\quad H_\tf(\theta)\\ 
            &=\quad \left\{h\in H\;\middle\vert\;
            \theta[s]> 0
            \text{ for all non-zero sub-objects }s\into h\right\}\\
            &=\quad \left\{h\in H\;\middle\vert\;
            \theta_0[s]> \alpha_0^\ast[s]
            \text{ for all non-zero sub-objects }s\into h\right\}\\
            &\subseteq \quad \left\{h\in H\;\middle\vert\;
            \theta_0[s]> 0
            \text{ for all non-zero sub-objects }s\into h\right\}\\
            &= \quad H_\tf(\theta_0),
        \end{align*}
        but \(H_\tf(\theta_0)\) is evidently contained in the torsion-free class
        \(\coh(X)\cap H\) since \(\theta_0\) lies in \(\CC(\coh X)\). The
        desired inequality of t-structures can then be read off from the above
        containment of torsion-free classes.
    \end{proof}
\end{proposition}

\subsection{Semistable categories and numerical
hearts}\label{subsec:semistablecategories}
The sub-poset of
numerical tilts of \(H\) is by \cref{thm:allheartsonacone} the union
\(\bigcup_\sigma \left[H_\tt\sigma,H^\tt\sigma\right]\) of intervals defined by
non-zero heart cones \(\sigma\) in the heart fan of \(H\). Following
\cref{cor:heartconesinsiltingfan,prop:heartconeofpartialperverse} we thus see
that all numerical tilts of \(H=\zeroper\XZ\) can be enumerated as
\[
    \{\Psi_\nu H\;|\; \nu \text{ a }\JJ\text{-path}\}\;\cup\;
    \{\Phi_\nu H[-1]\;|\; \nu \text{ a }\JJ\text{-path}\}\;\cup
    \smashoperator[r]{\bigcup_{\sigma\in\Arr(\Delta,\JJ)\smallsetminus 0}}\;
    [H_\tt\sigma,H^\tt\sigma]
\]
and it remains to address the question of when (semi-)geometric hearts are
maximal occupants of their heart cones. 

\paragraph{Classic, mixed, and reversed geometric hearts} Considering the heart
cones \(\sigma=\nu\Chamb_\JJ^0\) associated to geometric hearts \(\Psi_\nu\coh(\nu
X)\), we show that the interval \([H_\tt\sigma,H^\tt\sigma]\) is a boolean
lattice in which the geometric heart is maximal. For this a computation of all 
\(\sigma\)-semistable objects in \(H\) is necessary, this has been accomplished
by Karmazyn \cite[\S 5.2]{karmazynQuiverGITVarieties2017} and Wemyss \cite[\S
5]{wemyssFlopsClustersHomological2018} via geometric invariant theory.

\begin{lemma}
    \label{prop:semistablesonacone}
    For a chamber \(\nu \Chamb_\JJ^0\in \Chambers(\Delta,\JJ)\) and a generic
    vector \(\theta\in \nu\Chamb_\JJ^0\), an object of \(H=\zeroper\XZ\) is
    \(\theta\)-stable if and only if it is equal to a perverse point sheaf
    \(\Psi_\nu\OO_p\) for some closed point \(p\) in \(\nu X\). In particular,
    subcategory of semistable objects associated to \(\nu\Chamb_\JJ^0\) is given
    by the extension closure of such objects, i.e.\
    \[
        H_\ss\left(\nu\Chamb_\JJ^0\right)=\left\langle \Psi_\nu \OO_p\;|\; p\in
        \nu X(\mathbb{C}) \right\rangle.
    \]
\begin{proof} 
    By \cite[theorem 2.15]{wemyssFlopsClustersHomological2018} the category
    \(\zeroPer\XZ \cong \rmod\Lambda\) can be identified with a subcategory of
    representations of a quiver with vertices \(\ul\Delta\setminus\JJ\). Viewing
    elements of \(\mathfrak{h}(\ul\Delta\setminus\JJ)\) as dimension
    vectors for this quiver, the identification is such that the dimension
    vector of the representation corresponding to \(x\in \zeroper\XZ\) is
    precisely its \(\KK\)-theory class. 

    Under the above identifications, elements of the dual vector space
    \(\Theta(\ul\Delta\setminus\JJ)\) naturally give stability parameters (in
    the sense of \cite{kingModuliRepresentationsFinite1994}) on \(\zeroPer\XZ\)
    and by \cite[theorem 5.12]{wemyssFlopsClustersHomological2018}, the flop
    functor \(\Psi_\nu:\Db (\nu X)\to \Db X\) restricts to an equivalence of the
    associated semistable categories
    \[
        \left\{
            x\in \zeroPer\pair{\nu X}{Z} \;|\; x \text{ is }
            \varphi_\nu\smashdual \theta\text{-semistable }
        \right\} \longrightarrow
        \left\{
            x\in \zeroPer\pair{X}{Z} \;|\; x \text{ is }
            \theta\text{-semistable }.
        \right\}
    \]

    It suffices to show that for a generic
    \(\theta\in\Chamb_\JJ^0\), the \(\theta\)-stable objects in \(\zeroPer\XZ\)
    are precisely all the skyscraper sheaves at closed points---in our complete
    local setting, such points automatically lie in the \(\pi\)-exceptional
    fiber and thus the skyscrapers are contained in \(\zeroper\XZ\). Given this,
    the description of \(H_\ss\left(\nu\Chamb_\JJ^0\right)\) follows by
    recalling that any semistable object admits a finite filtration by stable
    ones.

    Now if an object is \(\theta\)-semistable, it remains so upon
    small perturbations of \(\theta\) within the hyperplane
    \(\{\delta_\JJ=0\}\) and thus it must have dimension vector
    \(n\delta_\JJ\) for some integer \(n>0\). The case \(n=1\) is addressed by
    constructing the coarse moduli space \(\mathcal{M}(\theta,\delta_\JJ)\),
    which Karmazyn shows is isomorphic to the variety \(\nu X\) in a way that
    each closed point
    \(p\in \nu X\) corresponds to the `perverse point sheaf' \(\Psi_\nu\OO_p\in
    \zeroPer\XZ\). Thus the collection of \(\theta\)-stable objects in
    \(\zeroper\XZ\) that have \(\KK\)-theory class \(\delta_\JJ\) is precisely
    \(
        \left\{\OO_p\;|\; p\subset X(\mathbb{C}) \right\}
    \).

    We show the case \(n\geq 2\) does not occur by replicating the argument
    of \cite[lemma 4.12]{vangarderenDonaldsonThomasInvariantsLength2022}.
    Indeed suppose for the sake of contradiction that \(x\in\zeroper\XZ\) is
    \(\theta\)-stable and \([x]=n\delta_\JJ\) for \(n\geq 2\). Since any
    skyscraper sheaf \(\OO_p\in\zeroPer\XZ\) is \(\theta\)-stable and
    necessarily distinct from \(x\), we have \(\Hom(x,\OO_p)=0\) for all \(p\in
    X\). Combining this with the fact that \(x\) is a two-term complex of
    coherent sheaves, we see that the sheaf \(\HH^0_{\coh X}(x)\) has empty
    support and hence \(x[-1]\) is a sheaf on \(X\). Choosing a sufficiently
    ample bundle \(\LL\) on \(X\) then gives us
    \(\chi(\LL\otimes x[-1]) \geq 0\), but this is absurd since
    \(x\) is numerically equivalent to some skyscraper sheaf \(\OO_p^{\oplus n}\).
\end{proof}
\end{lemma}

The following is then immediate from \cref{lem:characterisationsofHtt} and
\cref{thm:allheartsonacone}.

\begin{proposition}
    \label{prop:whatismixedgeometric}
    Let \(\nu\) be a spherical \(\JJ\)-path, defining a birational model \(W=\nu
    X\) and a cone \(\sigma=\nu \Chamb_\JJ^0\). An intermediate heart \(K\in
    \tilt(\zeroper\XZ)\) has heart cone \(\sigma\) for some spherical
    \(\JJ\)-path \(\nu\) if and only if \(K\) is a tilt of \(\Psi_\nu \coh (W)\) in a torsion class 
    \(
        \langle \Psi_\nu \OO_p \;|\; p\in U \rangle
    \)
    associated to some subset of closed points \(U\subset W(\mathbb{C})\). The
    correspondence  
    \[
        \left\{K\in \tilt(\zeroper\XZ)\;|\; \CC K = \nu\Chamb_\JJ^0\right\} \longrightarrow
        \{U \subset \nu X(\mathbb{C})\}
    \]
    thus obtained is an anti-isomorphism of posets.
\end{proposition}

In particular \(\coh(X)\), corresponding to \(U=\varnothing\), is maximal with
heart cone \(\Chamb_\JJ^0\), while the minimal such heart \(\anticoh(X)\),
corresponding to \(U=X(\mathbb{C})\), is the category of two-term complexes
\[
    \{x\in \coh X\;|\; x \text{ pure of dimension
    }1\}\;\ast\; \{y[-1]\in \coh X[-1]\;|\; \dim\Supp(y)=0\}.
\]
Analogous categories considered in \cite[example
2.6]{broomheadHeartFanAbelian2024} were dubbed \emph{reversed geometric hearts}.

\begin{remark}
    \label{rmk:artinNoether-geometrichearts}
    It is a classical fact that \(\coh X\) is a Noetherian category. On the
    other hand, any heart \(K\in [\anticoh X, \coh X)\) is necessarily
    non-Noetherian, since \(K\) then contains some \(\OO_p[-1]\) \((p\in
    C_i\subset X)\) and hence the sequence of morphisms \[\OO_{C_i}\to
    \OO_{C_i}(1)\to \OO_{C_i}(2)\to \cdots ,\] each of which has cone \(\OO_p\in
    K[1]\), gives an infinite chain of proper surjections in \(K\). 

    Dually, the heart \(\anticoh X\) is the only Artinian heart in the interval
    \([\anticoh X, \coh X]\). Indeed if a morphism \(x\to y\) is injective in
    \(\anticoh X\) with non-zero cokernel \(y/x\), then considering cohomologies
    with respect to \(\coh X\) produces a long exact sequence
        \[
            0\to \underbrace{\HH^0(x)\to \HH^0(y)\to \HH^0(y/x)}_\text{pure
                sheaves with one-dimensional support}\!\!\longrightarrow
                \underbrace{\HH^1(x)\to \HH^1(y)\to
                \HH^1(y/x)}_\text{sheaves with zero-dimensional support}\to 0
        \]
    in \(\coh X\). In particular, either \(\HH^0 x\) has lower rank than
    \(\HH^0 y\) on some exceptional curve \(C_i\subset X\), or \(y/x\)
    has zero-dimensional support (i.e.\ \(\HH^0(y/x)=0\)) and the length
    of \(\HH^1x\) is lower than that of \(\HH^1 y\). Thus  \(\anticoh
    X\) has no infinite descending chains of inclusions \(\cdots\into x_2\into
    x_1\into x\).
\end{remark}

\subsection{Many flavours of semi-geometric hearts} We compute the semistable
categories associated to non-maximal cones \(\sigma\) in the hyperplane
\(\{\delta_\JJ=0\}\), which will subsequently furnish a description of the
interval \([H_\tt\sigma, H^\tt\sigma]\). 

\begin{lemma}\label{lem:partialperversesemisimples}
    Let \(W=\nu X\) be a birational model of \(X\), and \(\tau\colon W\to Y\) a
    crepant contraction of some (but not all) exceptional curves in \(W\).
    Writing \(K=\Psi_\nu \zeroper\WY\) for the corresponding heart in \(\Dm X\)
    and picking a generic vector \(\theta\in \CC K\), an object \(k\in K\)
    satisfies \(\theta[k]=0\) if and only if it is filtered by the simple
    objects of \(K\).

    \begin{proof}
        It suffices to address the case \(W=X\) since the rest are analogous.
        Thus supposing \(\tau\) contracts the locus \(C_I=\bigcup_{i\in I}C_i\)
        defined by \(I\subset \Delta\setminus \JJ\), we must show
        \(K=\zeroper\XY\) satisfies 
        \begin{equation}
           \label{eqn:partialperversesemisimples}
            K_\ss(\theta) =
            \zeropersupp{I}\XZ
            \;\;\oplus \bigoplus_{p\in C\setminus C_I} \big\langle \OO_p
            \big\rangle
        \end{equation}
        for \(\theta\) generic in \(\CC K=\Chamb_{\JJ\cup I}^0\). It is
        evident that \(\theta\) vanishes on every simple of \(K\), so we address
        the converse containment.

        Thus suppose \(k\in K\) satisfies \(\theta[k]=0\). By
        \cref{thm:structureofzeroperxy}, \(k\in \zeroper\XZ\) is an extension of
        the simples of \(\zeropersupp{I}\XZ\) by some coherent sheaf supported
        on the \(\tau\)-uncontracted curves, so we can reduce to the case where
        \(k\) lies in \(\coh(X)\) and is supported within the closure of
        \(C\setminus C_I\).  

        We show such \(k\) continues to satisfy \(\theta[k]=0\) if \(\theta\)
        perturbed from the wall \(\Chamb^0_{\JJ\cup I}\) into the interior of
        \(\Chamb^0_\JJ\), from whence the result follows from the observation
        \[
            \left\{ x\in \coh(X)\;\middle\vert\; \dim\Supp(x)=0
            \right\} 
            \; \overset{(\ref{prop:semistablesonacone})}{=} \;
            H_\ss(\Chamb_\JJ^0) \; \overset{(\ref{thm:allheartsonacone})}{=} \;
            \left\{x\in
            H^\tt(\Chamb_\JJ^0)\;\middle\vert\; 
                \theta'[x]=0 \text{for generic }\theta'\in \Chamb_\JJ^0\right\}.
        \]
        Thus consider the vector \(\theta'=\theta + \sum_{i\in
        I}\alpha_i^\ast\), where the vectors \((\alpha_i^\ast)_{i\in I}\) are
        defined as 
        \[
            \alpha_i^\ast(\delta_\JJ)=0, \qquad
            \alpha_i^\ast(\alpha_j)=\begin{cases}
                1,&j=i \\ 
                0,&j\in \Delta\setminus (\JJ\cup \{i\})
            \end{cases}.
        \]
        Evidently \(\theta'\) lies generically in \(\Chamb_\JJ^0\). Writing
        \(I^c=\Delta\setminus (\JJ\cup I)\) so that \(k\) is supported on
        \(C_{I^c}=\bigcup_{i\in I^c}C_i\), note that \(k\) lies in the thick
        subcategory generated by \(\zeropersupp{J}\XZ\). It follows that the class \([k]\in \KK X\) can
        be written as a linear combination of the classes \(\delta_\JJ\) and
        \(\alpha_i\) (\(i\in I^c\)), and consequently \(\theta'[k]=0\) as required.
    \end{proof}
\end{lemma}

Consider the crepant contraction \(\tau\colon X\to Y\) of exceptional curves
\(C_I=\bigcup_{i\in I}C_i\), defining the intermediate heart \(K=\zeroper\XY\)
and the category \(K_\ss(\theta)\) as in~\eqref{eqn:partialperversesemisimples}
for generic \(\theta\in \CC K\). The category \(K_\ss(\theta)\) entirely lies
within \(H=\zeroper\XZ\), and thus \cref{prop:semistablesonacone} shows \(K\) is
the maximal tilt of \(H\) with heart cone \(\Chamb_{\JJ\cup I}^0\). The minimal
such heart, which we call the \emph{reversed semi-geometric heart}
\(\antizeroper\XY\), is the Artinian heart obtained by tilting \(K\) in the
torsion class \(K_\ss(\theta)\). To spell this out, writing
\(I^c=\Delta\setminus(\JJ\cup I)\) and \(C_J=\bigcup_{i\in I^c}C_i\) we have
\begin{align*}
    \antizeroper\XY &= \left\{ x\in \coh X\;\middle\vert\; \begin{array}{l}
        \Supp(x)\subset C_{I^c}\\ 
        x\text{ pure of dim. 1}
    \end{array}\right\}\;\ast\; 
    \left(\smashoperator[r]{\bigoplus_{p\in C\setminus C_I}}\;
    \left\langle \OO_p[-1] 
    \right\rangle\right)
    \;\ast\; \zeropersupp{I}{\XZ}[-1] \\ 
    &= \big\langle
    \left\{x\in \anticoh(X)\;|\; \Supp(x)\subset C_{I^c}\right\}\; \cup \; 
    \left\{x\in \zeroper\XZ[-1]\;|\; \Supp(x)\subset C_I\right\}
    \big\rangle.
\end{align*}

The poset \(\{K\in\tilt(H)\;|\; \CC K\supset\Chamb_{\JJ\cup I}^0\}\) is
thus equal to the interval \([\antizeroper\XY,\antizeroper\XY]\). The structure
can be described explicitly---decomposing \(I=I_1\sqcup\cdots \sqcup I_k\) such
that each \(C_{J}=\bigcup_{i\in J}C_i\) (\(J=I_1,\ldots,I_k\)) is a connected
component of \(C_I\), we have the following.

\begin{lemma}
    \label{prop:partialperverseinterval}
    Let \(\tau\colon X\to Y\) be the crepant contraction of \(C_I=C_{I_1}\sqcup
    \cdots \sqcup C_{I_k}\) as above. Given \(K\in \left[\,\antizeroper\XY,
    \zeroper\XY\right]\), each restricted heart
    \[K\cap \Dm X_J = \{k\in K\;|\; \Supp(k)\subset C_J\}\qquad
    (J=I_1,\ldots,I_k)\] is a tilt of \(\zeropersupp{J}\XZ\). Further, these
    restrictions and the set \(\left\{p\in C\setminus C_I\;|\; \OO_p[-1]\in
    K\right\}\) uniquely determine \(K\), and this determines a poset
    isomorphism
    \[
        \left[\,\antizeroper\XY, \zeroper\XY\,\right]
        \quad\simeq\quad
        \smashoperator[r]{\prod_{J=I_1,...,I_k}}
        \tilt\left(\zeropersupp{J}\XZ\right) \;\times\;
        \textnormal{Bool}\left(C\setminus
        C_I\right).
    \]
\end{lemma}

Here \(\textnormal{Bool}(-)\) denotes the boolean lattice on closed points, and
the poset on the right hand side is given the product order i.e.\
\((a_i)\leq (b_i)\) if and only if \(a_i\leq b_i\) for every \(i\).

\begin{proof}
        By \cref{thm:allheartsonacone} and the discussion following it, hearts
        \(K\in \left[\antizeroper\XY,\zeroper\XY\right]\) are in
        bijection with torsion classes \(T\subset \zeroper\XY\) that satisfy
        \(T\subseteq K_\ss\left(\Chamb_{\JJ\cup I}^0\right)\). For such a \(T\),
        \cref{lem:partialperversesemisimples} induces the decomposition
        \[
            T = \smashoperator[r]{\bigoplus_{J=I_1,...,I_k}}\;T\cap
            \zeropersupp{J}\XZ \;\;\oplus\; \smashoperator[r]{\bigoplus_{p\in
            C\setminus C_I}}T\cap \big\langle \OO_p\big\rangle,
        \]
        and each summand in the above decomposition is closed under extensions
        and factors (in \(K_\ss\left(\Chamb_{\JJ\cup I}^0\right)\)) since \(T\)
        is so. It follows that each \(T_J= T\,\cap\,\zeropersupp{J}\XZ\) is a
        torsion class in \(\zeropersupp{J}\XZ\), and \(K\,\cap\,\Dm X_J\) is the
        tilt of \(\zeropersupp{J}\XZ=\zeroper\XY\,\cap\,\Dm X_J\) in
        \(T_J=T\,\cap\,\Dm X_J\). Further \(T\) is clearly determined the
        summands, i.e.\ the desired bijection follows.
\end{proof}

Note when \(C\) has two or more irreducible components, each non-zero
non-maximal cone \(\sigma\in \Arr(\Delta\setminus \JJ)\) is the heart-cone of
multiple hearts of the form \(\Psi_\nu\zeroper\pair{\nu X}{Y}\)---to be precise,
one for each birational model \(\nu X\) such that \(\sigma \subset
\nu\Chamb_\JJ^0\). Indeed each chamber \(\nu\Chamb_\JJ^0\) containing \(\sigma\) uniquely determines a subset \(I\subset \Delta\setminus\nu \JJ\) such that
\(\varphi_\nu\smashdual \sigma =\Chamb_{\nu\JJ\cup I}^0\), and the data \((\nu,
I)\) together determine a partial contraction \(\nu X \to Y\) such
that \(\CC\left( \Psi_\nu \zeroper\pair{\nu X}{Y}\right) =\sigma\).

\begin{definition}
    \label{def:sigmapositive}
    Given a non-zero cone \(\sigma\in \Arr(\Delta,\JJ)\), we say a birational
    model \(\nu X\in \Bir\XZ\) is \emph{\(\sigma\)-positive} if \(\sigma\) is a
    face of \(\nu\Chamb_\JJ^0=\CC(\Psi_\nu \coh(\nu X))\).
\end{definition}

Fixing the reference heart \(H=\zeroper\XZ\), and in particular the model \(X\),
endows the set of birational models \(\Bir\XZ\) with a partial order where
\(W\geq W'\) if and only if there is a minimal sequence of simple flops from
\(X\) to \(W'\) passing through \(W\). This order, which under the natural
bijections \(\Bir\XZ\simeq \Chambers(\Delta,\JJ) \simeq \MaxModGen R\) corresponds
to the weak order on chambers or the mutation-order on modifying modules, helps
discern when a semi-geometric heart \(\Psi_\nu \zeroper\pair{\nu X}{W}\) is the
maximal occupant of its heart cone.

\begin{theorem}
    \label{thm:nonmaximalheartcones}
    Consider a birational model \(W=\nu X\) of \(X\) and let \(\tau:W\to Y\) be
    the crepant contraction of a collection of exceptional curves \(C_I\subset
    W\). Writing \(\sigma\subset \Theta(\ul\Delta\setminus\JJ)\) for the heart
    cone of \(\Psi_\nu\zeroper\WY\), the following statements are equivalent.
    \begin{enumerate}
        \item[{\crtcrossreflabel{(1)}[item:nonmaximalcone1]}]
            The birational model \(W\) is maximal among all \(\sigma\)-positive
            elements of \(\Bir\XZ\).
        \item[{\crtcrossreflabel{(2a)}[item:nonmaximalcone2a]}]
            The heart \(\Psi_\nu\zeroper\WY\) is the supremum (in \(\tilt(H)\)) of \(\left\{\Psi_{\upsilon}\coh
            (\upsilon X)\;|\; \sigma\subset \upsilon\Chamb_\JJ^0\right\}\).
        \item[{\crtcrossreflabel{(2b)}[item:nonmaximalcone2b]}]
            The heart \(\Psi_\nu\zeroper\WY\) is equal to the numerically
            defined tilt \(H^\tt(\sigma)\).

        \item[{\crtcrossreflabel{(3)}[item:nonmaximalcone3]}]
            The reversed semi-geometric heart \(\Psi_\nu\antizeroper\WY\) is
            intermediate with respect to \(H\).
        \item[{\crtcrossreflabel{(3a)}[item:nonmaximalcone3a]}]
            The heart \(\Psi_\nu\antizeroper\WY\) is the infimum of
            \(\left\{\Psi_{\upsilon}\anticoh (\upsilon X)\;|\; \sigma\subset
            \upsilon\Chamb_\JJ^0\right\}\).
        \item[{\crtcrossreflabel{(3b)}[item:nonmaximalcone3b]}]
            The heart \(\Psi_\nu\antizeroper\WY\) is equal to the numerically
            defined tilt \(H_\tt(\sigma)\).
        \item[{\crtcrossreflabel{(4)}[item:nonmaximalcone4]}]
            For generic \(\theta\in \sigma\), an object of \(H\) is
            \(\theta\)-stable if and only if it is a simple of
            \(\Psi_\nu\zeroper\WY\).
    \end{enumerate}
    Further for any non-zero cone \(\sigma \in \Arr(\Delta,\JJ)\), there is a
    unique birational model \(W=\nu X\) and a unique partial contraction
    \(W\to Y\) for which \(\sigma=\CC(\Psi_\nu\zeroper\WY)\) and
    the above holds.
\end{theorem}

Before proving the various equivalences we sketch a proof of the final claim
i.e.\ we show that for any \(\sigma\in \Arr(\Delta,\JJ)\) there is a unique
birational model and partial contraction satisfying the said conditions. By the
preceding discussion, each \(\sigma\)-positive birational \(W=\nu X\) admits
precisely one crepant contraction \(W\to Y\) such that
\(\CC(\Psi_\nu\zeroper\WY)=\sigma\) and hence it suffices to show that exactly
one \(\sigma\)-positive \(W\in \Bir\XZ\) can be maximal. But the corresponding
statement is evident for the weak order of cones---perturbing a generic vector
in \(\sigma\) `in the direction of' \(\Chamb_\JJ^0\) picks out the maximal
chamber surrounding \(\sigma\), see for example~\cite[lemma
4.2.12]{bjornerOrientedMatroids1999}.

Returning to the setup of \cref{thm:nonmaximalheartcones}, we now establish the
equivalence of statements
\ref{item:nonmaximalcone1}--\ref{item:nonmaximalcone4}. The following
equivalences arise from the interplay between numerical tilts and semistable
categories.

\begin{proof}
    [Proof of
    \ref{item:nonmaximalcone2b}\(\iff\)\ref{item:nonmaximalcone4}\(\iff\)\ref{item:nonmaximalcone3b}]
    Considering the heart \(K=\Psi_\nu\zeroper\WY\) and a generic parameter
    \(\theta\in\sigma\), note that the statement \ref{item:nonmaximalcone4} is
    equivalent to \(H_\ss\theta=K_\ss\theta\) by
    \cref{lem:partialperversesemisimples}. But we know \(K\) is intermediate
    with respect to \(H\) and has heart cone \(\sigma\), so
    \cref{lem:characterisationsofHtt} shows that we have
    \(H_\ss\theta=K_\ss\theta\) if and only if \(K\) is equal to
    \(H^\tt\theta\), thus showing
    \ref{item:nonmaximalcone2b}\(\iff\)\ref{item:nonmaximalcone4}.a

    If \ref{item:nonmaximalcone2b} (and hence also \ref{item:nonmaximalcone4})
    is true, then it is immediate that \ref{item:nonmaximalcone3b} holds, since
    \(\Psi_\nu\antizeroper\WY\) (i.e.\ the tilt of \(K\) in \(K_\ss\theta\))
    must coincide with \(H_\tt\theta\) (i.e.\ the tilt of \(H^\tt\theta\) in
    \(H_\ss\theta\)). Conversely suppose
    \(\Psi_\nu\antizeroper\WY=H^\tt\theta\). Since \(K\) lies in
    \([H_\tt\theta,H^\tt\theta]\), we see that \(K\cap H_\tt\theta[1]\) is a
    torsion-free class in \(H_\ss\theta\) and in particular is contained in
    \(H\). But we have \[K\cap H_\tt\theta[1]=\Psi_\nu\zeroper\WY \cap
    \Psi_\nu\antizeroper\WY[1])=K_\ss\theta,\] so in particular
    \(K_\ss\theta\subset H\) and hence \cref{lem:characterisationsofHtt} shows
    that both \ref{item:nonmaximalcone2b} and \ref{item:nonmaximalcone4} hold.
\end{proof}

Likewise, the equivalences
\ref{item:nonmaximalcone2a}\(\iff\)\ref{item:nonmaximalcone2b} and
\ref{item:nonmaximalcone3a}\(\iff\)\ref{item:nonmaximalcone3b} only rely on the
geometry of heart cones. We show how to prove the former, the latter being
similar.

\begin{proof}
    [Proof of \ref{item:nonmaximalcone2a}\(\iff\)\ref{item:nonmaximalcone2b}]
        It suffices to show that \(H^\tt(\sigma)\) is the supremum of the given
        collection of geometric hearts. If \(\sigma\) is a face of
        \(\upsilon \Chamb_\JJ^0=\CC(\Psi_\upsilon \coh (\upsilon X))\), then
        we clearly have \(\Psi_\upsilon \coh(\upsilon X)\leq H^\tt(\sigma)\).

        Conversely, choose generic vectors \(\theta_\upsilon \in
        \upsilon\Chamb_\JJ^0\) for each \(\upsilon\Chamb_\JJ^0\) containing
        \(\sigma\) as a face. Since the fan \(\Arr(\Delta,\JJ)\) is induced
        from a simplicial hyperplane arrangement, there is a tuple of positive
        reals \((\lambda_\upsilon)\) such that the weighted average
        \(\theta=\sum \lambda_\upsilon \theta_\upsilon\) sits
        generically in the face \(\sigma\). By the construction of numerical
        torsion theories, we see that if \(h\) lies in each torsion class
        \(H_\text{tr}(\theta_\upsilon)=\Psi_\upsilon\coh(\upsilon X)[1]\cap H\),
        then we also have \(h\in H_\text{tr}(\theta)=H^\tt(\sigma)[1]\cap H\).
        It follows that we have \(H^\tt(\sigma) \leq \sup
        \left\{\Psi_\upsilon\coh(\upsilon X)\;|\;\sigma\subset
        \upsilon\Chamb_\JJ^0\right\}\), and hence equality holds.
\end{proof}

Note that the
implication~\ref{item:nonmaximalcone3b}\(\Rightarrow\)\ref{item:nonmaximalcone3}
is immediate from definitions, so the following establishes that the
statements~\ref{item:nonmaximalcone2a}--\ref{item:nonmaximalcone4} are all
equivalent.

\begin{proof}
    [Proof of
    \ref{item:nonmaximalcone3}\(\Rightarrow\)\ref{item:nonmaximalcone4}]
    Suppose we
    have \(\Psi_\nu\antizeroper\WY\subset H[-1,0]\). Since
    \(K=\Psi_\nu\zeroper\WY\) also lies in \(H[-1,0]\), we see
    that for generic \(\theta\in\sigma\) the category
    \[K_\ss(\theta)=\Psi_\nu\zeroper\WY \cap \Psi_\nu\antizeroper\WY[1]\] lies in
    \(H\), and hence by \cref{lem:characterisationsofHtt} and
    \cref{cor:simplesofpartialperverse} we conclude that every object in
    \(H_\ss(\theta)\) is filtered by the simples of \(K\) i.e.\
    \ref{item:nonmaximalcone4} holds.
\end{proof}

It remains to show these are equivalent to \ref{item:nonmaximalcone1}, which
asserts that the shortest \(\JJ\)-path \(\mu\) satisfying \(\sigma\subset
\mu\Chamb_\JJ^0\) gives a sequence of flops from \(X\) to \(W\) (i.e.\ \(W=\mu
X\)). Unsurprisingly, we leverage our control over compositions of atomic
mutations.

\begin{proof}
    [Proof of \ref{item:nonmaximalcone1}\(\Rightarrow\)\ref{item:nonmaximalcone2a}]
    First consider the case when the shortest such path \(\mu\) is given by the
    empty word \(\emptyset\), i.e.\ \(W=X\) and \(\sigma=\Chamb_{\JJ\cup I}^0\)
    is a face of \(\Chamb_\JJ^0\). Then for \(K=\zeroper\XY\) and
    \(\theta\in\sigma\) generic, the category \(K_\ss(\theta)\) (as computed in
    \cref{lem:partialperversesemisimples}) clearly lies in \(H\) and thus
    \cref{lem:characterisationsofHtt} shows we have \(K=H^\tt(\theta)\), i.e.\
    statement \ref{item:nonmaximalcone2b} (and hence also
    \ref{item:nonmaximalcone2a}) holds.

    Now suppose \(W,Y,\sigma\) are as in the general situation and \(\mu\) is
    the shortest spherical \(\JJ\)-path satisfying \(\sigma\subset
    \mu\Chamb_\JJ^0\), i.e.\ \(\varphi_\nu\smashdual\sigma = \Chamb_{\mu\JJ\cup
    I}\) for some \(I\subset \mu\JJ\). If \ref{item:nonmaximalcone1} holds, then
    we have \(W=\mu X\) and thus \(W\to Y\) is the crepant contraction of
    \(C_I\subset W\). The reasoning above then shows
    \[
        \zeroper\WY = \sup \;\left\{\Psi_\upsilon \coh(\upsilon W)\;|\;
        \Chamb_{\mu\JJ\cup I}^0\subset \upsilon\Chamb_{\mu\JJ}^0\right\},
    \]
    and it suffices to restrict the above expression to spherical
    \(\mu\JJ\)-path \(\upsilon\) that are atomic. But if \(\upsilon\) is
    atomic then so is \(\upsilon\mu\), since every hyperplane crossed by
    \(\upsilon\) must contain \(\sigma\) while the choice of \(\mu\) ensures
    that no hyperplane crossed by it contains \(\sigma\). Thus have
    \begin{align*}
        \Psi_\mu\zeroper\WY
        &= \sup \;\left\{\Psi_\mu \circ \Psi_\upsilon\coh(\upsilon W)\;|\;
        \Chamb_{\mu\JJ\cup I}^0\subset \upsilon\Chamb_{\mu\JJ}^0\right\} \\
        &= \sup\; \left\{\Psi_{\upsilon\mu}\coh(\upsilon\mu X)\;|\; \sigma\subset
        \upsilon\mu\Chamb_\JJ^0\right\}
    \end{align*}
    To see this is equivalent to the required statement
    \ref{item:nonmaximalcone2b}, observe that for ny spherical \(\mu \JJ\)-path
    \(\upsilon\) we have \(\Chamb^0_{\mu\JJ\cup I}\subset
    \upsilon\Chamb_{\mu\JJ}^0\) if and only if \(\sigma\subset
    \upsilon\mu\Chamb_\JJ^0\), and every maximal chamber containing \(\sigma\)
    as a face can be realised in this way since
    \(\varphi_\mu\smashdual\colon \Arr(\ul\Delta,\JJ)\to \Arr(\ul\Delta,\mu\JJ)\) is
    an isomorphism of fans.
\end{proof}

\begin{proof}
    [Proof of
    \ref{item:nonmaximalcone4}\(\Rightarrow\)\ref{item:nonmaximalcone1}]
    Let \(W'=\mu X\) be the maximal \(\sigma\)-positive birational model,
    admitting the partial contraction \(W'\to Y'\) such that
    \(\CC\left(\Psi_\mu\zeroper\pair{W'}{Y'}\right)=\sigma\). The statement
    \ref{item:nonmaximalcone1} (and hence \ref{item:nonmaximalcone4}) holds for
    this partial contraction by construction, so we see that an object of \(H\)
    is \(\theta\)-stable if and only if it is a simple object of
    \(\Psi_\mu\zeroper\pair{W'}{Y'}\).

    Choosing an atomic sequence of flops \(\upsilon\) from \(W'\) to \(W\), we
    note that the path \(\upsilon\) from \(\mu\Chamb_\JJ^0\) to
    \(\upsilon\mu\Chamb_\JJ^0=\nu\Chamb_\JJ^0\) only crosses hyperplanes
    containing \(\sigma\) while the path \(\mu\) by definition never passes
    through a chamber containing \(\sigma\). Thus the composite path
    \(\upsilon\mu\) is also atomic and \cref{cor:mutationdecomp} allows us to
    write \(\Psi_\nu=\Psi_\mu\circ\Psi_\upsilon\).

    Suppose \ref{item:nonmaximalcone4} also holds for \(W\to Y\), i.e.\ an object of
    \(H\) is \(\theta\)-stable if and only if it is a simple of
    \(\Psi_\nu\zeroper\WY\). Thus we have the equality of sets
    \[
        \left\{\Psi_\mu^{-1}h
                    \;|\; h\in H \text{ \(\theta\)-stable}\right\}
        = \left\{s
                    \;|\; s\in \zeroper\pair{W'}{Y'}\text{ simple}\right\}
        = \left\{\Psi_\upsilon t
                    \;|\; t\in \zeroper\WY\text{ simple}\right\}.
    \]
    We need to show \(W=W'\), or equivalently that
    \(\upsilon\) is the empty path. If not, it has length \(\geq 1\) so can be
    written as \(\upsilon=\nu_i\cdot\upsilon'\) for some \(\mu\JJ\)-path
    \(\upsilon'\) and \(i\in \Delta\setminus\upsilon'\mu\JJ\). In other words,
    the final curve flopped by \(\upsilon\) is \(C_i\subset \upsilon'W'\), which
    has proper transform \(C_{\iota(i)}\subset W\).

    Further the wall \(\upsilon'\Chamb_\JJ^0\cap \nu_i\upsilon'\Chamb_\JJ^0\)
    contains the cone \(\sigma\), so we must have \(\iota(i)\in I\). Thus the
    curve \(C_{\iota(i)}\subset W\) is contracted by \(\tau:W\to Y\).

    In particular the object \(t=\OO_{C_{\iota(i)}}(-1)\) is simple in
    \(\zeroper\WY\), and tracking it across the flop functors using
    \cite[lemma 4.15]{wemyssFlopsClustersHomological2018} we see that
    \begin{align*}
        \Psi_\upsilon(t)
        \;=\;\Psi_{\upsilon'}\circ \Psi_i(\OO_{C_{\iota(i)}}(-1))
        \;=\;\Psi_{\upsilon'}\left(\OO_{C_i}(-1)\right)[-1]
        \;\in\; \zeroper\pair{W'}{Z}[-2,-1].
    \end{align*}
    But all simples of \(\zeroper\pair{W'}{Y'}\) must lie in
    \(\zeroper\pair{W'}{Z}\), this furnishes the desired contradiction.
\end{proof}

This concludes our analysis of the non-maximal heart cones in
\(\HFan(\zeroper\XZ)\), and thus of all numerical tilts of \(\zeroper\XZ\). In
the next section we examine local and global properties of the partial order to
conclude that all tilts of \(\zeroper\XZ\) are numerical.

\section{Covering and limiting relations}\label{sec:covlimiting}

With the goal of showing that every heart in \(\tilt(\zeroper\XZ)\) is detected
by \(\KK\)-theory---thus proving~\cref{thmC}~\ref{item:thmC1}---we now examine
covering relations in the poset of t-structures and the dynamics of \(\Pic(\nu X)\)-actions. 

A happy consequence of the study of covering relations is the classification of
bricks in \(\zeroper\XZ\). To explain this, as well as to point out the peculiar
behaviours that (non-)widely generated torsion classes in \(\zeroper\XZ\)
exhibit (see discussion following \cref{rmk:qNqA}), we first recall standard
facts in the lattice theory of torsion classes in an algebraic category.

\subsection{The lattice theory of torsion classes} \label{subsec:posets}
Let \(\mathscr{T}\) be a triangulated category with a
t-structure whose heart \(H\subset \mathscr{T}\) is algebraic, i.e.\ is an Artinian
and Noetherian Abelian category with finitely many simple objects. We have the
poset \(\tilt(H)\) of t-structures in \(\mathscr{T}\) intermediate with respect
to \(H\), the poset \(\tors(H)\) of torsion classes in \(H\), and the poset
\(\torf(H)\) of torsion-free classes in \(H\)---these are related by the
order-preserving (or order-reversing) bijections~\eqref{eqn:posetisoms1}. 

The hypotheses on \(H\) guarantee that each poset is a \textit{complete
lattice}, i.e.\ admits arbitrary infima (meets, or greatest lower bounds) and suprema
(joins, or least upper bounds). Indeed a subcategory of an Artinian heart is a
torsion-free class if and only if it is closed under taking extensions and
sub-objects, and thus \(\torf(H)\) has infima since an arbitrary intersection of
torsion-free classes remains a torsion-free class. The existence of suprema is
dual.

\paragraph{Brick labelling} Recall that we say \(a\) \textit{covers} \(b\)
(written \(a\gtrdot b\))  in a partially ordered set if \(a>b\) and there is no
element \(c\) satisfying \(a>c>b\). The \emph{Hasse quiver} of a poset has
vertices given by the elements, and an arrow \(a\to b\) for each covering
relation \(a\gtrdot b\). The join and meet operations on lattices of torsion
classes are semidistributive so their Hasse quivers are naturally labelled by
certain indecomposable objects \cite{barnardMinimalInclusionsTorsion2019}, we
now describe the construction.

\begin{definition}
    An object \(b\in H\) is a \emph{brick} if all non-zero endomorphisms of
    \(b\) are invertible. A \emph{semibrick} \(S\subset H\) is a set of bricks that are
    pairwise orthogonal, i.e.\ each \(b\in S\) is a brick, and whenever
    \(b_1,b_2\in S\) are non-isomorphic bricks, we have \(\Hom(b_1,b_2)=0\).
\end{definition}

By Schur's lemma, every simple object in \(H\) is a brick and the collection of all
simples forms a semibrick. It is not hard to see that if \(s\in
H\) is simple, then \(T=\langle s\rangle\) is a minimal non-zero torsion
class in \(H\) (i.e.\ the relation \(0\subsetdot T\) is covering in \(\tors(H)\)),
and any torsion-free class covering \(0\) must be of this form. Likewise every
maximal proper torsion class (i.e.\ a torsion class \(U\) such that
the relation \(U\subsetdot H\) is covering) is precisely of the form
\(\orth s\) for some simple object \(s\in H\). Intermediate hearts arising from
such torsion theories are called \emph{simple tilts}, and all covering relations
arise in this way.

\begin{theorem}
    [{\cite[theorems 3.3--3.4]{demonetLatticeTheoryTorsion2023}}]
    \label{thm:bricklabelsimpletilt}
    If \(K'\lessdot K\) is a covering relation in \(\tilt(H)\), then \((K\cap
    H)\setminus K'\) contains a unique brick \(b\) which we call the
    \emph{brick-label} of the covering relation. This brick is a simple object
    of \(K\), and the corresponding simple tilt in the torsion class \(\langle b
    \rangle \in \tors(K)\) is \(K'\).

    More generally if \(K'< K\) are intermediate t-structures with respect to
    \(H\) then the set of bricks in \((K\cap H)\setminus K'\) is non-empty and
    equals the set of brick-labels of covering relations in the interval
    \([K',K]\subset \tilt(H)\).
\end{theorem}

Thus each arrow in the Hasse quiver of \(\tilt(H)\) is labelled by a unique
brick in \(H\), and each brick in \(H\) arises as the brick-label of at least
one such arrow. 

The posets \(\tors(H)\) and \(\torf(H)\) inherit this labelling. In particular
if the covering relation \(T'\supsetdot T\) in \(\tors(H)\) has brick-label
\(b\), then \(b\) is the unique brick in \(T'\setminus T\) and we have \(T'\cap
T\orth = \langle b \rangle\), \(T=(\orth b) \,\ast\, T'\), and \(T'=T\ast \langle
b \rangle\). The analogous statement holds for covering relations of
torsion-free classes.

\paragraph{Widely generated torsion theories}
We review how the lattice of torsion classes can be employed to
study the set \(\sbrick(H)\) of semibricks in \(H\). Note that if \(K\subset
H[-1,0]\) is an intermediate heart, then every simple object of \(K\) lies in a
single cohomological degree with respect to \(H\). Thus writing \(\simp(K)\) for
the set of simple objects of \(K\), we see that \(\simp(K)[1]\cap H\) is a
semibrick in \(H\), and this semibrick is contained in the torsion class
associated to \(K\). This defines a map \(\simp(-)[1]\cap H:\tilt(H)\to
\sbrick(H)\), equivalently a map \(\tors(H)\to \sbrick(H)\) which maps any
torsion class to a semibrick inside it.

On the other hand each \(S\in \sbrick(H)\) defines a torsion-free class \(F=
S\orth\), with corresponding torsion class \(T=\orth (S\orth)\) characterised by
the property of being minimal among all torsion classes containing \(S\). We say
the torsion class \(T\) in this case is \emph{generated by \(S\)}. This defines
a map \(\sbrick(H)\to \tors(H)\), equivalently a map \(\sbrick(H)\to \tilt(H)\)
which we show is a section of \(\simp(-)[1]\cap H\).

\begin{lemma}
    \label{lem:sbrick-from-tilt}
    Suppose \(H=T\ast F\) is a torsion pair such that the torsion class \(T\) is
    generated by a semibrick \(S\subseteq H\). Writing \(K=F\ast T[-1]\) for the
    corresponding tilt, then an object \(k\in K\cap H[-1]\) is simple in \(K\) if
    and only if it lies in \(S[-1]\).
    \begin{proof}
        Given \(k \in S[-1]\), suppose there is an injection
        \(k'\into k\) in \(K\), with \(k'\neq 0\). Since \(T[-1]\subset K\) is a
        torsion-free class, \(k'\) also lies in \(T[-1]\). Thus the
        object \(k'[1]\in T\) is filtered by objects in \(S\) and their factors
        in \(H\) \cite[lemma 3.1]{marksTorsionClassesWide2017}. In particular, there is some
        \(s\in S\) with a non-zero morphism \(s\to k'[1]\).

        By injectivity of \(k'\into k\), the composite map \(s[-1]\to k' \to k\)
        is non-zero. Since \(S\) is a semibrick, it follows that \(s[-1]\cong k\)
        and this is a splitting of the injection \(k'\into k\). But \(k\) is
        indecomposable (since \(k[1]\) is a brick), so the map \(k'\to k\) is an
        isomorphism. Thus \(k\) has no non-trivial sub-objects in \(K\), i.e.\
        \(k\) is simple as required.

        Conversely suppose \(k\in K\cap H[-1]\) is simple in \(K\). It follows
        that \(k[1]\) is the quotient (in \(H\)) of some \(s\in S\). In other
        words, there is an \(h\in H\) and an exact triangle \(h\to s \to
        k[1]\to h[1]\). Claim \(h=0\), so that \(k[1]=s\) lies in \(S\) as
        required.

        To prove the claim, first note that the triangle \(s[-1]\to k\to h\to s\)
        shows that \(h\) cannot be a non-zero object in \(K\). Thus if \(h\) is
        non-zero, then considering the torsion part of \(h\) shows that there is
        a non-zero composite map \(s'\to h \to s\) for some \(s'\in S\). But
        \(S\) is a semibrick, so a similar argument as above shows that the map
        \(h\to s\) is an isomorphism and \(k=0\), which is a contradiction.
    \end{proof}
\end{lemma}

It follows that assigning a semibrick \(S\subseteq H\) to the the minimal
torsion class it generates gives an injective map \(\sbrick(H)\to \tors(H)\).
The following proposition provides lattice-theoretic and homological
characterisations of the image of this map, and provides alternative ways to
read off a semibrick \(S\) from the torsion class it generates.

\begin{proposition}
    \label{prop:sbrickqAcorrespondence}
    Given a torsion pair \(H=T\ast F\) with corresponding tilt \(K=F\ast
    T[-1]\), the following are equivalent.
    \begin{enumerate}[(1)]
        \item The torsion class \(T\) is generated by a semibrick \(S\subseteq H\).
            \label{item:sbrickqA1}
        \item The interval \([0,T]\subset \tors(H)\) is \emph{coatomic}, i.e.\
            for every \(U\in [0,T)\) there is a \(U'\in [0,T)\) with
            \(U\subseteq U' \subsetdot T\).
            \label{item:sbrickqA2}
        \item In \(K\), every non-zero object has a simple sub-object.
            \label{item:sbrickqA3}
    \end{enumerate}
    If the above statements hold, then the semibrick \(S\) is uniquely determined as
    \begin{align*}
        S &= \{b\in H \;|\; b\text{ is the brick-label of a relation
          \(U\subsetdot T\) in }\tors(H)\}\\
          &= \{b\in H \;|\; b[-1]\text{ is a simple object of }K\}.
    \end{align*}
\end{proposition}

We have seen (\cref{lem:sbrick-from-tilt}) that if~\ref{item:sbrickqA1} holds,
then the semibrick \(S\) which generates \(T\) is unique and is given by
\(\simp(K)[1]\cap H\). The
equivalence~\ref{item:sbrickqA1}\(\iff\)\ref{item:sbrickqA2} is the content of~\cite[theorem 7.2]{asaiWideSubcategoriesLattices2022}, which also shows that in this case
\(T\) is generated by the set of brick labels of relations \(U\subsetdot T\).
Since this set is a semibrick (see theorem 4.2 \emph{ibid.}), it must coincide
with \(S\).

It thus remains to show the equivalence~\ref{item:sbrickqA1}\(\iff\)\ref{item:sbrickqA3}, which we now do.

\begin{proof}[Proof of
    \cref{prop:sbrickqAcorrespondence} \ref{item:sbrickqA1}\( \Rightarrow
    \)\ref{item:sbrickqA3}]
     Suppose \(T\) is generated by a semibrick \(S\), and \(k\in K\) is a
     non-zero object with no simple sub-object. In particular any sub-object
     \(k'\into k\) shares this property (i.e.\ \(k'\) has no simple sub-object),
     and there is at least one such proper non-zero sub-object \(k'\).

     Now no object of \(S\subset \simp(K)\) can map to \(k\), so \(k\) lies
     in \(T[-1]\orth\). Thus neither \(k\) nor \(k'\) lie in \(T[-1]\),
     so passing to sub-objects if necessary, we may in fact assume \(k'\) and
     \(k\) lie in \(F\). But \(F\subset K\) is a torsion class, so the
     cokernel \(k/k'\) of this morphism  also lies in \(F\). It follows that the
     map \(k'\to k\) is also a proper injection in \(H\).

     But repeating the argument produces a chain of proper injections
     \(\cdots\into k''\into k'\into k\) in \(H\), which is a contradiction since
     \(H\) is Artinian. Thus every non-zero \(k\in K\) necessarily has a
     simple sub-object.
 \end{proof}

\begin{proof}[Proof of
    \cref{prop:sbrickqAcorrespondence} \ref{item:sbrickqA3}\( \Rightarrow
    \)\ref{item:sbrickqA1}]
     Suppose~\ref{item:sbrickqA3} holds. We show that any torsion-free
     class that is larger than \(F\) must contain some
     element of the semibrick \(S=\simp(K)[1]\cap H\). It follows that \(T\)
     is the minimal torsion class containing \(S\), i.e.\ \(T\) is generated
     by \(S\) as required.

     Thus consider any torsion pair \(H=T\ast F\) such that \(F\subsetneq
     F'\). Thus \(T\cap F'\) is a non-zero torsion-free class in the Abelian
     category \(K[1]= F[1]\ast (T\cap T') \ast (T\cap F')\), in particular
     \(T\cap F'\) is closed under taking sub-objects in this category. But by
     hypothesis on \(K\) (equivalently \(K[1]\)), any non-zero object in \(T\cap
     F'\) has a simple sub-object which therefore also lies in \(T\cap F'\).
     Thus \(F'\) has non-trivial intersection with the set \(\simp(K[1])\cap H
     = S\), as required.
\end{proof}

The obvious dual statements to
\cref{lem:sbrick-from-tilt,prop:sbrickqAcorrespondence} hold. In particular if
\(H=T\ast F\) is a torsion pair, then \(F\) is generated by a semibrick \(S\) if
and only if every non-zero object in the tilt \(K=F\ast T[-1]\) has a simple
factor, and in this case the semibrick generating \(F\) is determined as
\(S=\simp(K)\cap H\).

\begin{remark}\label{rmk:qNqA}
    If \(W\subseteq H\) is a \emph{wide} subcategory (i.e.\ \(W\) is closed
    under extensions, kernels, and cokernels, and is therefore Abelian),
    then the set of simples of \(W\) is evidently a semibrick of \(H\).
    Ringel \cite[\S 1.2]{ringelRepresentationsKspeciesBimodules1976} shows that
    every wide subcategory of \(H\) is in fact the extension-closure of its
    simples, and conversely the extension closure of any semibrick is a wide
    subcategory. Thus torsion(-free) classes in \(H\) that are generated by a
    semibrick are called \emph{widely generated}~\cite{asaiSemibricks2020,asaiWideSubcategoriesLattices2022,barnardMinimalInclusionsTorsion2019,marksTorsionClassesWide2017}.
\end{remark}

By \cref{prop:sbrickqAcorrespondence} and its dual, any torsion
(torsion-free) class in \(H\) associated to an Artinian (resp.\
Noetherian) tilted heart is widely generated. On the other hand a torsion pair
\(H=T\ast F\) is such that both \(T\) and \(F\) are widely generated then it is
straightforward to show that the tilt \(K=F\ast T[-1]\) is Artinian if and only if it is
Noetherian. Thus a torsion pair \(H=T\ast F\) can, in theory, exhibit seven
behaviours with regards to whether or not \(T\) and \(F\) are widely generated,
and whether \(K=F\ast T[-1]\) is Artinian or Noetherian or neither. 

All seven behaviours can be observed in the setting of a flopping contraction
\(\pi\colon X\to Z\) and its standard algebraic heart \(H=\zeroper\XZ\). Indeed
by \cref{rmk:artinNoether-geometrichearts} we have the following.

\begin{center}
\begin{tabular}
    {|p{0.4\textwidth}
     *{2}{>{\centering\arraybackslash}p{0.1\textwidth}}
     *{2}{>{\centering\arraybackslash}p{0.1\textwidth}}|}
    \hline
    \multirowcell{2}[0ex][l]{Tilted heart}
                 & \multicolumn{2}{c}{Chain conditions}
                 & \multicolumn{2}{c|}{Widely generated}\\
                 & Artin.?
                 & Noeth.?
                 & tors.?
                 & torf.?\\
    \hline
    Any algebraic, e.g.\ \(H\)
                 &\cmark
                 &\cmark
                 &\cmark
                 &\cmark \\
    Any geometric, e.g.\ \(\coh X\)
                 &\xmark
                 &\cmark
                 &\xmark
                 &\cmark \\
    Any reversed-geometric, e.g.\ \(\anticoh X\)
                 &\cmark
                 &\xmark
                 &\cmark
                 &\xmark\\
    \hdashline
\end{tabular}
\end{center}

Now consider a heart \(K\in [\anticoh X, \coh X]\) associated to a set of
closed points \(U\subseteq X(\mathbb{C})\) as in
\cref{prop:whatismixedgeometric}. All simples of \(K\) arise as skyscrapers at
closed points, hence the simple sub-objects of any \(k\in K\) are given as
\(\{\OO_p[-1]\;|\; p\in \Supp(k)\cap U\}\) while its simple quotients are
\(\{\OO_p\;|\; p\in \Supp(k)\setminus U\}\). Thus if we additionally assume that
the \(\pi\)-exceptional fiber has \(n\geq 3\) integral components
\(C_1,C_2,\ldots,C_n\) (indexed such that \(C_1\cap C_n=\emptyset\)), then
choosing closed points \(p_i\in C_i\) (\(i=1,\ldots,n\)) allows the construction
of hearts in \(H[-1,0]\) which exhibit the remaining four behaviours.

\begin{center}
\begin{tabular}
    {|p{0.4\textwidth}
     *{2}{>{\centering\arraybackslash}p{0.1\textwidth}}
     *{2}{>{\centering\arraybackslash}p{0.1\textwidth}}|}
    \hdashline
    Tilt of \(\coh X\) in... &&&&\\
    \qquad\(\langle \OO_{p_1},\ldots,\OO_{p_n} \rangle\)
                 &\xmark
                 &\xmark
                 &\cmark
                 &\cmark \\
    \qquad\(\langle \OO_{p_1}\rangle\)
                 &\xmark
                 &\xmark
                 &\xmark
                 &\cmark \\
    \qquad \(\langle \OO_p\;|\; p\in C_1 \,\cup \,\{p_2,\ldots,p_n\}
    \rangle\)
                 &\xmark
                 &\xmark
                 &\cmark
                 &\xmark \\
    \qquad\(\langle \OO_{p}\;|\; p\in C_1\rangle\)
                 &\xmark
                 &\xmark
                 &\xmark
                 &\xmark\\
   \hline
\end{tabular}
\end{center}

\subsection{Simple tilts of length hearts} \label{subsec:simplemutations} 
Following~\cite[\S 5]{hiranoFaithfulActionsHyperplane2018}, we will show that
simple tilts of algebraic hearts in \(\tilt(\zeroper\XZ)\) arise precisely from
simple mutation functors. 

Recall that the algebraic tilts of \(\zeroper\XZ\) are
all given by \(\Psi_\nu H\) or \(\Phi_\nu H[-1]\), where we write
\(H=\rflmod({}_\nu \Lambda_\nu)\) for the standard heart in the category
\(\Dfl({}_\nu \Lambda_\nu)\) associated to a modifying \(R\)-module \(M=\nu N\).
Each indecomposable summand \(M_i\subset M\) (\(i\in
\ul\Delta\setminus\mathbb{J}(M)\)) then gives an indecomposable projective
\((\End_RM)\)-module \(P_i=\Hom_R(M,M_i)\), and hence a dual simple module
\(S_i\in H\). 

We show that the intermediate hearts \(\Phi_i H[-1]\) and
\(\Psi_i H\) are precisely the simple tilts of \(H\) in \(S_i\).

\begin{lemma}
    \label{lem:simpletorsiontheories}
    Let \(M=\nu N\) and \(H=\rflmod({}_\nu\Lambda_\nu)\) be as above. For any
    \(i\in\ul\Delta\setminus \nu\JJ\), we have \(\langle S_i\rangle=\Psi_i H[1]
    \cap H = \Phi_i H[-1]\cap H\). In other words, \(\Psi_i H\) is the tilt of
    \(H\) in the torsion class \(\langle S_i\rangle\), while \(\Phi_i H[-1]\) is
    the tilt of \(H\) in the torsion-free class \(\langle S_i\rangle\).

    \begin{proof} 
        It suffices to prove that \(\Phi_i H[-1]\cap H=\langle S_i\rangle\); the
        analogous statement for \(\Psi_i\) follows by applying the functor
        \(\Phi_{i}^{-1}=\Psi_{\iota(i)}\) (\cref{lem:mutationcomposition}) and
        noting that \(\Psi_{\iota(i)}S_i=S_{\iota(i)}[-1]\)~\cite[lemma 5.3]{hiranoFaithfulActionsHyperplane2018}. For ease of notation we state the proof for \(M=N\). 

        Write \(H=T_i\ast
        F_i\) for the torsion pair defined by \(\Phi_i H\in H[-1,0]\), noting
        the torsion and torsion-free classes are defined in terms of Hom- or
        Ext-orthogonality with respect to the tilting module
        \({}_{\nu_i}\Lambda\) as in \eqref{eqn:mutatedtorsion}.
        
        Now \({}_{\nu_i}\Lambda\) is a mutation of \(\Lambda=\bigoplus_{i\in
        \ul\Delta\setminus \JJ}P_j\) at an indecomposable summand, so there is an
        exact sequence
        \[
            0\to \Lambda \to P'\to{}_{\nu_i} \Lambda\to 0
        \]
        with \(P'\in\mathop\text{add}\{P_j\;|\;j\neq i\}\). Since
        \(\Hom(P_j,S_i)=0\) for \(j\neq i\), the exact sequence gives us
        \(\Hom({}_{\nu_i}\Lambda,S_i)=0\) showing \(S_i \in T_i\). On the other hand if \(x\in \orth\{S_i\}\) then there is a
        projective cover \(P\to x\to 0\) with \(P\in
        \mathop\text{add}\{P_j\;|\;j\neq i\}\), so noting that
        \({}_{\nu_i}\Lambda\) has projective dimension \(\leq 1\) we get an
        exact sequence \(\Ext^1({}_{\nu_i}\Lambda,P)\to
        \Ext^1({}_{\nu_i}\Lambda,x)\to 0\).  But \(P\in
        \mathop\text{add}({}_{\nu_i}\Lambda)\) and \({}_{\nu_i}\Lambda\) has
        no self-extensions, so in fact all the terms in the
        sequence must vanish and we have \(x\in T_i\). Thus we have shown
        \(\orth\{S_i\}\subseteq T_i\) and \(\langle S_i\rangle\subseteq F_i\),
        and it is clear that equality must hold.
    \end{proof}
\end{lemma}

The following is then immediate from \cref{thm:bricklabelsimpletilt} and
\cref{lem:mutationcomposition}.

\begin{proposition}
    \label{cor:algebraiccoveringrelations}
    If \(\nu\) is a \(\JJ\)-path and \(i\in \ul\Delta\setminus\nu\JJ\) is such
    that \(\nu N > \nu_i \nu N\), then the poset \(\tilt(\zeroper\XZ)\) has
    relations 
    \[
        \Psi_\nu H > \Psi_{\nu_i \nu}H, \qquad \Phi_{\nu_i\nu} H[-1] > \Phi_\nu
        H[-1].
    \] 
    These are covering and have brick labels \(\Phi_\nu S_i\) and \(\Psi_\nu S_i\) respectively.
    Further, every covering relation in \(\tilt(\zeroper\XZ)\) which involves an algebraic
    heart is of the above form.
\end{proposition}

It follows that the sub-posets of \(\tilt(H)\) given by 
\begin{align*}
    \tilt^+(H)&=\{\Psi_\nu H\phantom{[-1]}\;|\; \nu \text{ a positive \(\JJ\)-path}\},\\
    \tilt^-(H)&=\{\Phi_\nu H[-1]\;|\; \nu \text{ a positive \(\JJ\)-path}\}
\end{align*}
are respectively isomorphic and anti-isomorphic to the poset of chambers
\(\Chambers(\ul\Delta,\JJ)\), and inherit its properties outlined in
\cref{lem:atomicproperties}. Further, it is evident that any interval in
\(\tilt(H)\) whose end points lie in \(\tilt^+ (H)\) (\(\tilt^-(H)\)) must entirely lie in
\(\tilt^+(H)\) (resp.\ \(\tilt^-(H)\))---that is to say if \(K\) is a heart
satisfying \(\Psi_\nu H\geq K \geq \Psi_{\nu'}H\) for some \(\JJ\)-paths
\(\nu,\nu'\), then there is a \(\JJ\)-path \(\nu''\) such that
\(K=\Psi_{\nu''}H\).

\subsection{Dynamics of the nef monoid}
\label{subsec:nefmonoid-dynamics}
The Hasse quiver of algebraic hearts, in principle, should say
all one wants to hear about the poset \(\tilt(H)\). Lattice theoretic analyses
of affine weak orders however remain elusive in practice. Fortunately, our
understanding of the global structure of \(\tilt(H)\) is greatly aided by the
presence of the natural actions of various Picard groups---these play a role
identical to that of coweight lattices in the study of Bruhat orders.

To spell this out, note that tensoring by any line bundle \(\LL\in \Pic X\)
naturally gives an autoequivalence of \(\Dm X\). We focus on the induced action on the poset
of t-structures, where a line bundle \(\LL\in \Pic X\) thus twists a heart \(K\subset\Dm X\) to
\[
    \LL\otimes K = \left\{\LL \LTensor_X k \;|\; k\in K\right\}.
\]
In what follows we study the orbit of \(H=\zeroper\XZ\) in relation to the poset
\(\tilt(H)\), setting up notation and stating key results. Proofs will
follow (and rely on) a discussion of how these actions are observed on
\(\KK\)-theory.

Certain twists of \(H\) can be described explicitly---namely, if we write
\(\LL_i=\det(\mathscr{N}_i)\) (\(i\in \Delta\setminus\JJ\)) for elements of the
standard basis of \(\Pic X\) with indexing as in
\cref{lem:indexing-of-mmg-summands}, then we have the following.

\begin{proposition}
    \label{prop:lbtwistgeneration}
    For any \(i\in \Delta\setminus \JJ\) and \(n>0\), the following statements hold.
    \begin{enumerate}[(1)]
        \item The heart \((\LL_i\dual)^{\otimes n} \otimes H\) is intermediate
            with respect to \(H\), and is the
            tilt of \(H\) in the smallest torsion class containing
            \(\OO_{C_i}(-n-1)[1]\).
        \item The heart \(\LL_i^{\otimes n}\otimes H[-1]\) is intermediate with
            respect to \(H\), and is the tilt of \(H\) in the smallest
            torsion-free class containing \(\OO_{C_i}(n-2)\).
    \end{enumerate}
\end{proposition}

If we denote by \(\tau\colon X\to X_i\) the crepant partial contraction of
\(\bigcup_{j\in\Delta\setminus\{i\}}C_j\) (i.e.\ of all exceptional curves
\emph{except} \(C_i\)), then it is clear from the above description that we have
\[
    \begin{array}{rclclclcl}
        H\phantom{[-1]} &>& \LL_i\dual \otimes H \phantom{[-1]}&>& (\LL_i\dual)^{\otimes 2}\otimes H \phantom{[-1]}
        &>&\cdots&>& \zeroper\pair{X\phantom{{}_i}}{X_i}, \\ 
        H[-1] &<& \LL_i\otimes H[-1] &<& \LL_i^{\otimes 2}\otimes H[-1]
        &<&\cdots&<& \antizeroper\pair{X\phantom{{}_i}}{X_i}.
    \end{array}
\]

\begin{proposition}
    \label{prop:rankoneneflimit}
    For any \(i\in \Delta\setminus \JJ\), the torsion class \(\zeroper
    \pair{X\phantom{{}_i}}{X_i}[1]\cap H\in \tors(H)\) is generated by  \(\{\OO_{C_i}(n)[1]\;|\; n\leq -2\}\) while the torsion-free class \(\antizeroper \pair{X\phantom{{}_i}}{X_i}\cap H\in
    \torf(H)\) is generated by \(\{\OO_{C_i}(n)\;|\; n\geq -1\}\). In other
    words, in the poset \(\tilt(H)\) we have
    \[
        \zeroper\pair{X\phantom{{}_i}}{X_i}
        = \inf \;\left\{(\LL_i\dual)^{\otimes n}\otimes H\;|\; n\geq 0\right\},
        \qquad
        \antizeroper\pair{X\phantom{{}_i}}{X_i}
        = \sup \;\left\{\LL_i^{\otimes n}\otimes H[-1]\;|\; n\geq 0\right\}.
    \]
\end{proposition}

In fact we realise many more geometric and semi-geometric hearts as limits of
twists by line bundles on birational models. To explain this, note that if
\(W=\nu X\) is a birational model then the natural action of \(\Pic W\) on
\(\Dm W\) can be translated across the flop functor \(\Psi_\nu\) to an action on
\(\Dm X\). For brevity we continue to denote this action
\(\Pic W\circlearrowright\Dm X\) by \(\;\otimes\;\), i.e.\ for \(\LL\in \Pic
(W)\) and \(h\in \Dm X\) we write
\[
    \LL\otimes h = \Psi_\nu \left(\LL\LTensor_W \Psi_\nu^{-1}h\right).
\]

Now \(\Pic W\) sits as the lattice of integral points in the vector space
\(\RealPic(W)\), which we identify (as in \cref{subsec:movablecones})  with the
hyperplane \(\{\theta\in \Theta(\ul\Delta\setminus\JJ)\;|\;
\theta(\delta_\JJ)=0\}\). Cones in this hyperplane cut submonoids of
\(\Pic W\), in particular we can consider the various monoid actions \((\Pic W \cap
\sigma)\circlearrowright\Dm X\) determined by the choice of a birational model
\(W\in \Bir\XZ\) and a cone \(\sigma\in \Arr(\Delta,\JJ)\).

The birational model \(W\) is \(\sigma\)-positive (\cref{def:sigmapositive})
precisely when \(\Pic W\cap \sigma\) lies in the nef cone.
For such \(\sigma\) and \(W\) we can consider the \((\Pic W \cap \sigma)\)-orbit of
the reference heart \(H=\zeroper\XY\). The following result shows that the orbit
is independent of the choice of \(\sigma\)-positive birational model, lies in
\(\tilt(H)\), and limits to the maximal (semi-)geometric heart determined by
\(\sigma\).

\begin{theorem}
    \label{thm:nefintermediacy}
    Given a sequence of flops \(\nu\) from \(X\) with corresponding
    birational model \(W=\nu X\), the following statements hold in \(\Dm X\).
    \begin{enumerate}[(1)]
        \item  \label{item:nefintermediacy1}
            For any \(\LL\in \Pic W\), we have \(H\geq \LL\,\dual
            \otimes H\) in the poset of t-structures on \(\Dm X\) if and only if
            \(\LL\) is nef, if and only if \(\LL \otimes H[-1]\geq H[-1]\). When
            these conditions hold, we have
            \[H\;\geq\; \LL\,\dual \otimes H \;>\; \Psi_\nu \coh(W) \;>\;\Psi_\nu\anticoh W \;>\; \LL\otimes H[-1] \;\geq\; H[-1].\]
        \item  \label{item:nefintermediacy2}
            Let \(\sigma\in \Arr(\Delta,\JJ)\) be a a cone such that \(W\) is
            \(\sigma\)-positive. Then
            \begin{align*}
                H^\tt(\sigma) &= \inf\,\left\{\LL\,\dual\otimes H\phantom{[-1]}\;|\; \LL\in
                \Pic W \cap \sigma\right\}, \\
                H_\tt(\sigma) &= \sup\,\left\{\LL\otimes H[-1]\;|\; \LL\in
                \Pic W \cap \sigma\right\}.
            \end{align*}
        \item  \label{item:nefintermediacy4}
            Suppose \(W'=\nu' X\) is another birational model, and the line
            bundles \(\LL\in \Pic W\) and \(\LL'\in \Pic W'\) coincide as
            elements of \(\Theta(\ul\Delta\setminus \JJ)\).
            If both \(\LL\) and \(\LL'\) are nef, then
            \(\LL'\,\dual\otimes H = \LL\,\dual\otimes H\) and \(\LL'\otimes
            H[-1]=\LL\otimes H[-1]\).
        \item  \label{item:nefintermediacy3}
            Let \(K\in \tilt(H)\) be of the form \(\LL'\,\dual\otimes H\) or
            \(\LL'\otimes H[-1]\) for some line bundle \(\LL'\) on some
            birational model \(W'\). If we have an inequality \(K\leq \coh(W)\)
            or \(K\geq \anticoh(W)\), then there is an \(\LL\in \Pic W\) such
            that \(K=\LL\,\dual\otimes H\) or \(K=\LL\otimes H[-1]\).
    \end{enumerate}
\end{theorem}

\begin{remark}
    For birational models \(W\) and \(W'\) of \(X\), the line bundles \(\LL\in \Pic W\) and
    \(\LL'\in \Pic W'\) coincide as elements of
    \(\Theta(\ul\Delta\setminus\JJ)\) if and only if they are proper transforms
    of each other across the birational map \(W\dashrightarrow W'\). The
    decomposition of \(\{\delta_\JJ=0\}\subset \Theta(\ul\Delta\setminus\JJ)\)
    into chambers \(\nu\Chamb_\JJ^0\) and their faces then coincides with the
    decomposition of \(\RealPic(X)\) into nef cones of birational models, the
    so-called \emph{decomposition of the movable cone} that recovers the
    Kawamata--Koll\'ar--Mori--Reid decomposition
    \cite{kawamataCrepantBlowingUp3Dimensional1988, kollarFlops1989,
    moriThreefoldsWhoseCanonical1982, reidMinimalModelsCanonical}  when \(X\) is
    a minimal model (i.e.\ \(\mathbb{Q}\)-factorial).
\end{remark}

Given any non-zero cone
\(\sigma\in \Arr(\Delta,\JJ)\), \cref{thm:nonmaximalheartcones} gives a
birational model \(W=\nu X\) and a partial contraction \(W\to Y\) such that
\(H^\tt(\sigma)=\Psi_\nu\zeroper\WY\). Then the monoid \(\Pic W \cap \sigma\)
can be naturally identified with the nef cone of \(Y\), and the limit of the
twists \(\LL\,\dual\otimes H\) ranging over \(\LL\in \Pic W \cap \sigma\) is
precisely the heart \(\zeroper\WY\) which is geometric on the curves `detected
by' \(\Pic W \cap \sigma\). In particular the limit of all twists
\(\{\LL\,\dual\otimes H\;|\; \LL\in \Pic W\text{ nef}\}\) is the geometric
heart \(\Psi_\nu\coh W\), a `fixed point' of the action of \(\Pic W\).

\paragraph{Global order on Picard orbits} Suppose \(W=\nu X\) is a birational model of \(X\) with
exceptional curves \(C_i\) \((i\in \Delta\setminus\nu \JJ)\). The
intersection pairing with these \(1\)-cycles gives an isomorphism
\(\deg \colon  \Pic W \to \mathbb{Z}^{(\Delta\setminus\nu \JJ)}\).

We endow \(\mathbb{Z}^{(\Delta\setminus\nu\JJ)}\) with
the product order (i.e.\ for tuples \((a_i),(b_i)\in
\mathbb{Z}^{(\Delta\setminus\nu\JJ)}\) we have \((a_i)\leq (b_i)\) if and only if
\(a_i\leq b_i\) for each \(i\in\Delta\setminus\nu\JJ\)) and write
\(\ul{0},\ul{1}\in \mathbb{Z}^{(\Delta\setminus\nu\JJ)}\) for the tuples whose
entries are all \(0\) or all \(1\) respectively. Thus for example a line bundle
\(\LL\in\Pic W\) is nef if and only if \(\ul{0}\leq \deg (\LL) \), and two line
bundles \(\LL_1,\LL_2\in \Pic W\) satisfy \(\deg (\LL_1) \leq \deg (\LL_2)\) if and
only if \(\LL_2\otimes \LL_1\,\dual\) is nef.

The partial order on various \(\Pic W\)-twists and shifts of \(H=\zeroper\XZ\) can then be described as
follows.

\begin{corollary}
    \label{cor:partialorderoftwists}
    Given line bundles \(\LL_1,\LL_2\) on some birational model \(W\in\Bir\XZ\) and
    integers \(n_1,n_2\in \mathbb{Z}\), we have \(\LL_1\otimes H[n_1]\geq \LL_2\otimes
    H[n_2]\) in the poset of t-structures on \(\Dm X\) if and only if either
    \(n_1> n_2\), or \(n_1=n_2\) and \(\deg(\LL_1) \geq \deg (\LL_2)\). Thus we have an isomorphism of posets \[\{\LL\otimes H[n]\;|\; \LL\in \Pic W, \; n\in
    \mathbb{Z}\}\quad \simeq \quad \mathbb{Z}\times
    \mathbb{Z}^{(\Delta\setminus \JJ)}\] 
    where \(\bbZ\times \bbZ^{(\Delta\setminus \JJ)}\) is ordered
    lexicographically on its two factors.

    \begin{proof}
        If \(n_1=n_2\), then it is clear that \(\LL_1\otimes H[n_1]\geq \LL_2\otimes
        H[n_1]\) if and only if \(H\geq \LL_1\dual\otimes \LL_2 \otimes H\), which
        by \cref{thm:nefintermediacy} \ref{item:nefintermediacy1} occurs if and
        only if \(\deg (\LL_1) \geq \deg(\LL_2)\).

        Thus to conclude it suffices to show \(\LL_1\otimes H > \LL_2\otimes
        H[-1]\) for all \(\LL_1,\LL_2\in \Pic W\). Now we may write \(\LL_1=\LL_-\dual\otimes \LL_+\) where the line bundles \(\LL_+,\LL_-\in
        \Pic W\) are nef. \Cref{thm:nefintermediacy} \ref{item:nefintermediacy1}
        shows there are inequalities \(\LL_+\otimes H[-1]\geq H[-1]\) and
        \(\LL_-\dual \otimes H > \Psi_\nu\coh W\).  Applying the functor
        \(\LL_-\dual\otimes (-)[1]\) to the first inequality thus yields
        \(\LL_1\otimes H>\Psi_\nu\coh W\), and a similar argument gives
        \(\Psi_\nu\anticoh W > \LL_2\otimes H[-1]\).
    \end{proof}
\end{corollary}

Restricting to the interval of intermediate hearts \(\left[H[-1],H\right]\), we
see that
\[
    \left\{K\in \tilt(H)\;\middle\vert\; \begin{array}{l}K=\LL\otimes H \text{
    or }\LL\otimes H[-1]\text{ for}\\\text{some }W\in \Bir\XZ,\;\LL\in
    \Pic W\end{array}\right\}  
\] 
decomposes into a union of sub-posets isomorphic to \((\bbZ_{\geq
0})^{|\Delta\setminus \JJ|}\sqcup (\bbZ_{\leq 0})^{|\Delta\setminus \JJ|}\)
indexed over \(\Bir\XZ\). This highly regular sub-poset of \(\tilt(H)\) provides
anchors to track the remaining intermediate algebraic hearts.

\begin{theorem}
    \label{thm:nefcomparisons}
    For any algebraic heart \(K\in \tilt(H)\), there is a birational model \(W\) of \(X\) and line bundles \(\LL_1,\LL_2\in \Pic W\) with degrees 
    \(\ul{0} \;\leq\; \deg(\LL_1) \;<\; \deg (\LL_2)\;\leq\; \deg (\LL_1) + \ul{1}\)
    such that one of the inequalities below holds:
    \[
        \LL_1\dual\otimes H \;\geq\; K \;>\; \LL_2\dual \otimes H, \quad
        \text{or}\quad 
        \LL_1\otimes H[-1] \;\leq\; K \;<\; \LL_2 \otimes H[-1].   
    \]
\end{theorem}

In order to prove
\cref{prop:lbtwistgeneration,prop:rankoneneflimit,thm:nefintermediacy,thm:nefcomparisons},
we now examine the actions of Picard groups on t-structures, heart cones, and
modifying modules. This builds upon the work of Hirano--Wemyss \cite[\S
7]{hiranoStabilityConditions3fold2023} on the subject.

\paragraph{Actions on the heart fan} The action of \(\Pic X\) on the
Grothendieck group \(\KK X\) is straightforward to analyse after choosing
an appropriate basis---choosing the classes \(\delta_\JJ=[\OO_p]\) and \(\alpha_i=[\OO_{C_i}(-1)]\) (\(i\in\Delta\setminus\JJ\)), the action \(\Pic X \circlearrowright \KK X\) is given by
\begin{equation}
    \label{eqn:picactiononktheory}
    \LL_i\otimes \delta_\JJ = \delta_\JJ, \qquad
    \LL_i\otimes \alpha_j =
    \begin{cases}
        \alpha_j+\delta_\JJ, &j=i\\
        \alpha_j, &\text{otherwise}
    \end{cases}.
\end{equation}
This clearly preserves the root system, i.e.\ the transposed action on
\(\Theta(\ul\Delta\setminus \JJ)\) preserves the intersection
arrangement \(\Arr(\ul\Delta,\JJ)\) and in particular takes chambers to
chambers. This can be visualised by noting that the functionals
\(\{\alpha_i\;|\;i\in \Delta\setminus\mathfrak{I}\}\) give coordinates on the
level set \(\{\delta_\JJ=1\}\) and the action of \(\Pic X\) to the level set
restricts to translations along the lattice of integral points.

If \(W=\nu X\) is a different birational model, then flop functor \(\Psi_\nu:\Dm
W\to \Dm X\) gives an isomorphism of \(\KK\)-theory \(\KK W \to \KK X\) which
maps \(\delta_{\nu\JJ}\mapsto \delta_\JJ\). Writing \(\beta_i=\left[\Psi_\nu
\OO_{C_i}(-1)\right]\) (\(i\in \Delta\setminus\nu\JJ\)) for the classes of 
the simples of \(\Psi_\nu\zeroper\WZ\), we see that \(\{\delta_\JJ, \beta_i\;|\; i\in
\Delta\setminus\mathfrak{I}\}\) is a basis for \(\KK X\) and the action \(\Pic W
\circlearrowright \KK X\) in this basis is given by a formula analogous to
\eqref{eqn:picactiononktheory}. In particular the intersection arrangement is
again preserved, and there is an induced action on the set of chambers.

For a bounded heart \(K\subset \Dm X\) and \(\LL\in \Pic W\), the chamber
\(\LL\cdot \CC K\) is the heart cone of \(\LL\,\dual\otimes K\). The following
preliminary lemma shows that this action plays well with the partial order on
chambers.

\begin{lemma}
    \label{lem:antinef-order}
    Given a birational model \(W=\nu X\) and nef line bundles \(\LL_1,\LL_2\in\Pic
    W\), we have \( \LL_1 \Chamb_\JJ^+ \geq \LL_2 \Chamb_\JJ^+\) in the weak
    order if and only if \(\deg (\LL_1) \leq \deg (\LL_2)\).
    \begin{proof}
        If \(\LL\in\Pic W\) is nef and \(i\in\Delta\setminus\nu\JJ\), then the hyperplane
        \(\{n\delta_\JJ-\beta_i=0\}\) separates \(\Chamb_\JJ^+\) and
        \(\LL\,\Chamb_\JJ^+\) if and only if \(0 < n \leq \deg(\LL|_{C_i})\). It
        follows that if \(\LL_1,\LL_2\) are nef such that \(\LL_1
        \Chamb_\JJ^+\geq \LL_2 \Chamb_\JJ^+\), then we necessarily have \(\deg
        (\LL_1) \leq \deg (\LL_2)\).

        Conversely consider any hyperplane
        \(\{\alpha=0\}\) defined by a restricted root \(\alpha\). Since \(\alpha\) is simply
        the projection of a root in an affine Dynkin root system, we may write
        \[
            \alpha=n\delta_\JJ +\sum_{i\in \Delta\setminus \JJ}{a_i}\beta_i,
        \] 
        for some \(n\in \bbZ\) and some tuple of positive integers \(\ul{a}=(a_i)\).
        Now for generic \(\theta\in \Chamb_\JJ^+\) normalised so that
        \(\theta(\delta_\JJ)=1\), and for any \(\LL\in \Pic W\), we have 
        \[
            (\LL\,\theta)(\alpha) =\theta(\alpha)+ 
                \langle \ul{a}, \deg\LL\rangle 
        \]
        where \(\langle -,-\rangle\) denotes the standard (Euclidean) inner
        product between two tuples in \(\bbZ^{(\Delta\setminus \nu\JJ)}\).
        In particular if \(\LL\) is nef then \(\langle\ul{a},\deg \LL\rangle\geq
        0\), so the hyperplane \(\{\alpha=0\}\) separates \(\Chamb_\JJ^+\) and
        \(\LL\,\Chamb_\JJ^+\) if and only if 
        \[
            \langle \ul{a},\deg\LL\rangle \geq\quad  -\theta(\alpha)
            \quad\geq\quad 0.
        \] 
        It follows that if \(\LL_1,\LL_2\in \Pic W\) are such that
        \(\deg(\LL_2)\geq \deg(\LL_1)\geq \ul{0}\) then any hyperplane
        separating \(\Chamb_\JJ^+\) from \(\LL_1\,\Chamb_\JJ^+\) also separates
        it from \(\LL_2\,\Chamb_\JJ^+\) as required.
    \end{proof}
\end{lemma}

\paragraph{Actions on modifying modules} Given a birational model \(\pi\colon W\to
Z\) and a line bundle \(\LL\in \Pic W\), the pushforward \(\pi_\ast \LL\) gives
a reflexive \(R\)-module and this determines an injective homomorphism \(\Pic
(W) \to \mathop\text{Cl}(R)\) into the class group of
\(R\). It is clear that if \(W'\) is another birational model, then the inclusions of
\(\Pic W\) and \(\Pic W'\) into \(\mathop\text{Cl}(R)\) are compatible with the
isomorphism \(\Pic W \cong \Pic W'\) induced by taking proper transforms,
and we write \(\cl(R)\subseteq \mathop\text{Cl}(R)\) for the common image of all
such inclusions. We remark that the equality \(\cl(R)=\text{Cl}(R)\) holds if
and only if \(X\) (and hence every \(W\in \Bir\XZ\)) is a minimal model of
\(Z\), this follows from \cite[proposition 9.1]{iyamaTitsConeIntersections} and
relies on the assumption that the singularity of \(Z\) is isolated.

Iyama and Wemyss \cite[\S 9.4]{iyamaTitsConeIntersections} extensively study the
action of \(\cl(R)\) on the set \(\MaxMod R\), given by  \[L\cdot M =
(L^\ast \otimes_R M)^{\ast\ast}\] for \(L\in \cl(R)\), \(M\in \MaxMod R\), and
\((-)^\ast = \Hom(-,R)\). In particular they show that the action is compatible
with that on the intersection arrangement under the natural bijection
\(\CC\colon \MaxMod R\to \Chambers(\ul\Delta,\JJ)\) of
\cref{thm:MMRDynkinlabelling}, we translate the result into a form suitable for our
purposes (see also \cite[lemma 7.2]{hiranoStabilityConditions3fold2023}).

\begin{proposition}
    \label{prop:compatiblelinebundleaction}
    If \(\pi:W\to Z\) is a birational model of \(X\), then for every \(\LL\in
    \Pic W\) and \(M\in \MaxMod R\) we have \(\LL\cdot \CC(M) = \CC(\pi_\ast\LL
    \cdot M )\).
    \begin{proof}
        It suffices to consider \(W=X\), and to prove the statement for the line
        bundles \(\LL_i = (\det \mathscr{N}_i)\,\dual\)
        (\(i\in\Delta\setminus\JJ\)) which generate
        \(\Pic X\). Further, we may assume the basic modifying module \(N\) is
        \emph{maximal modifying}, i.e.\
        \[
            \mathop\text{add}\left(N\right)= \left\{M\in \rmod R\text{
            reflexive}\;|\; \Ext^1\left(N,
            M\right)=\Ext^1\left(M,N\right)=0\right\}.
        \]
        Indeed if not, then note that by \cite[theorem
        9.1(1)]{iyamaTitsConeIntersections} and \cite[corollary
        4.18]{iyamaMaximalModificationsAuslander2014} there is a modifying
        module \(N^c\) such that \(N\oplus N^c\) is maximal modifying. This
        choice can clearly be made such that \(N\oplus N^c\) is basic, and
        further \cite[theorem 9.5]{iyamaTitsConeIntersections} shows that this
        \(N\oplus N^c\) has at most \(|\ul\Delta|\) indecomposable summands so
        we can index the summands of \(N^c\) as \(N^c=\oplus_{i\in \II}N_i\) for
        some \(\II\subset \JJ\). Then we can work with the pair \((N\oplus N^c,
        \JJ\setminus \II)\) instead of \((N,\JJ)\) in what follows, using
        \cite[theorem 8.15]{iyamaTitsConeIntersections} to translate the results
        back to \(N\).

        Consider as in the proof of \cref{thm:MMRDynkinlabelling} a general elephant \(\Spec \ol{R}\into \Spec R\) and the
        corresponding reduction functor \(\bbF(-)=(-)\otimes_R \ol{R}\). Writing
        \(\rmodif(R)\) for the set of modifying \(R\)-modules,
        Iyama--Wemyss \cite[definition 9.7]{iyamaTitsConeIntersections} define
        a map \(\text{ind}\colon \rmodif(R)\to \Theta(\ul\Delta\setminus\JJ)\)
        as the composite
        \begin{equation}
            \rmodif(R)
            \xrightarrow{\;[\Hom_R(N,-)]\;}
            \KK_\text{split}(\rproj\Lambda)
            \xrightarrow{\quad\mathbb{F}(-)\quad}
            \KK_\text{split}(\rproj\bbF\Lambda)
            \longrightarrow
            \Theta(\ul\Delta\setminus\JJ)
        \end{equation}
        where the final map sends class of the indecomposable projective
        \(\bbF\Lambda_i=\bbF\Hom_R(N,N_i)\) to the vector \(\theta_i\in
        \Theta(\ul\Delta\setminus\JJ)\) dual to \(\alpha_i\in
        \mathfrak{h}(\ul\Delta\setminus\JJ)\). By construction, for a maximal
        modifying module \(M\) the cone \(\CC(M)\) is generated by the vectors
        \(\text{ind}(M_i)\) associated to indecomposable summands \(M_i\subset
        M\). Now \cite[proposition 9.10 (1) and theorem 9.23
        (2)]{iyamaTitsConeIntersections} show that for any \(M\in \rmodif(R)\)
        we have
        \begin{align*}
            \mathop\text{ind}(\pi_\ast \LL_i \cdot M)
            &= \mathop\text{ind} (\det N_i \otimes_R M)^{\ast\ast}\\
            &= \mathop\text{ind}M +
            \delta_\JJ(\mathop\text{ind}{M})\cdot
            (\mathop\text{ind}(\det N_i) - \mathop\text{ind} R)\\
            &= \mathop\text{ind}M +
            \delta_\JJ(\mathop\text{ind}{M})\cdot
            (\mathop\text{ind} N_i - \mathop\text{rk}N_i\cdot \mathop\text{ind}
            N_0) \\
            &= \mathop\text{ind} M \circ (\LL_i\otimes - )
        \end{align*}
        where the penultimate equality uses the identity
        \(\mathop\text{ind}N_i-\mathop\text{ind}(\det N_i)=
        (\mathop\text{rk}N_i-1)\cdot \mathop\text{ind}R\) which can be deduced,
        for example, from \cite[corollary 9.22]{iyamaTitsConeIntersections}.
        Thus
        \(
            \mathop\text{ind}(\pi_\ast \LL_i \cdot M)
            = \LL_i\cdot \mathop\text{ind}M
        \)
        for every \(M\in \rmodif(R)\), and the result follows.
    \end{proof}
\end{proposition}

In what follows we suppress the map \(\pi_\ast\) from notation,
directly considering the action \(\Pic W \circlearrowright \MaxMod R\). Given
\(M\in\MaxMod R\) and \(\LL\in\Pic W\), we write \(\LL\cdot (\End_R M) = \End_R(\LL M)\), in particular \(\LL\Lambda\) is the endomorphism algebra of \(\LL N\). 
Note there is natural isomorphism of algebras \(\varepsilon\colon \LL\cdot
\End_R M\to \End_R M\) and thus an equivalence
\[
    \varepsilon\colon  \Dfl\left(\LL\,\End_R  M\right)
    \to  \Dfl \left(\End_RM\right)
\]
which in particular identifies the standard hearts, i.e.\ \(\varepsilon(H)=H\).

\paragraph{Actions on t-structures} Given \(\LL\in \Pic W\), we examine
the intermediacy of \(\LL\,\dual \otimes H\) (with respect to \(H\)) by
comparing it with the tautologically intermediate heart associated with
\(\LL N \in \MaxMod R\), i.e.\ the image of the Brenner--Butler map
\[
    \rflmod(\LL\Lambda) \xrightarrow{\quad\RHom(\Hom_R(\LL\cdot N,
    N),-)\quad} \Dfl \Lambda.
\]
Choosing an atomic sequence of mutations \(\lambda\) from \(N\) to \(\LL\cdot
N\), the above map is simply \(\Psi_\lambda\) and the resulting intermediate
heart is \(\Psi_\lambda H\). Evidently this heart agrees with
\(\LL\,\dual\otimes H\) in \(\KK\)-theory, and we now show that they are in fact
equal when \(\LL\) is nef.

\begin{proposition}
    \label{prop:lbdiagram}
    Let \(\LL\in \Pic W\) be a nef line bundle and \(\lambda\) an atomic path
    from \(N\) to \(\LL\,N\). Then the following diagram commutes on objects.
    \begin{equation}
        \label{eqn:linebundlemutation}
        \begin{tikzcd}[column sep=huge, row sep=large]
            \Dm X \dar{\textnormal{VdB}} \arrow[rr,"\LL\,\dual\otimes (-)"]
            &   & \Dm X \dar{\textnormal{VdB}} \\
            \Dfl \Lambda \rar{\varepsilon}
            & \Dfl \;\LL\Lambda \rar{\Psi_{\lambda}}
            & \Dfl \Lambda
            \end{tikzcd}
        \end{equation}
    \begin{proof}
        First consider the case when \(\LL\) is a nef bundle on
        \(X\), i.e.\ when the sequence of flops \(\nu\) is trivial. We induct on \(\deg \LL\), noting that the statement trivially holds for
        \(\LL=\OO_X\).

        So suppose the diagram \eqref{eqn:linebundlemutation}
        commutes for a nef \(\LL\in \Pic X\), and consider the bundle \(\LL_i\otimes
        \LL\) for some \(i\in \Delta\setminus \JJ\). By \cref{lem:antinef-order,lem:atomicproperties}, we can choose positive paths \(\lambda\), \(\lambda_i\) whose composite
        \(\lambda\lambda_i\) remains atomic such that \[\LL_i\, 
        \Chamb_\JJ^+=\lambda_i \Chamb_\JJ^+,\qquad (\LL_i\otimes\LL)\,
        \Chamb_\JJ^+ = \lambda\lambda_i \Chamb_\JJ^+.\]
        Further since \(\lambda\) is a reduced path from \(\LL_i\,
        \Chamb_\JJ^+\) to \((\LL_i\otimes\LL)\, \Chamb_\JJ^+\), translating
        it along \(\LL_i\,\dual\) gives a reduced (hence atomic) path
        \(\lambda'\) such that 
        \[\LL\,\Chamb_\JJ^+ = \lambda'\Chamb_\JJ^+.\]
        By \cref{prop:compatiblelinebundleaction}, we have \(\LL_i N =
        \lambda N\), \((\LL_i\otimes \LL) N = \lambda\lambda_i N\), and \(\LL N
        = \lambda' N\), and each path appearing in these expressions is atomic.
        Then the corresponding diagram \eqref{eqn:linebundlemutation} for
        \(\LL_i\otimes\LL\) is precisely the outer circuit of the diagram below,
        and it suffices to prove that each component circuit (i)--(v) commutes.
        \[
            \begin{tikzcd}[column sep=normal, row sep=normal]
                \Dm X
                \arrow[rr,"\LL\,\dual\otimes (-)"]
                \arrow[dd,"\text{VdB}"]
                & & \Dm X
                \arrow[rr,"\LL_i\dual\otimes (-)"]
                \arrow[dd,"\text{VdB}"]
                & & \Dm X
                \arrow[dd,"\text{VdB}"]
                \\
                & \text{\small{(i)}} && \text{\small{(ii)}} &\\
                \Dfl\Lambda
                \arrow[r, "\varepsilon"]
                \arrow[rrdd, "\varepsilon"', bend right=30]
                & \Dfl\; \LL\Lambda
                \arrow[r, "\Psi_{\lambda'}"]
                \arrow[rdd, "\varepsilon"', bend right=15]
                & \Dfl \Lambda
                \arrow[r, "\varepsilon"]
                & \Dfl\, \LL_i \Lambda
                \arrow[r, "\Psi_{\lambda_i}"]
                & \Dfl\Lambda \\
                &
                \hspace{-2em}\text{\small{(iii)}}&
                \text{\small{(iv)}}&
                \hspace{2em}\text{\small{(v)}}&\\
                & & \Dfl\, (\LL_i\otimes\LL)\Lambda
                \arrow[ruu, "\Psi_{\lambda}"', bend right=15]
                \arrow[rruu, "\Psi_{\lambda\lambda_i}"', bend right=30]& &
            \end{tikzcd}
        \]
        The pentagon (i) commutes by induction hypothesis, and (ii)
        is precisely the diagram shown to commute in
        \cite[theorem 7.4]{hiranoStabilityConditions3fold2023}. The triangle (iii)
        commutes since there is a natural isomorphism of functors
        \[
            (\pi_\ast\LL_i\,\dual \otimes_R -)^\ast{}^\ast
            \circ (\pi_\ast\LL\,\dual \otimes_R -)^\ast{}^\ast
            \quad\cong\quad
            (\pi_\ast(\LL_i\otimes\LL)\,\dual \otimes_R -)^\ast{}^\ast,
        \]
        so that the isomorphism of algebras
        \(\Lambda\xrightarrow{\;\varepsilon\;} (\LL_i\otimes\LL)\Lambda\) is the
        composite \(\Lambda \xrightarrow{\;\varepsilon\;} \LL\,\Lambda
        \xrightarrow{\;\varepsilon\;} (\LL_i\otimes\LL)\Lambda\). The
        commutativity of (iv) is the content of \cite[lemma
        7.3]{hiranoStabilityConditions3fold2023}, while the triangle (v)
        commutes by \cref{cor:mutationdecomp} since \(\lambda\lambda_i\) is an
        atomic path.  Thus the whole diagram commutes as required.

        Now suppose we are in the general case, i.e.\ \(\LL\) is a nef
        bundle on the flop \(W=\nu X\) and \(\lambda\) is an atomic path from
        \(N\) to \(\LL\, N\). Write \(M=\nu N\) for the modifying \(R\)-module
        generator associated to \(W\), and choose an atomic path \(\lambda'\)
        from \(M\) to \(\LL\, M\). Further we may assume the sequence of flops
        \(\nu\) is atomic, so that translating \(\nu\) by \(\LL\) as above gives
        an atomic path \(\nu'\) from \(\LL\, N\) to \(\LL\, M\).

        We can thus construct a diagram
        \[\begin{tikzcd}[row sep=normal, column sep=normal]
            \Dm X
                \arrow[rr, "\text{flop}^{-1}"]
                \arrow[dd,"\text{VdB}"] &&
            \Dm W
                \arrow[rr, "\LL\,\dual\otimes(-)"]
                \arrow[dd,"\text{VdB}"] &&
            \Dm W
                \arrow[rr, "\text{flop}"]
                \arrow[dd,"\text{VdB}"] &&
            \Dm X
                \arrow[dd,"\text{VdB}"]\\
            &\text{\small{(i)}}&&\text{\small{(ii)}}&&\text{\small{(iii)}}&\\
            \Dfl \Lambda
                \arrow[rr, "\Psi_\nu^{-1}"]
                \arrow[ddrrr, bend right=20, "\varepsilon"] &&
            \Dfl{}_\nu\Lambda_\nu
                \arrow[r, "\varepsilon"] &
            \Dfl\;\LL\,{}_\nu\Lambda_\nu
                \arrow[r, "\Psi_{\lambda'}"]
                \arrow[dd, "\Psi_{\nu'}"] &
            \Dfl{}_\nu\Lambda_\nu
                \arrow[rr, "\Psi_\nu"] &&
            \Dfl\Lambda\\
            &&\text{\small{(iv)}}&& \text{\small{(v)}}&&\\
            &&&
            \Dfl\;\LL\,\Lambda
                \arrow[rrruu, bend right=20, "\Psi_{\lambda}"] &&&
        \end{tikzcd}\]
        where the squares (i) and (iii) commute by \cref{thm:flopmutation}, while the
        circuit (ii) commutes from the above discussion since \(\LL\) is a
        nef bundle on \(W\) (the trivial flop of \(W\)). The commutativity
        of (iv) is again the content of \cite[lemma
        7.3]{hiranoStabilityConditions3fold2023}. To show (v) commutes we use 
        \cref{cor:mutationdecomp}, observing that the composite paths
        \(\lambda'\nu\) and \(\nu'\lambda\) are both atomic---indeed any
        hyperplane crossed by the atomic path \(\nu\) must pass through the
        cones
        \[
            \Chamb_\JJ^+\cap \nu\Chamb_\JJ^+\quad\supseteq\quad \smashoperator[r]{\bigcap_{i\in
            \Delta\setminus\nu\JJ}}\;\{\beta_i=0\}
        \]
        and hence is of the form \(\left\{\sum_{i\in \Delta\setminus\nu\JJ}
        b_i\beta_i =0\right\}\) for some tuple of non-negative integers \(b_i\).
        But then any generic vector \(\theta\in \nu\Chamb_\JJ^+\) and its
        translate \(\LL\,\theta\in \lambda'\nu \Chamb_\JJ^+\) both evidently lie on
        the same side of such a hyperplane, which thereby cannot be crossed by
        the atomic path \(\lambda'\). Similar reasoning shows the hyperplanes
        crossed by \(\lambda\) and \(\nu'\) are distinct.

Thus the whole diagram commutes, and
        reading the outer circuit yields the required diagram
        \eqref{eqn:linebundlemutation}.
    \end{proof}
\end{proposition}

\paragraph{Intermediacy of twists} We now prove parts
\ref{item:nefintermediacy1}, \ref{item:nefintermediacy4},
\ref{item:nefintermediacy2}, and \ref{item:nefintermediacy3} of
\cref{thm:nefintermediacy} in that order, noting that each successive proof
relies on the previous ones.

\begin{proof}[Proof of \cref{thm:nefintermediacy} \ref{item:nefintermediacy1}]
    If \(\LL\in \Pic W\) is nef, then \cref{prop:lbdiagram} shows that the heart
    \(\LL\,\dual \otimes H\) coincides with \(\Psi_\lambda H\) for some atomic
    path \(\lambda\). In particular, it is
    intermediate with respect to \(H\).

    On the other hand given any \(\LL\in \Pic W\) such that \(H\geq
    \LL\,\dual\otimes H\), we may express \(\LL\) as a
    difference of nef bundles \(\LL=\LL_-\dual\otimes\LL_+\) to obtain the
    inequality \(\LL_-\dual\otimes H \geq \LL_+\dual\otimes H\) in \(\tilt(H)\).
    Noting \(\LL_-\dual\otimes H = \Psi_\nu H\) and \(\LL_+\dual\otimes
    H=\Psi_{\nu'}H\) for some paths \(\nu,\nu'\), by
    \cref{cor:algebraiccoveringrelations} we thus have an inequality of chambers
    \(\LL_-\, \Chamb_\JJ^+\geq \LL_+\, \Chamb_\JJ^+\) in the weak order. Then
    \cref{lem:antinef-order} shows \(\deg (\LL_-)\leq \deg (\LL_+)\), i.e.\
    \(\LL\) is necessarily nef.

    One can likewise show that the nef-ness of \(\LL\)
    is also equivalent to the inequality of hearts \(\LL\otimes H[-1]\geq H[-1]\), by first
    inverting the diagram \eqref{eqn:linebundlemutation} to show \(\LL\otimes
    H[-1]=\Phi_{\lambda}H[-1]\) for some \(\JJ\)-path \(\lambda\) and
    then proceeding as above.

    Lastly we show \(\LL\,\dual\otimes H > \Psi_\nu\coh (W)\) when
    \(\LL\in\Pic W\) is nef. Consider the \(\KK\)-theory classes
    \(\beta_i=[\Psi_\nu \OO_{C_i}(-1)]\) for \(i\in\Delta\setminus\nu\JJ\), and the vector
    \(\beta_0^\ast\in\Theta(\ul\Delta\setminus\JJ)\) determined as
    \[
        \beta_0^\ast(\beta_i)=\begin{cases}
            1, &i=0 \\ 
            0, &i\in \Delta\setminus \nu \JJ.
        \end{cases}
    \]
    One checks that \(\beta_0^\ast(\delta_\JJ)=1\) and
    \(\beta_0^\ast(\alpha_i)=0\) for all \(i\in \Delta\setminus\JJ\), so we have
    \(\beta_0^\ast\in\Chamb_\JJ^+\) and
    thus the vector \(\LL\cdot \beta_0^\ast\) lies in \(\LL\,
    \Chamb_\JJ^+=\CC(\LL\,\dual\otimes H)\).  Further since \(\LL\) is nef, we
    see that \(\LL\cdot\beta_0^\ast = \beta_0^\ast + \theta^0\) for some
    \(\theta^0\in \nu\Chamb_\JJ^0\). It follows that
    \begin{align*}
        \LL\,\dual\otimes H[1]\cap H \;&\subseteq\; H^\tr(\LL\cdot
        \beta_0^\ast)\\
        &= \{h\in H\;|\; \LL\cdot\beta_0^\ast[f]\leq 0 \text{ for all
            factors } h\onto f\} \\
        &= \{h\in H\;|\; \theta^0[f]\leq -\beta_0^\ast[f] \text{ for all
            factors } h\onto f\} \\
        &\subseteq \{h\in H\;|\; \theta^0[f] < 0 \text{ for all
            non-zero factors } h\onto f\neq 0\} \\
        &= H_\tr(\theta_0) \\ \;&\subseteq\; \Psi_\nu\coh W[1]\cap H
    \end{align*}
    and thus \(\LL\,\dual\otimes H \geq \Psi_\nu \coh W\), with the inequality
    being necessarily strict since \(\Psi_\nu\coh W\) is not algebraic.
    The proof of the corresponding statement for \(\LL\otimes H[-1]\) is
    analogous.
\end{proof}

\begin{proof}[Proof of \cref{thm:nefintermediacy} \ref{item:nefintermediacy4}]
    Given birational models \(W=\nu X\) and \(W'\in \nu' X\), we have seen that
    the action of any \(\LL\in \Pic W\) on \(\KK X\) coincides with that of its
    proper transform \(\LL'\in \Pic W'\). Thus in particular \(\LL\,\dual \,
    \Chamb_\JJ^+ = \LL'\,\dual \, \Chamb_\JJ^+\), i.e.\ \(\CC(\LL \otimes H)=
    \CC(\LL'\otimes H)\).

    Now if \(\LL\) and \(\LL'\) are both nef, then \(\LL \otimes H\) and
    \(\LL'\otimes H\) are both intermediate with respect to \(H\) by \cref{thm:nefintermediacy}
    \ref{item:nefintermediacy1}. But
    the assignment \(\CC\colon \tilt(H)\to \HFan(H)\) restricts
    to a bijection between algebraic intermediate hearts and full dimensional
    heart cones, so \(\LL\,\dual\otimes H = \LL'\,\dual \otimes H\). Likewise,
    the hearts \(\LL\otimes H [-1]\) and \(\LL'\otimes H[-1]\) are equal.
\end{proof}

\begin{proof}[Proof of \cref{thm:nefintermediacy} \ref{item:nefintermediacy2}]
    Suppose \(\sigma\) and \(W\) are as given, and consider a line bundle
    \(\LL\in \Pic W\cap \sigma\). If \(\upsilon X\) is another
    \(\sigma\)-positive birational model, then the proper transform of \(\LL\)
    is nef on \(\upsilon X\) and thus \cref{thm:nefintermediacy}
    \ref{item:nefintermediacy1} and \ref{item:nefintermediacy4} give us
    \(\LL\,\dual\otimes H > \Psi_\upsilon \coh (\upsilon X)\). Considering all such
    inequalities together, we obtain
    \[
        \inf\,\left\{\LL\,\dual\otimes H\;|\; \LL\in \Pic W\cap \sigma\right\}
        \;\geq\;
        \sup\,\left\{\Psi_\upsilon\coh (\upsilon X) \;|\; \sigma\subset \upsilon
        \Chamb_\JJ^0\right\} \;\overset{(\ref{thm:nonmaximalheartcones})}{=}
        H^\tt(\sigma).
    \]

    Conversely suppose \(h\in H\) lies in \(\LL\,\dual\otimes H\) for every
    \(\LL\in \Pic W \cap \sigma\), i.e.\ \(H\) lies in the torsion-free class
    associated to the infimum on the left-hand side above.
    Pick a line bundle \(\LL_0\) that lies \emph{generically} in \(\Pic
    (W) \cap \sigma\), i.e.\ \(\LL_0\) does not lie in \(\Pic W \cap \sigma'\) for
    any proper face \(\sigma'\subset \sigma\). Thus \(\LL_0\) determines a generic
    vector \(\theta_0\in \sigma\). Considering the vector \(\theta=\LL_0^{\otimes n}\cdot
    \alpha_0^\ast\) where \(n > 0\) is some integer and \(\alpha_0^\ast\) is the vector given by
    \eqref{eqn:deltastar}, one calculates \[\theta=\alpha_0^\ast + n \cdot
    \theta_0.\]
    Moreover \(\theta\) clearly lies in \(\LL_0^{\otimes n}\,
    \Chamb_\JJ^+ = \CC((\LL_0\dual)^{\otimes n}\otimes H)\), and thus we have
    \(h\in H_\tf(\theta)\). In other words, \(\alpha_0^\ast[s] + n \cdot\theta_0[s]>0\)
    for every non-zero sub-object \(s\into h\). But \(n>0\) was arbitrary, so we
    must in fact have \(\theta_0[s]\geq 0\) for every sub-object \(s\into h\) --
    this is precisely the condition \(h\in H^\tf(\theta_0)\), and
    thus giving the opposite inequality
    \[\inf\,\left\{\LL\,\dual\otimes
    H\;|\; \LL\in \Pic W\cap \sigma\right\} \;\leq\; H^\tt(\theta_0) \;=\;
    H^\tt(\sigma).\]

    The statement for \(H_\tt(\sigma)\) can be proved analogously.
\end{proof}

\begin{proof}[Proof of \cref{thm:nefintermediacy} \ref{item:nefintermediacy3}]
    Suppose \(K\) is as given, in particular the heart cone \(\CC K\) is of the
    form \(\LL\,\Chamb_\JJ^+\) or \(\LL\,\Chamb_\JJ^-\) for a unique \(\LL\in
    \Pic W\)---without loss of generality the former. If this \(\LL\) were nef,
    then by \cref{thm:nefintermediacy} \ref{item:nefintermediacy1} we would have
    that \(K\) and \(\LL\,\dual\otimes H\) are elements of \(\tilt(H)\) with
    identical full-dimensional heart cones, whence they must be equal as
    required. 

    Noting \(\CC K\) lies in \(\{\delta_\JJ\geq 0\}\), we see that the inequality \(K\leq
    \anticoh(W)\) cannot occur. Therefore we must show
    that \(\LL\) is necessarily nef if \(K\geq \coh(W)\).
    By~\cref{thm:nefintermediacy} \ref{item:nefintermediacy2}, this inequality
    can be rewritten as the containment of torsion classes
    \begin{equation}
        \label{eqn:kinsup}
        K[1]\cap H \;\subseteq\;  \sup\,\{\LL_0\dual\otimes H[1]\cap H \;|\; {\LL_0\in \Pic W\text{ nef}}\}.
    \end{equation}
    Now for any two nef bundles \(\LL_1,\LL_2\in \Pic W\), from
    \cref{thm:nefintermediacy} \ref{item:nefintermediacy1} we see that both
    \(\LL_1\,\dual\otimes H\) and \(\LL_2\,\dual\otimes H\) are greater than
    \((\LL_1\otimes\LL_2)\,\dual\otimes H\) in the poset of t-structures. In
    other words, we have the containmen of torsion classes  
        \begin{equation}
            \label{eqn:lbtorsionclassdomination}
        \left(\LL_1\,\dual\otimes H[1]\cap H\right)\;\cup\;
        \left(\LL_2\,\dual\otimes H[1]\cap H\right) \quad\subseteq\quad
        (\LL_1\otimes \LL_2)\,\dual\otimes H[1]\cap H.
\end{equation} 
    It follows that the supremum in~\eqref{eqn:kinsup} is in fact a union. The
    torsion class \(K[1]\cap H\) is finitely generated (for example by the
    finite set of brick labels in the interval \([K, H]\subset \tilt(H)\)), so a
    sufficiently ample \(\LL_0\in \Pic W\) determines a torsion class
    \(\LL_0\otimes H[1]\cap H\) in this union that contains all the generators
    of \(K[1]\cap H\), and hence contains \(K[1]\cap H\).

    Thus we have \(H\geq K \geq \LL_0\dual \otimes H\), and consequently (by \cref{cor:algebraiccoveringrelations}) the
    inequality of heart cones \[\Chamb_\JJ^+\geq \LL\,\Chamb_\JJ^+\geq \LL_0\,\Chamb_\JJ^+\]
    for some for some nef \(\LL_0\in
    \Pic W\). The nef-ness of \(\LL\) is then evident, e.g.\ by considering
    the hyperplanes \(\left\{[{\Psi_\nu\OO_{C_i}(n)}]=0\right\}\) defined by curves \(C_i\subset
    W\) and \(n\in \bbZ\).
\end{proof}

The proofs of
\cref{prop:lbtwistgeneration,prop:rankoneneflimit,thm:nefcomparisons}
now follow.

\begin{proof}[Proof of \cref{prop:lbtwistgeneration}]
    Given the line bundle \(\LL=\LL_i^{\otimes n}\) (\(i\in \Delta\setminus
    \JJ\), \(n > 0\)), \cref{thm:nefintermediacy} \ref{item:nefintermediacy1}
    shows that \(\LL\,\dual\otimes H\) is intermediate with respect to \(H\). Further since it
    is an Artinian tilt, \cref{prop:sbrickqAcorrespondence} shows that the
    corresponding torsion class \(\LL\,\dual\otimes H[1]\cap H\) is generated by the
    semibrick
    \begin{align*}
        S
        &=\{b\in H\;|\; b[-1]\text{ is a simple object of }\LL\,\dual \otimes H\}\\
        &=\{b\in H\;|\; b\text{ is the brick label of a covering relation
        }K\gtrdot \LL\,\dual\otimes H \text{ in }\tilt(H)\}.
    \end{align*}
    In particular, the size of \(S\) is equal to the number of hearts in
    \(\tilt(H)\) covered by \(\LL\,\dual\otimes H\). Equivalently
    (\cref{cor:algebraiccoveringrelations}), this is the number of covering
    relations of the form \(\sigma \lessdot \LL\,\Chamb_\JJ^+\) in
    \(\Chambers(\ul\Delta,\JJ)\). 

    Such covering relations in the weak order correspond to hyperplanes that
    support (i.e.\ contain a codimension one face of) \(\LL\,\Chamb_\JJ^+\) and
    separate it from \(\Chamb_\JJ^+\). But the hyperplanes supporting
    \(\LL\,\Chamb_\JJ^+\) are defined by the roots
    \[
        \LL\,\dual\otimes \alpha_i = \alpha_i -n\delta_\JJ, \quad
        \LL\,\dual\otimes \alpha_j = \alpha_j \quad (j\in \Delta\setminus
        \JJ, \; j\neq i), \quad\text{and}\quad \LL\,\dual\otimes \alpha_0 =
        \alpha_0 + n\delta_i\cdot \delta_\JJ.
    \]
    Of these only \(\alpha_i-n\delta_\JJ\) separates \(\Chamb_\JJ^+\) from
    \(\LL\,\Chamb_\JJ^+\). Thus the semibrick \(S\) has precisely one element,
    and examining the simples of \(\LL\,\dual\otimes H\) shows that this element
    must be \[\LL\,\dual\otimes \OO_{C_i}(-1)[1] = \OO_{C_i}(-n-1)[1].\]

    The statement for \(\LL\otimes H[-1]\) can be proved analogously.
\end{proof}

\begin{proof}[Proof of \cref{prop:rankoneneflimit}]
    Immediate from \cref{prop:lbtwistgeneration} and
    \cref{thm:nefintermediacy} \ref{item:nefintermediacy2}.
\end{proof}

\begin{proof}[Proof of \cref{thm:nefcomparisons}]
    We prove the statement for \(K=\Psi_\nu H\) for some positive \(\JJ\)-path
    \(\nu\), the case for \(\Phi_\nu H[-1]\) being analogous. Now given a
    positive \(\JJ\)-path \(\nu\), there is a unique spherical \(\JJ\)-path
    \(\mu\) such that
    \begin{align*}
        \nu\Chamb_\JJ^+ &\;\subseteq\; \mu\Chamb_\JJ^0 \;+\; \mathbb{R}_{\geq 0}\cdot
        \alpha_0^\ast \\
              &\;=\; \left\{\theta_0+t\alpha_0^\ast\;|\; \theta_0\in
              \nu\Chamb_\JJ^0,\; t \geq 0\right\},
    \end{align*}
    where \(\alpha_0^\ast\) is as in \eqref{eqn:deltastar}. Indeed, the half
    space \(\{\delta_\JJ\geq 0\}\) (which contains \(\Arr^+(\ul\Delta,\JJ)\))
    can be written as the union of all such regions and they intersect pairwise in
    subsets of positive codimension. In particular for each class
    \(\beta_i=[\Psi_\nu \OO_{C_i}(-1)]\) \((i\in \Delta\setminus\nu \JJ)\),
    we have \[\nu\Chamb_\JJ^+\subseteq \{\theta\;|\; \theta(\beta_i) \geq 0\}.\]
    A similar argument involving a decomposition of
    \(\{\delta_\JJ \geq 0\}\) shows that for each \(i\in \Delta\setminus\nu\JJ\)
    there is a unique integer \(d_i\) such that
    \begin{equation}
        \label{eqn:boundsonCK}
        \nu\Chamb_\JJ^+ \;\subseteq\; \left\{\theta\;|\; (d_i +
        1)\cdot \theta(\delta_\JJ)  \geq \theta(\beta_i) \geq d_i\cdot \theta(\delta_\JJ)\right\},
    \end{equation}
    and clearly we must have \(d_i\geq 0\) in this case. Then the vector
    \((d_i)\in \mathbb{Z}^{|\Delta\setminus\nu\JJ|}\) determines line bundles
    \(\LL_1,\LL_2\in \Pic W\) with degrees \((d_i)\), \((d_i+1)\) respectively
---evidently we have
    \[\ul{0} \;\leq\; \deg (\LL_1) \;<\; \deg (\LL_2) \;=\; \deg (\LL_1) +
    \ul{1}.\] 
    Moreover, by construction we have the inequalities of chambers 
    \begin{equation}
        \label{eqn:conebounds}
        \LL_1\,\Chamb_\JJ^+\;\geq\; \nu\Chamb_\JJ^+\;\geq\;
        \LL_2\,\Chamb_\JJ^+.
    \end{equation}
    To see why, suppose the hyperplane \(\{\alpha=0\}\) separates
    \(\Chamb_\JJ^+\) and \(\LL_1\Chamb_\JJ^+\), where (up to replacing
    \(\alpha\) by its negative) we may write 
    \[
        \alpha= n\delta_\JJ + \sum_{i\in \Delta\setminus \mu\JJ} a_i \beta_i 
    \]
    for some \(n\in\bbZ\) and a tuple of non-negative integers \((a_i)\). For
    generic \(\theta\in \Chamb_\JJ^+\), the only way in which \(\theta(\alpha)\) and
    \((\LL_1\theta)(\alpha)=\theta(\alpha)+\sum_ia_i d_i\,\theta(\delta_\JJ)\) can have opposite signs
    is if \(\theta(\alpha)\leq 0\) and \((\LL_1\theta)(\alpha)\geq 0\). In
    particular choosing the vector \(\theta=\alpha_0^\ast\) as in \eqref{eqn:deltastar}
    which vanishes on all \(\beta_i\) (\(i\in \Delta\setminus \mu\JJ\)), we obtain the inequality 
    \[
        n + \sum_{i\in \Delta\setminus \mu\JJ} a_i d_i \geq 0.
    \]
    It follows that for generic \(\theta\in \nu\Chamb_\JJ^+\) normalised so that
    \(\theta(\delta_\JJ)=1\), we have 
    \[
        \theta(\alpha) \;=\; n+\sum_{i\in\Delta\setminus\mu\JJ}
        a_i\theta(\beta_i)\;\geq\; n+\sum_{i\in\Delta\setminus\mu\JJ} a_i d_i \;\geq\; 0.
    \]
    In other words, the hyperplane cut by \(\alpha\) must also separate
    \(\Chamb_\JJ^+\) and \(\nu\Chamb_\JJ^+\), showing \(\LL_1\,\Chamb_\JJ^+\geq
    \nu\Chamb_\JJ^+\) as required. The inequality
    \(\nu\Chamb_\JJ^+>\LL_2\,\Chamb_\JJ^+\) is analogous.
    
    The desired inequality on \(K\) then follows by translating the inequality
    of heart cones~\eqref{eqn:conebounds} into inequalities of intermediate
    algebraic hearts via \cref{cor:algebraiccoveringrelations}.
\end{proof}

\subsection{Every torsion theory on the perverse heart is numerical}
\label{subsec:mainargument}
The above analysis of \(\tilt(H)\) paves the way for this result, showing that
the heart fan of \(H\) detects \emph{every} intermediate t-structure and thus an
arbitrary heart in \(\tilt(H)\) must be one of the algebraic, geometric, or
semi-geometric t-structures described in \cref{sec:numtors}.

\begin{theorem}
    \label{thm:everyheartisnumerical}
    Let \(K\) be an arbitrary tilt of \(H\). Then the heart cone
    \(\CC K\) is non-zero.
\end{theorem}

We prove \cref{thm:everyheartisnumerical} over the course of this subsection.
Thus, fix \(K\in \tilt(H)\) as above. If \(K\) is algebraic, then \(\CC K\)
is full-dimensional (in particular, non-zero) by
\cref{thm:heartfandefn}. So we may restrict to the case when \(K\) is
non-algebraic. By \cref{cor:heartconesinsiltingfan}, \(\CC K\) is thus a
(possibly zero) cone in \(\Arr(\Delta,\JJ)\). The following lemma is then
the key tool we exploit.

\begin{lemma}
    \label{lem:boundonbothsides}
    Given \(K\in\tilt(H)\), suppose there is a non-zero cone
    \(\sigma\in\Arr(\Delta,\JJ)\) and a \(\sigma\)-positive birational model
    \(W\in \Bir\XZ\) such that \(\LL\,\dual\otimes H \;\geq\; K \;\geq\; \LL\otimes H[-1]\) for
    every \(\LL\in \Pic W \cap \sigma\). Then the heart cone \(\CC K\) contains
    \(\sigma\), and is in particular non-zero.

    \begin{proof}
        The conditions, combined with \cref{thm:nefintermediacy}
        \ref{item:nefintermediacy2} show \(H^\tt(\sigma)\geq K\geq
        H_\tt(\sigma)\) so the result follows
        from \cref{thm:allheartsonacone}.
    \end{proof}
\end{lemma}

It is relatively straightforward to obtain the bound on one side for
\cref{lem:boundonbothsides}.

\begin{lemma}
    \label{lem:boundononeside}
    If the heart \(K\in \tilt(H)\) is not algebraic, then there is a non-zero cone
    \(\sigma\in\Arr(\Delta,\JJ)\) and a \(\sigma\)-positive birational model
    \(W\in \Bir\XZ\) such that \(\LL\,\dual\otimes H \geq K\) for every \(\LL\in \Pic W
    \cap \sigma\).

    Likewise, there is a (possibly different) non-zero cone \(\sigma'\in
    \Arr(\Delta,\JJ)\) and a \(\sigma'\)-positive birational model
    \(W'\in\Bir\XZ\) such that \(K\geq \LL'\otimes H[-1]\) for every \(\LL'\in
    \Pic W' \cap \sigma'\).

    \begin{proof}
        We prove the first statement, the second being analogous. Since
        \(K\) is non-algebraic, \cref{cor:algebraiccoveringrelations} shows that
        \(K\) is not covered by an algebraic heart. In particular if
        \(K_n\in\tilt(H)\) is algebraic and satisfies \(K_n>K\) then there is an algebraic
        heart \(K_{n+1}\) such that \(K_n\gtrdot
        K_{n+1}>K\). Starting with the tautological relation
        \(H=K_0>K\), this produces an infinite chain \(K_0 > K_1 >
        K_2 > \cdots >K\) approaching \(K\) using algebraic tilts of
        \(H\), each of the form \(\Psi_\nu H\) for some \(\JJ\)-path \(\nu\).

        By \cref{thm:nefcomparisons}, for each \(K_n\) there is a birational
        model \(W_n\) and line bundles \(\LL_n, \LL_n'\) in \(\Pic W_n\) such that
        \[\LL_n\dual\otimes H \;\geq\; K_n \;\geq\; \LL'_n\dual\otimes H, \quad
        \text{and}\quad\ul{0}\;\leq\; \deg (\LL_n) < \deg(\LL_n') \leq  \deg
        (\LL_n) + \ul{1}.\]
        But there are only finitely many birational models of \(X\), so we can
        pass to a subsequence if necessary to assume all birational
        models \(W_n\) are equal. Thus there is a \(W=\nu X\) and a sequence of
        nef line bundles \(\LL_1,\LL_2,\ldots\in \Pic W\) such that
        \(\LL_n\dual\otimes H > K\) for each \(n\).

        Claim there is some integral exceptional curve \(C_i\subset W\) such that
        \(\deg(\LL_n|_{C_i})\) attains arbitrarily large magnitude. If not, then
        the tuples \(\deg(\LL_n)\) are bounded and we can
        choose a line bundle \(\LL_\infty\in \Pic W\) such that
        \(\deg\LL_\infty\geq \deg \LL_n + \ul{1}\) for each \(n\).
        Thus \(\deg\LL_n'\leq \deg \LL_\infty\) for each \(n\), which shows
        \(\LL_n'\dual\otimes H \geq \LL_\infty\dual \otimes H\) and hence \(K_n
        \in [\LL_\infty\dual\otimes H, H]\) for each \(n\).

        Noting
        \(\LL_\infty\dual\otimes H=\Psi_\nu H\) for some \(\JJ\)-path \(\nu\),
        we see that the interval \([\Psi_\nu H, H]\subset\tilt(H)\) contains
        infinitely many elements. But we have seen
        (\cref{cor:algebraiccoveringrelations}) that there is a poset
        isomorphism 
        \[
            \tilt(H)\supset [\Psi_\nu H, H]\quad\simeq \quad
            [\nu\Chamb_\JJ^+,\Chamb_\JJ^+]\subset \Chambers(\ul\Delta,\JJ)
        \]
        and the latter interval is finite by \cref{lem:atomicproperties}---a
        contradiction.

        Therefore writing \(\LL\in \Pic W\) for the line bundle which has degree \(1\) on \(C_i\) and is trivial elsewhere, we see that for every
        \(N>0\) there is an \(n\) such that \((\LL\,\dual)^{\otimes N} \otimes H >
        \LL_n\dual\otimes H\) and thus 
        \((\LL\,\dual)^{\otimes N}\otimes H > K\) for all \(N\). 

        Choosing \(\sigma\in\Arr(\Delta,\JJ)\) to be the ray generated by (the
        vector corresponding to) \(\LL\) then yields the result.
    \end{proof}
\end{lemma}

As an immediate consequence, we see that any heart bounded by a geometric
one is numerical.

\begin{lemma}
    \label{lem:boundedbycoh}
    If \(K\in\tilt(H)\) is such that \(K\geq \Psi_\nu \coh W\) or \(\Psi_\nu
    \anticoh W \geq K\) for some birational model \(W=\nu X\), then \(\CC K\)
    is non-zero.

    \begin{proof}
        Again we may assume \(K\) is non-algebraic and lies in \([\Psi_\nu \coh
        W, H]\), so that \cref{lem:boundononeside} gives a non-zero cone
        \(\sigma\in \Arr(\Delta,\JJ)\) and a \(\sigma\)-positive birational
        model \(W'\) such that \(\LL'\,\dual\otimes H \geq K \geq \Psi_\nu\coh
        W\) for every \(\LL'\in \Pic W' \cap \sigma\). 

        But then~\cref{thm:nefintermediacy}~\ref{item:nefintermediacy3} shows that for
        each \(\LL'\in \Pic W'\cap \sigma\), we may choose a nef bundle \(\LL\in \Pic
        (W)\) (given by the proper transform of \(\LL'\)) such that \(\LL'\,\dual\otimes H = \LL\,\dual\otimes H\). In particular, \(W\) is \(\sigma\)-positive and
        we have \(\LL\,\dual\otimes H \geq K\) for every \(\LL\in \Pic W\cap
        \sigma\).  

        On the other hand we also have the inequality \(K\geq
        \Psi_\nu\coh W \geq \LL \otimes H[-1]\) for every \(\LL\in \Pic W \cap
        \sigma\), so \cref{lem:boundonbothsides} yields the conclusion.
    \end{proof}
\end{lemma}

We deduce analogous bounds on hearts \(K\in \tilt(H)\) which contain skyscraper
sheaves \(\OO_p\) (for closed points \(p\in X\)) or shifts thereof, noting that a priori each
\(\OO_p\) is only a two-term complex in \(K[0,1]\).

\begin{lemma}
    \label{lem:skyscrapertrick}
    Suppose \(C_i\subset X\) is an integral exceptional curve with a closed
    point \(p\in C_i\), and \(\LL_i\in \Pic X\) is the line bundle which has
    degree \(1\) on \(C_i\) and is trivial elsewhere. Writing \(\sigma\in
    \Arr(\Delta,\JJ)\) for the cone generated by \(\LL_i\), the following
    statements hold for any \(K\in\tilt(H)\).
    \begin{enumerate}[(1)]
        \item If \(\OO_p\in K\), then \(K\geq \LL\otimes H[-1]\) for
            every \(\LL\in\Pic X \cap \sigma\).
            \label{item:skyscraperbounds1}
        \item If \(\OO_p\in K[1]\), then \(\LL\,\dual\otimes H \geq K\) for
            every \(\LL\in\Pic X \cap \sigma\).
            \label{item:skyscraperbounds2}
    \end{enumerate}
    In particular if there are closed points \(p,q\in C_i\) such that \(\OO_p\in
    K\) and \(\OO_q\in K[1]\), then \(\CC K\) is non-zero.

    \begin{proof}
        If \(p\in C_i\) is a closed point
        such that \(\OO_p\in K\), then \(\OO_p\) lies in the torsion-free class
        \(F=K\cap H\). The exact triangle \[\OO_{C_i}(-1)\to
        \OO_p\to \OO_{C_i}(-2)[1]\to \OO_{C_i}(-1)[1]\] shows \(\OO_{C_i}(-1)\)
        is a sub-object (in \(H\)) of \(\OO_p\), thus we also have
        \(\OO_{C_i}(-1)\in F\). Further, inductively considering extensions
        \[\OO_{C_i}(n)\to \OO_{C_i}(n+1)\to \OO_p\to \OO_{C_i}(n)[1]\] shows
        \(\OO_{C_i}(n)\in F\) for every \(n\geq -1\), and thus \(F\) contains
        each torsion class \(\LL_i^{\otimes n}\otimes H[-1] \cap H\) by
        \cref{prop:lbtwistgeneration}. Statement \ref{item:skyscraperbounds1}
        follows, and the proof of \ref{item:skyscraperbounds2} is analogous.
    \end{proof}
\end{lemma}

Thus given an arbitrary non-algebraic t-structure \(K\), we show \(\CC K\) is
non-zero by comparing \(K\) against some appropriately chosen geometry. The
following lemma informs this choice.

\begin{lemma}
    \label{lem:anyKlivesonX}
    Given \(K\in \tilt(H)\), there is a sequence of flops \(\nu\) from \(X\)
    with corresponding birational model \(W=\nu X\) such that \(K\) lies in
    \(\Psi_\nu \zeroper\WZ[-1,0]\) and
    contains the objects \(\Psi_\nu \OO_{C_i}(-1)\) for each exceptional curve
    \(C_i\subset W\).

    \begin{proof}
        Consider the sub-poset of \(\tilt(H)\) given by
        \[
            [K, H]\;\cap\;
            \{\Psi_\nu \zeroper\pair{\nu X}{Z} \;|\; \nu\text{ a spherical \(\JJ\)-path}\}.
        \]
        This poset is non-empty (it contains \(H\)) and finite, so has a minimal
        element determined by some spherical \(\JJ\)-path \(\nu\). This
        determines our sequence of flops. 

        Continuing to write \(\Psi_\nu H\) for
        the heart \(\Psi_\nu \zeroper\pair{\nu X}{Z} = \Psi_\nu \rflmod({}_\nu
        \Lambda_\nu)\), we evidently have the chain of inequalities 
        \(H \;\geq\; \Psi_\nu H \;\geq\; K \;\geq\; H[-1] \;\geq\; \Psi_\nu H
        [-1]\), 
        so that \(K\) lies in \(\Psi_\nu H [-1,0]\). 

        Suppose \(i\in \Delta\setminus \nu\JJ\) is such that \(\Psi_\nu
        S_i=\Psi_\nu\OO_{C_i}(-1)\) lies outside \(K\), and hence outside the
        torsion-free class \(K\cap \Psi_\nu H\) in \(\Psi_\nu H\). Being a simple of
        \(\Psi_\nu H\), this object must then lie in the torsion class
        \(K[1]\cap \Psi_\nu H\). Noting \(\langle \Psi_\nu S_i\rangle \) is
        precisely the torsion class \( \Psi_\nu \Psi_i H [1] \cap \Psi_\nu H\)
        (\cref{lem:simpletorsiontheories}), we thus have
        \( \Psi_\nu H \;\geq\; \Psi_\nu \Psi_i H \;\geq\; K\).

        It follows that \(\Psi_\nu\Psi_i H\) lies in the interval
        \([K,H]\), and in particular is an algebraic heart in \(\tilt(H)\). Since its heart cone
        coincides with that of \( \Psi_{\nu_i\nu}\zeroper\pair{\nu_i \nu X}{Z}\)
        we conclude that the two hearts must be equal (\cref{thm:heartfandefn}).
        But this contradicts the minimality of \(\Psi_\nu H\), so we conclude
        that \(K\) must contain every simple \(\Psi_\nu S_i\) for \(i\in
        \Delta\setminus\nu\JJ\).
    \end{proof}
\end{lemma}

\begin{definition}
    \label{def:livesonW}
    Given a birational model \(W=\nu X\), we say a heart \(K\subset \Dm X\)
    \emph{lives on \(W\)} if it is intermediate with respect to \(\Psi_\nu
    \zeroper\WZ\) and contains the object \(\Psi_\nu\OO_{C_i}(-1)\) for every
    integral exceptional curve \(C_i\subset W\), i.e.\ if there are inclusions 
    \[ 
        \langle\Psi_\nu\OO_{C_i}(-1)\;|\; i\in \Delta\setminus \nu\JJ\rangle
        \;\subseteq\; K \;\subseteq\; \Psi_\nu\zeroper\WZ[-1,0]. 
    \] 
\end{definition}
Then \cref{lem:anyKlivesonX}
explains that every non-algebraic heart in \(\tilt(H)\) lives on \emph{some}
birational model, so up to the application of a flop functor we may assume \(K\)
lives on \(X\). This additional hypothesis gives us tremendous control on the torsion pair
associated to \(K\), in particular we show that every skyscraper on \(X\) is
either torsion or torsion-free.

\begin{lemma}
    \label{lem:everyskyscraperistorsion}
    If \(K\in \tilt(H)\) lives on \(X\) and \(p\in X(\mathbb{C})\) is a closed point in the
    exceptional locus, then the sheaf \(\OO_p\) lies in a single
    cohomological degree with respect to \(K\) (i.e.\ we have \(\OO_p\in K\) or
    \(\OO_p\in K[1]\).)

    \begin{proof}
        The intermediate heart \(K\) induces the torsion pair \(H=(K[1]\cap
        H)\ast (K\cap H)\), and thus an exact triangle  \(t\to \OO_p \to f \to
        t[1]\) with \(t\in K[1]\cap H\), \(f\in K\cap H\). In particular,
        considering \(\KK\)-theory classes 
        \[[t]=\sum_{i\in \ul\Delta\setminus\JJ}t_i\alpha_i, \qquad
        [\OO_p]= \sum_{i\in \ul\Delta\setminus\JJ}\delta_i \alpha_i,
        \qquad[f]=\sum_{i\in\ul\Delta\setminus\JJ}f_i\alpha_i,\] 
        we see that for
        each \(i\) the \(t_i,\delta_i,f_i\) are non-negative integers satisfying
        \(t_i+f_i=\delta_i\), and the coefficients \(\delta_i\) correspond to
        ranks of the summands of \(\VV\XZ\) as in
        \cref{lem:ktheoryofskyscraper}. In particular we have \(\delta_0=
        \text{rk}(\OO_X) = 1\), so we must have \((t_0,f_0)=(0,1)\) or
        \((t_0,f_0)=(1,0)\).

        If \(f_0=0\), then \(f\) is filtered by the simples
        \(S_i=\OO_{C_i}(-1)\) for \(i\in \Delta\setminus \JJ\), and is in
        particular a pure sheaf. But this implies \(\Hom(\OO_p,f)=0\), thus the map \(t\to
        \OO_p\) is an isomorphism and \(\OO_p\) lies in \(K[1]\).

        On the other hand if \(t_0=0\), then \(t\) lies in the extension-closure
        of the sheaves \(\OO_{C_i}(-1)\), and hence in \(K\) (since \(K\) lives
        on \(X\)). It follows that \(t\in K[1]\cap K\) must be the zero object,
        so that \(\OO_p=f\in K\).
    \end{proof}
\end{lemma}

We are equipped to show that an arbitrary t-structure in \(\tilt(H)\) is
numerical.

\begin{proof}[Proof of \cref{thm:everyheartisnumerical}]
    We may assume \(K\) is non-algebraic and, applying a flop functor if
    necessary, lives on \(X\). Thus by \cref{lem:everyskyscraperistorsion}, for
    every closed point \(p\) in the exceptional locus of \(X\), either \(\OO_p\) or
    \(\OO_p[-1]\) lies in \(K\).

    If there is some curve \(C_i\) with two points \(p,q\in C_i\) such that
    \(\OO_p\in K\) and \(\OO_q[-1]\in K\), then we obtain the conclusion using
    \cref{lem:skyscrapertrick}.

    Otherwise, since the exceptional locus of \(X\) is connected, either \(K\)
    or \(K[1]\) contains the subcategory \(\{\OO_p\;|\; p\in X(\mathbb{C})\}\).
    If \(\OO_p\in K\) for every closed point \(p\), then (by
    \cref{lem:skyscrapertrick}) we have \(K\geq \LL\otimes H[-1]\)
    for every nef bundle \(\LL\in \Pic X\) and thus \(K\geq \anticoh X\). Further, the
    torsion-free part \(K\cap H\) contains
    \(H_\ss(\Chamb_\JJ^0)=\langle\OO_p\;|\; p\in C\rangle\) so in fact \(K\geq
    \coh X\). By \cref{lem:boundedbycoh}, we conclude \(\CC K\neq 0\). 

    The case when \(\OO_p[-1]\in K\) for every \(p\) is handled likewise.
\end{proof}

\subsection{Classification of bricks}\label{subsec:brickclassification} We show
that every brick in \(H=\zeroper\XZ\) must arise from an algebraic or geometric
simple.

\begin{theorem}
    \label{thm:brickclassification}
    Let \(b\in H\) be a brick. Then either
    \begin{enumerate}
        \item[(1)] there is a spherical \(\JJ\)-path \(\nu\) and a closed point
            \(p\in \nu X(\mathbb{C})\) such that \(b=\Psi_\nu \OO_p\) , or
        \item[(2)] there is a \(\JJ\)-path \(\nu\) and an index
            \(i\in\ul\Delta\setminus\nu\JJ\) such that \(b\in \left\{\Psi_\nu
            S_i[1], \Phi_\nu S_i[-1]\right\}\).
            \end{enumerate}
    In particular, the \(\KK\)-theory class \([b]\in
    \mathfrak{h}(\ul\Delta\setminus\JJ)\) is a primitive restricted root.
\end{theorem}

We prove the theorem over the course of this subsection. Thus fix a
brick \(b\in H\), and consider the torsion pair \(H=T\ast
F\) where \(T\) is the minimal torsion class containing \(b\) and accordingly
\(F=\{h\in H\;|\; \Hom(b,h)=0\}\). 

By \cref{lem:anyKlivesonX}, the corresponding
tilted heart \(K=F\ast T[-1]\) lives on some birational model \(W=\upsilon X\).
We show that \(b\) then lies in a single cohomological degree with respect to
\(\Psi_\upsilon \coh W\), i.e.\ \(\Psi_\upsilon^{-1}b\) or
\(\Psi_\upsilon^{-1}b[-1]\) is a sheaf on \(W\). To this end, we establish the
following preliminary lemma.

\begin{lemma}
    Let \(\mathcal{F}\) be a coherent sheaf on a variety, with
    \(\dim\Supp(\mathcal{F})\leq 1\). If \(p\) is a closed point such that 
    \(\Hom(\mathcal{F},\OO_p)=0\) or
    \(\Hom(\OO_p,\mathcal{F})=\Ext^1(\OO_p,\mathcal{F})=0\), then \(p\) lies
    outside \(\Supp(\mathcal{F})\).

    \begin{proof}
        The conclusion is immediate from the first hypothesis. On the other hand
        if \(\Hom(\OO_p,\mathcal{F})=\Ext^1(\OO_p,\mathcal{F})=0\), then in
        particular we have \(\curlyHom(\OO_p,\mathcal{F})=0\) where we note that
        the homomorphism sheaf has zero dimensional support and is therefore
        determined by its module of global sections. Likewise we can deduce
        (e.g.\ from the Grothendieck spectral sequence associated to
        \(\Gamma\circ\mathbf{R}\curlyHom(\OO_p,-)\)) that
        \(\Gamma(\curlyExt^1(\OO_p,\mathcal{F}))=\Ext^1(\OO_p,\mathcal{F})=0\)
        and hence \(\curlyExt^1(\OO_p,\mathcal{F})=0\).

        Taking stalks at \(p\), we thus have
        \(\Hom(\OO_p,\mathcal{F}_p)=\Ext^1(\OO_p,\mathcal{F}_p)=0\) where
        \(\mathcal{F}_p\) is a finitely generated module over the local ring at
        \(p\), and \(\OO_p\) is the unique simple module. It follows that
        \(\mathcal{F}_p\) is either zero (i.e.\ \(p\notin\Supp(\mathcal{F})\))
        or \(\text{depth}(\mathcal{F}_p)\geq 2\)---but the latter situation is
        absurd since the depth of a module is bounded above by the dimension of
        its support.
    \end{proof}
\end{lemma}

We also require the following lemma, which follows from a more general
observation due to Sota Asai regarding torsion classes generated by bricks.

\begin{lemma}

    Let \(b,T\) be as above, and suppose \(b'\) is another brick that lies in \(T\). Then
    either \(b=b'\), or \(\Hom(b',b)=0\).

    \begin{proof}
        By~\cite[lemma 1.7]{asaiSemibricks2020}, if \(\Hom(b',b)=0\) then there
        is a surjection \(b'\onto b\) in \(H\), and in particular, \(b\) lies in
        the minimal torsion class \(T'\) containing \(b'\). Thus we must have
        \(T=T'\). But the set of torsion classes in \(H\) generated by bricks (i.e.\
        \emph{completely join irreducible} torsion classes) are naturally in
        bijection with the set of bricks in \(H\)~\cite[theorem
        1.5]{barnardMinimalInclusionsTorsion2019} so in fact \(b=b'\) as
        required.
    \end{proof}
\end{lemma}

By \cref{lem:anyKlivesonX}, \(K\) lives on some birational
model \(W=\upsilon X\) and in particular, for every closed point \(p\in
W(\mathbb{C})\) we either have \(\Psi_\upsilon \OO_p\in T\) or \(\Psi_\upsilon
\OO_p\in F\) (\cref{lem:everyskyscraperistorsion}). By the above lemma, for each
\(p\in W(\mathbb{C})\) we thus have 
\[
    \Hom(\Psi_\upsilon\OO_p,b)=0 \qquad\text{or}\qquad
    b=\Psi_\upsilon\OO_p \qquad\text{or}\qquad
    \Hom(b,\Psi_\upsilon\OO_p)=0
\]
where the first two cases occur when \(\Psi_\upsilon\OO_p\in T\), while the
last is equivalent to \(\Psi_\upsilon\OO_p\in F\).

\begin{proposition}
    Let \(b\in H\) be a brick, \(H=T\ast F\) the torsion pair whose torsion
    class \(T\) is generated by \(b\), and \(W=\upsilon X\) a birational model
    such that \(K=F\ast T[-1]\) lives on \(W\). Then one of the
    following (mutually exclusive) possibilities is true. 
    \begin{enumerate}[(1)]
        \item Either \(b=\Psi_\upsilon\OO_q\) for some closed point \(q\in
            W(\mathbb{C})\), 
        \item or \(b=\Psi_\upsilon\mathcal{F}[1]\) for a pure sheaf
            \(\mathcal{F}\in \coh(W)\) with one-dimensional support,
        \item or \(b=\Psi_\upsilon \mathcal{G}\) for a pure sheaf
            \(\mathcal{G}\in \coh(W)\) with one-dimensional support.
    \end{enumerate}
    \begin{proof}
        If \(b=\Psi_\upsilon\OO_p\) for some \(p\) then there is nothing to prove,
        so we assume that doesn't happen. Thus from the above discussion, for every \(p\) we have
        \(\Hom(\Psi_\upsilon\OO_p,b)=0\) or \(\Hom(b,\Psi_\upsilon\OO_p)=0\).

        Considering the torsion pair \(H=(\Psi_\upsilon\coh W[1]\cap H)\ast
        (\coh W\cap H)\) shows \(b\) lies in an exact triangle 
        \begin{equation} 
            \label{eqn:bsheaves} 
            \Psi_\upsilon \mathcal{F}[1]\to b\to \Psi_\upsilon \mathcal{G}\to
            \Psi_\upsilon \mathcal{F}[2] 
        \end{equation} 
        where \(\mathcal{F}\in \coh(W)\cap \Psi_\upsilon^{-1}H[-1]\) and
        \(\mathcal{G}\in \coh(W)\cap \Psi_\upsilon^{-1}H\) are the cohomology
        objects of \(\Psi_\upsilon^{-1}b\) seen as a complex of sheaves on
        \(W\). 

        For any \(p\in W(\mathbb{C})\), note that \(\OO_p\) lies in the heart
        \(\Psi_\upsilon^{-1}H\) and thus \(\Hom(\OO_p,\mathcal{F})= 0\). 
        On the other hand, considering long exact sequences shows that
        \(\Hom(b,\Psi_\upsilon\OO_p)=\Hom(\mathcal{G},\OO_p)\) and
        \(\Ext^1(\OO_p,\mathcal{F})\subset \Hom(\Psi_\upsilon \OO_p, b)\). 

        It
        follows that if \(\Hom(b,\Psi_\upsilon\OO_p)=0\) then \(p\notin\Supp(\mathcal{G})\),
        while if \(\Hom(\Psi_\upsilon \OO_p,b)=0\) then \(p\notin \Supp(\mathcal{F})\). Thus
        \(\Supp(\mathcal{F})\cap\Supp(\mathcal{G})=\varnothing\), so that the
        sequence~\eqref{eqn:bsheaves} splits. But \(b\) is a brick (in
        particular indecomposable), so one of \(\mathcal{F}\) or \(\mathcal{G}\)
        must vanish.

        In the case \(b=\Psi_\upsilon \mathcal{F}[1]\), the fact that
        \(\mathcal{F}\) is pure of one-dimensional support follows from the
        observation that \(\Hom(\OO_p,\mathcal{F})=0\) for all \(p\).

        On the other hand if \(b=\Psi_\upsilon\mathcal{G}\) and \(\mathcal{G}\) has zero-dimensional support, then the condition \(\End(\mathcal{G})=\mathbb{C}\) ensures
        \(\mathcal{G}=\OO_q\) for some \(q\in W(\mathbb{C})\) and we are reduced
        to the first case. If \(\mathcal{G}\) has one-dimensional support and \(p\) is a point in
        \(\Supp(\mathcal{G})\), then we must have \(\Hom(\OO_p,\mathcal{G})=0\)
        so that \(\mathcal{G}\) is pure.
    \end{proof}
\end{proposition}

The proof of \cref{thm:brickclassification} is then reduced to an analysis of
each of the cases arising in the proposition above. If \(b=\Psi_\upsilon\OO_q\)
for a closed point then there is nothing to prove. The other two cases are
addressed below.

\begin{lemma}
    \label{lem:smallbrickclassification}
    If \(b=\Psi_\upsilon\mathcal{F}[1]\) for some sheaf
    \(\mathcal{F}\in \coh(W)\cap \Psi_\upsilon^{-1}H[-1]\), then the heart \(K\)
    is algebraic of the form \(\Psi_\nu H\) for some \(\JJ\)-path \(\nu\).
    Consequently, \(b=\Psi_\nu S_i[1]\) for some \(i\in
    \ul\Delta\setminus\nu\JJ\).

    \begin{proof}
        By (the dual statement to) \cref{prop:sbrickqAcorrespondence} the
        torsion-free class \(\Psi_\upsilon\coh(W)\cap H\) is generated by the
        semibrick \(\{\Psi_\upsilon\OO_p\;|\; p\in W(\mathbb{C})\}\). The
        constraints on \(b\) guarantee that every \(\Psi_\upsilon\OO_p\) lies in
        \(F\), and thus \(\Psi_\upsilon\coh(W)\cap H\subset F\). Equivalently,
        we have the containment of torsion classes
        \[
            T\;\subseteq\; \Psi_\upsilon\coh(W)[1]\cap H \;=\;
            \sup\left\{\LL\,\dual\otimes H[1]\cap H \;\middle\vert\; {\LL\in
            \Pic W\text{ nef}} \right\}.
        \]
        The supremum on the right hand side is an infinite union (since
        any two elements in the set are dominated by a third as
        in~\eqref{eqn:lbtorsionclassdomination}). It follows that \(b\), and
        hence also \(T\), is contained in the torsion class \(\LL\,\dual\otimes
        H[1] \cap H\) for some \(\LL\in \Pic W\) nef. In other words, \(K\) is
        contained in the interval \([\LL\,\dual\otimes H, H]\subset
        \tilt^+(H)\). The result follows.
    \end{proof}
\end{lemma}

\begin{lemma}
    If \(b=\Psi_\upsilon\mathcal{G}\) for some purely one-dimensional sheaf
    \(\mathcal{G}\in \coh(W)\cap \Psi_\upsilon^{-1}H\), then \(b\) is a
    simple object in an algebraic heart \(\Phi_\nu H[-1]\) determined by some
    \(\JJ\)-path \(\nu\). Consequently, \(b=\Phi_\nu S_i[-1]\) for some \(i\in
    \ul\Delta\setminus\nu\JJ\).

    \begin{proof}
        Consider the torsion pair \(H=T'\ast F'\) where \(F'\) is the smallest
        torsion-free class  in \(H\) that contains \(b\). By (the dual statement
        to) \cref{prop:sbrickqAcorrespondence}, \(b\) is a simple in the tilted
        heart \(K'=F'\ast T'[-1]\). 

        Since \(\mathcal{G}\) is pure of one
        dimensional support, we have
        \(\Hom(\Psi_\nu\OO_p,b)=\Hom(\OO_p,\mathcal{G})=0\) for all \(p\in
        W(\mathbb{C})\), and hence each \(\Psi_\nu\OO_p\) lies in \(T'\). It
        follows that so does the torsion class \(\Psi_\nu\anticoh(W)[1]\cap H\),
        which is generated by the semibrick \(\{\Psi_\upsilon\OO_p\;|\; p\in
        W(\mathbb{C})\}\). We can then conclude in a similar fashion as the
        proof of \cref{lem:smallbrickclassification}.
    \end{proof}
\end{lemma}

\fakesection{References}
\printbibliography
\end{document}